\documentclass{amsart}

\usepackage{latexsym,amsthm,mathrsfs,amsmath,amscd,enumerate,enumitem, amscd,color,dsfont, cancel, mathtools, textcomp, float}

\usepackage[utf8x]{inputenc}

\usepackage[cspex,bbgreekl]{mathbbol}
\usepackage{amsfonts}
\usepackage{amssymb}             

\DeclareSymbolFontAlphabet{\mathbbl}{bbold}
\DeclareSymbolFontAlphabet{\mathbb}{AMSb}%

 \newtheorem{thm}{Theorem}[section]

\theoremstyle{definition}

 \theoremstyle{remark}

\usepackage{hyperref}
\usepackage{tikz}
\usetikzlibrary{plotmarks}

\newcommand{\supp}{\mathop{\mathrm{supp}}}

\newcommand{\loc}{\textup{loc}}
\newcommand{\glob}{\textup{glob}}
\makeatletter
\DeclareRobustCommand\widecheck[1]{{\mathpalette\@widecheck{#1}}}
\def\@widecheck#1#2{%
    \setbox\z@\hbox{\m@th$#1#2$}%
    \setbox\tw@\hbox{\m@th$#1%
       \widehat{%
          \vrule\@width\z@\@height\ht\z@
          \vrule\@height\z@\@width\wd\z@}$}%
    \dp\tw@-\ht\z@
    \@tempdima\ht\z@ \advance\@tempdima2\ht\tw@ \divide\@tempdima\thr@@
    \setbox\tw@\hbox{%
       \raise\@tempdima\hbox{\scalebox{1}[-1]{\lower\@tempdima\box
\tw@}}}%
    {\ooalign{\box\tw@ \cr \box\z@}}}
\makeatother

\usepackage{caption}
\usepackage{subcaption}

\usepackage[left=2.2cm,top=2.5cm,right=2.2cm,bottom=2.5cm]{geometry} 

\numberwithin{equation}{section}
\allowdisplaybreaks

\begin{document}

\title[]{Littlewood-paley functions associated with general Ornstein-Uhlenbeck semigroups}

\author[V. Almeida]{V\'{\i}ctor Almeida}

\author[J.J. Betancor]{Jorge J. Betancor}

\author[J.C. Fari\~na]{Juan C. Fari\~na}

\author[P. Quijano]{Pablo Quijano}

\author[L. Rodr\'{\i}guez-Mesa]{Lourdes Rodr\'{\i}guez-Mesa}

\address{V\'{\i}ctor Almeida, Jorge J. Betancor, Juan C. Fari\~na and Lourdes Rodr\'{\i}guez-Mesa\newline
	Departamento de An\'alisis Matem\'atico, Universidad de La Laguna,\newline
	Campus de Anchieta, Avda. Astrof\'isico S\'anchez, s/n,\newline
	38721 La Laguna (Sta. Cruz de Tenerife), Spain}
\email{valmeida@ull.edu.es, jbetanco@ull.es, jcfarina@ull.edu.es,
lrguez@ull.edu.es
}

\address{Pablo Quijano\newline
	Instituto de Matemática Aplicada del Litoral, Santa Fe-Argentina}
\email{pabloquijanoar@gmail.com}

\thanks{The authors are partially supported by grant PID2019-106093GB-I00 from the Spanish Government}



\date{\today}

\begin{abstract}
In this paper we establish $L^p(\mathbb{R}^d,\gamma_\infty)$-boundedness properties for square functions involving time and spatial derivatives of Ornstein-Uhlenbeck semigroups. Here $\gamma_\infty$ denotes the invariant measure. In order to prove the strong type results for $1<p<\infty$ we use $R$-boundedness. The weak type (1,1) property is established by studying separately global and local operators defined for the square Littlewood-Paley functions. By the way we prove $L^p(\mathbb{R}^d,\gamma _\infty )$-boundedness properties for maximal and variation operators for Ornstein-Uhlenbeck semigroups.

\end{abstract}

\maketitle

\section{Introduction}\label{S1}
In this paper we establish $L^p(\mathbb R^d,\gamma_\infty)$-boundedness properties of Littlewood-Paley functions associated with Ornstein-Uhlenbeck semigroups.

Suppose that $Q$ is a real, symmetric and positive definite $d \times d$ matrix and that $B$ is a real $d \times d$-matrix having all its eigenvalues with negative real parts. $Q$ and $B$ are usually named the covariance and the drift matrix, respectively. For every $t \in (0,\infty]$ we define the symmetric and positive measure matrix $Q_t$ by
$$
Q_t= \int_0^t e^{sB} Q e^{sB^*} ds
$$
and the normalized measure $\gamma_t$ in $\mathbb R^d$ by
$$
d\gamma_t(x)=(2\pi)^{-\frac{d}{2}}({\rm det}\;Q_t)^{-\frac{1}{2}} e^{-\frac{1}{2}\langle Q^{-1}_tx,x\rangle}dx.
$$
By $C_b(\mathbb R^d)$ we denote the space of bounded continuous functions in $\mathbb R^d$. We consider the semigroup of operators $\{\mathcal H_t\}_{t>0}$ where, for every $f\in C_b(\mathbb R^d)$,
$$
\mathcal H_t(f)(x) = \int_{\mathbb R^d} f(e^{tB}x-y)\,d\gamma_t(y),\quad x\in \mathbb R^d.
$$
$\{\mathcal H_t\}_{t>0}$ is called the Ornstein-Uhlenbeck semigroup defined by $Q$ and $B$. $\gamma_\infty$ is the unique invariant measure with respect to $\{\mathcal H_t\}_{t>0}$.

After some manipulations we can write, for every $f\in C_b(\mathbb R^d)$,
$$
\mathcal H_t(f)(x)= \int_{\mathbb R^d}h_t(x,y) f(y) d\gamma_\infty(y),\quad x \in \mathbb R^d \mbox{ and } t
>0,
$$
where
\begin{equation}\label{(1.0)}
h_t(x,y) = \left( \frac{{\rm det}\,Q_\infty}{{\rm det}\,{Q_t}}\right)^\frac{1}{2} e^{R(x)}\exp\left[-\frac{1}{2}\langle (Q_t^{-1}-Q_\infty^{-1})(y-D_tx),y-D_tx\rangle\right],\quad x,y\in \mathbb R^d \mbox{ and }t>0,
\end{equation}
being $D_t=Q_\infty e^{-tB^*}Q_\infty^{-1}$, $t>0$, and $R(x)=\frac{1}{2} \langle Q_\infty^{-1}x, x\rangle$, $x\in \mathbb R^d$. When $Q=I$ and $B=-I$ the semigroup $\{\mathcal H_t\}_{t>0}$ reduces to the symmetric Ornstein-Uhlenbeck semigroup.

Let $1 \leq p < \infty$. $\{\mathcal H_t\}_{t>0}$ extends to a positivity preserving semigroup of contractions in $L^p(\mathbb R^d,\gamma_\infty)$. By denoting  by $-\mathfrak{L}_p$ the infinitesimal generator of  $\{\mathcal H_t\}_{t>0}$ in $L^p(\mathbb R^d, \gamma_\infty)$, the space $C_c^\infty(\mathbb R^d)$ of smooth and compactly supported functions in $\mathbb R^d$ is a core for $\mathfrak{L}_p$ and, for every $f \in C_c^\infty(\mathbb R^d)$,
$$
\mathfrak{L}_p(f)(x) = -\frac{1}{2}{\rm div}(Q\nabla f)(x) - \langle \nabla f(x),Bx\rangle, \quad x\in \mathbb R^d,
$$
(\cite[Chapter 5]{LB}). If $1< p <\infty$ $\{\mathcal H_t\}_{t>0}$ extends to an analytic contraction semigroup in $L^p(\mathbb R^d,\gamma_\infty)$ in the sector $S_{\theta_p}$ of angle $\theta_p=\frac{\pi}{2} -\varphi_p$ where
$$
\varphi_p=\arctan\frac{\sqrt{(p-2)^2 + p^2(\tan \varphi_{\mbox{\tiny A}})^2}}{2\sqrt{p-1}}
$$
and $\varphi=\arctan \|Q^{-\frac{1}{2}}(B-B^*)Q^{-\frac{1}{2}}\|$ (\cite[Theorem 2 and Remark 6]{CFMP}).

The operator $\mathcal L_p$ is sectorial but it is not one to one being the kernel of $\mathcal L_p$ the subspace $N(\mathcal L_p)$ consisting of all the constant functions (\cite[Theorem 8.1.17]{LB}). $\mathcal L_p$ is a sectorial one to one operator with dense range on $L^p_0(\mathbb R^d,\gamma_\infty) = \{f\in L^p(\mathbb R^d,\gamma_\infty): \int_{\mathbb R^d} f d\gamma_\infty=0\}$. In \cite[Theorem 10]{CD} it was established that if $1<p<\infty$, $\mathcal L_p$ has bounded holomorphic functional calculus whose sharp angle is $\theta_p$.

The study of harmonic analysis in the Ornstein-Uhlenbeck setting was began by Muckenhoupt \cite{Mu1} who considered the one dimensional case. In the higher dimension the situation is quite different and new arguments and ideas are needed. E. M. Stein, in his celebrated monography \cite{Ste}, developed a general theory for harmonic analysis associated with symmetric diffusion semigroup. The symmetric Ornstein-Uhlenbeck semigroup is a special case of Stein symmetric diffusion semigroup. Sj\"ogren (\cite{Sj}) extended the Muckenhoupt result proving that the maximal operator defined by $\{\mathcal H_t\}_{t>0}$ is bounded from $L^1(\mathbb R^d,\gamma_\infty)$ into $L^{1,\infty}(\mathbb R^d,\gamma_\infty)$ for every $d \geq 1$. Twenty years later, higher order Riesz transforms in the Ornstein-Uhlenbeck setting were studied in \cite{G} and \cite{GST} where it is proved that they are bounded operators from $L^p(\mathbb R^d,\gamma_\infty)$ into itself, for every $1<p<\infty$. Riesz transforms are not bounded from $L^1(\mathbb R^d,\gamma_\infty)$ into itself. Harmonic analysis operators associated with the symmetric Ornstein-Uhlenbeck have been studied in the last two decades. Some of these operators are the following ones: maximal operators (\cite{GCMMST1} and \cite{MPS2}), Littlewood-Paley functions (\cite{HTV}, \cite{P2} and \cite{PS}), spectral multipliers (\cite{GCMST2}), singular integrals (\cite{AFS}, \cite{GCMST1} and \cite{P1}), operators with $H^\infty$-functional calculus (\cite{GCMMST2}) and variation operators (\cite{HMMT}). All those operators are bounded from $L^p(\mathbb R^d,\gamma_\infty)$ into itself, for every $1<p<\infty$, and from $L^1(\mathbb R^d,\gamma_\infty)$ into $L^{1,\infty}(\mathbb R^d,\gamma_\infty)$.

In contrast with the symmetric case, harmonic analysis operators in the general nonsymmetric setting had not been very studied. Mauceri and Noselli considered the maximal operator (\cite{MN1}) and Riesz transforms of the first order (\cite{MN2}) where $Q=I$ and $B=-\lambda(I+R)$ being $\lambda >0$ and $R$ generates a periodic group $\{e^{Rt}\}_{t>0}$ of rotations. Recently, Casarino, Ciatti and Sj\"ogren have studied maximal operators (\cite{CCS2} and \cite{CCS1}),  Riesz transforms (\cite{CCS3}) and spectral multipliers (\cite{CCS4})  associated with general nonsymmetric Ornstein-Uhlenbeck operators.

One of our objective in this paper is to established the $L^p$-boundedness properties of some Littlewood-Paley functions, also called square functions, involving time and spatial derivatives of subordinated Ornstein-Uhlenbek semigroups in the nonsymmetric setting.

Let $\nu>0$. We consider the Poisson-like integral defined by
\begin{equation}\label{(1.1)}
P^\nu_t(f)(x)= \frac{t^{2\nu}}{4^\nu\Gamma(\nu)}\int_0^\infty e^{-\frac{t^2}{4u}}\mathcal H_u(f)(x)\frac{du}{u^{\nu+1}},\quad  x \in \mathbb R^d\mbox{ and } t>0,
\end{equation}
that is subordinated to the Ornstein-Uhlenbeck $\{\mathcal H_t\}_{t>0}$. Note that when $\nu=\frac{1}{2}$ we recover the Poisson semigroup generated by $-\sqrt{\mathcal L}$. According to \cite[(1.9)]{ST}, (\ref{(1.1)}) defines a solution of the initial value problem
$$
\left.\begin{array}{rcl}
\displaystyle \partial_t^2u + \frac{1-2\nu}{t} \partial_t u & = & \mathcal L u\\
u(0,x) & = & f(x)
\end{array}\right\}
$$
for suitable $f$, and this partial differential equation is connected with the fractional power $\mathcal L^\nu$ of $\mathcal L$.

Let $k \in \mathbb N$ and $\alpha=(\alpha_1,\ldots,\alpha_d)\in \mathbb N^d$ such that $k +\displaystyle \sum_{i=1}^d \alpha_i >0$. We define the Littlewood-Paley function $g_{k,\alpha}^\nu$ as follows
$$
g_{k,\alpha}^\nu(f)(x) = \left(\int_0^\infty \left|t^{k+\widehat\alpha}\partial_t^k\partial_x^\alpha P_t^\nu(f)(x)\right|^2 \frac{dt}{t}\right)^{\frac{1}{2}}, \quad x \in \mathbb R^d,
$$
where $\displaystyle\widehat \alpha=\sum_{i=1}^d\alpha_i$ and $\displaystyle\partial^\alpha_x=\partial_{x_1}^{\alpha_1}\ldots\partial _{x_d}^{\alpha_d}$.

We now state one of our main result.

\begin{thm}\label{Th1.1}
Let $\nu>0$, $k\in \mathbb N$ and $\alpha=(\alpha_1,\ldots,\alpha_d)$ such that $k+\widehat\alpha>0$. Then, $g^\nu_{k,\alpha}$ is bounded from $L^p(\mathbb R^d,\gamma_\infty)$ into itself, for every $1<p<\infty$. If in addition $\widehat\alpha=1$ and $d>1$ or $\widehat{\alpha}=0,2$ and $d>2$, $g_{k,\alpha}^\nu$ is bounded from $L^1(\mathbb R^d,\gamma_\infty)$ into $L^{1,\infty}(\mathbb R^d,\gamma_\infty)$. 
\end{thm}

Some comments about the proof of this theorem are in order. To see that $g^\nu_{k,\alpha}$ is bounded from $L^p(\mathbb R^d,\gamma_\infty)$ into itself, for every $1<p<\infty$, when $\widehat \alpha > 0$, we consider the square function $G_\alpha$ associated with $\{\mathcal H_t\}_{t>0}$ defined by
$$
G_\alpha(f)(x)= \left(\int_0^\infty \left|t^{\frac{\widehat \alpha}{2}}\partial_x^\alpha\mathcal H_t(f)(x)\right|^2\frac{dt}{t}\right)^{\frac{1}{2}}, \quad x \in \mathbb R^d.
$$
Let $1<p<\infty$ and $\alpha \in \mathbb N^d \setminus \{0\}$. The $\alpha$-order Riesz transform in this context is given by
$$
\mathcal R_\alpha (f)= \partial^\alpha_x \mathcal L^{-\frac{\widehat{\alpha}}{2}}\Pi_0(f),\quad f\in L^p_0(\mathbb R^d,\gamma_\infty),
$$
where $\Pi_0$ denotes the projection from $L^p(\mathbb R^d,\gamma_\infty)$ to $L^p_0(\mathbb R^d,\gamma_\infty)$ and the $-\frac{\widehat{\alpha}}{2}$ power of $\mathcal L$ is defined by
$$
\mathcal L^{-\frac{\widehat{\alpha}}{2}}g = \frac{1}{\Gamma(\frac{\widehat{\alpha}}{2})}\int_0^\infty t^{\frac{\widehat{\alpha}}{2}-1}\mathcal H_t(g) dt,\quad g\in L^p_0(\mathbb R^d,\gamma_\infty).
$$
Here $\partial_x^\alpha$ is understood as a distributional derivative. According to \cite[Proposition 2.3]{MN2}, $\mathcal R_\alpha$ is bounded from $L^p(\mathbb R^d,\gamma_\infty)$ into itself. We prove that the set $\{t^{\frac{\widehat{\alpha}}{2}}\partial_x^\alpha\mathcal H_t\}_{t>0}$ is $R$-bounded in $L^p(\mathbb R^d, \gamma_\infty)$ and then that the square function $G_\alpha$ is bounded from $L^p( \mathbb R^d,\gamma_\infty)$ into itself. We recall that a family $\mathcal A$ of bounded operators of $L^p(\mathbb R^d, \gamma_\infty)$ into itself is said to be $R$-bounded in $L^p(\mathbb R^d,\gamma_\infty)$ when there exists $C>0$ such that if $T_k \in \mathcal A$ and $f_k \in L^p(\mathbb R^d,\gamma_\infty)$, $k=1,\ldots,n$, with $n \in \mathbb N \setminus\{0\}$, then
$$
\mathbb E \left\| \sum_{k=1}^n r_kT_kf_k\right\|_{L^p(\mathbb{R}^d,\gamma_\infty)} \leq C \mathbb E\left\|\sum_{k=1}^n r_kf_k\right\|_{L^p(\mathbb R^d,\gamma_\infty)},
$$
where $\{r_k\}_{k=1}^\infty$ is a sequence of independent Rademacher variables and $\mathbb E$ denotes as usual the expectation. The main properties of the $R$-bounded sets of operators can be found in \cite{HVNVW}.

Finally we can conclude that $g^\nu_{k,\alpha}$ is bounded from $L^p(\mathbb R^d,\gamma_\infty)$ into itself because $G_\alpha$ has this property.

We now fix our attention in the case $\alpha=0$ and $k \in \mathbb N$, $k \geq 1$. Let $1<p<\infty$. It was mentioned that the Ornstein-Uhlenbeck operator $\mathcal L_p$ has bounded holomorphic functional calculus with sharp angle $\theta_p$. Let $\theta\in (\theta_p,\frac{\pi}{2})$. We denote $S_\theta$ the sector of angle $\theta$ and $\Psi(S_\theta)$ the function space that consists of all those analytic functions $m$ in $S_\theta$ such that for a certain $s>0$ the function $M$ defined by
$$
M(z)=\frac{(1+z)^{2s}}{z^s}m(z),\quad z\in S_\theta,
$$
is bounded in $S_\theta$. According to \cite[Corollary 6.7]{CDMY}, for every $m\in \Psi(S_\theta)$ there exists $C>0$ such that, for every $f\in L^p(\mathbb R^d,\gamma_\infty)$,
\begin{equation}\label{(1.2)}
\left\|\left(\int_0^\infty|m(t\mathcal L_p)f(\cdot)|^2\frac{dt}{t}\right)^{\frac{1}{2}}\right\|_{L^p(\mathbb R^d,\gamma_\infty)} \leq C\|f\|_{L^p(\mathbb R^d,\gamma_\infty)}.
\end{equation}
By defining $m(z)=ze^{-z}$, $z\in S_\theta$, \eqref{(1.2)} implies that the Littlewood-Paley function $G^1$ associated with $\{\mathcal H_t\}_{t>0}$ and defined by
$$
G^1(f)(x)=\left(\int_0^\infty |t\partial_t\mathcal{H}_t(f)(x)|^2 \frac{dt}{t}\right)^{\frac{1}{2}},\quad x \in \mathbb R^d,
$$
is bounded from $L^p(\mathbb R^d,\gamma_\infty)$ into itself, and then $g_{k,0}^\nu$ is bounded from $L^p(\mathbb R^d,\gamma_\infty)$ into itself (see Section~\ref{section2}).

We remark that in the symmetric case, that is, when $Q=I$ and $B=-I$, the $L^p(\mathbb R^d,\gamma_\infty)$-boundedness of $G^1$ can be deduced from \cite[Theorem 10, p. 111]{Ste}. Furthermore, in the proof of \cite[Theorem 4.2]{MN1} it was proved that $G^1$ is bounded from $L^2(\mathbb R^d,\gamma_\infty)$ into itself.

In order to see that $g_{k,\alpha}^\nu$ is bounded from $L^1(\mathbb R^d,\gamma_\infty)$ into $L^{1,\infty}(\mathbb R^d,\gamma_\infty)$ we can not use Littlewood-Paley functions defined by $\{\mathcal H_t\}_{t>0}$ as above. As far as we know, it has not been proved that $G^1$ is bounded from $L^1(\mathbb R^d,\gamma_\infty)$ into $L^{1,\infty}(\mathbb R^d,\gamma_\infty)$ even in the symmetric case. Furthermore, $g_{0,\alpha}^{\frac{1}{2}}$ is not bounded from $L^1(\mathbb R^d,\gamma_\infty)$ into $L^{1,\infty}(\mathbb R^d,\gamma_\infty)$ when $\widehat{\alpha}>2$ in the symmetric setting. 

We prove that $g^\nu_{k,\alpha}$ is of weak type $(1,1)$ with respect to $\gamma_\infty$ when $0\leq \widehat\alpha \leq 2$ by considering two operators called the local and the global parts of $g_{k,\alpha}^\nu$. This local-global method appears in the first time in the seminal Muckenhoupt's paper \cite{Mu1} and it is usual to study $L^p$-boundedness properties of harmonic analysis operators in the Ornstein-Uhlenbeck context. The local part of the original operator $g_{k,\alpha}^\nu$ can be seen as a vector valued singular integral while the global part can be controlled by a positive operator. Our proof does not use Calder\'on-Zygmund theory to establish the weak type $(1,1)$ for the local part. We exploit the fact that in the local region, close to the diagonal, the Ornstein-Uhlenbeck semigroup is in some sense a nice perturbation of the classical heat semigroup. In order to establish our result we need some estimates involving the integral kernel of the Ornstein-Uhlenbeck semigroup due to Casarino, Ciatti and Sj\"ogren (\cite{CCS2}, \cite{CCS1} and \cite{CCS3}). Our comparative procedure to deal with the local operators can be applied to study Littlewood-Paley functions associated to generalized
Ornstein-Uhlenbeck operators, where the norm $L^2((0,\infty),dt/t)$ is
replaced by $L^q((0,\infty),dt/t)$, $1<q<\infty$, in Banach valued
setting and also to obtain new characterizations for the uniformly convex
and smooth Banach spaces (see \cite{BFRST}, \cite{HTV}, \cite{MTX3} and \cite{Xu}).

We consider the maximal operator $P_{*,k,\alpha}^\nu$ defined by
$$
 P_{*,k,\alpha}^\nu(f)=\sup_{t>0}\left|t^{k+\widehat\alpha} \partial^k_t\partial^\alpha_xP_t^\nu(f)\right|
$$
where $k\in \mathbb N$ and $\alpha \in \mathbb N^d$. As a  consequence of Theorem \ref{Th1.1} we can establish the following result that extends \cite[Theorem 1.1]{CCS3}.

\begin{thm}\label{Th1.2}
Let $\nu > 0$, $k\in \mathbb N$ and $\alpha \in \mathbb N^d$. Then, $P_{*,k,\alpha}^\nu$ is bounded from $L^p(\mathbb R^d,\gamma_\infty)$ into itself, for every $1 <p<\infty$. Furthermore, $P_{*,k,\alpha}^\nu$ is bounded from $L^1(\mathbb R^d,\gamma_\infty)$ into $L^{1,\infty}(\mathbb R^d,\gamma_\infty)$ provided that $\widehat{\alpha} \leq 1$ and $d>1$ or when $\widehat{\alpha}=2$ and $d>2$.
\end{thm}

Let $\rho >2$. Suppose that $g$ is a complex function defined in $(0,\infty)$. The $\rho$-variation $V_\rho(g)$ is defined
$$
V_\rho(g) = \sup_{\substack{0<t_k<t_{k-1}<\ldots <t_1\\ k\in \mathbb N}} \left(\sum_{j=1}^{k-1} |g(t_{j+1})-g(t_{j})|^\rho\right)^{\frac{1}{\rho}}.
$$
For the exponent $\rho = 2$ the oscillation associated to a particular sequence $\{t_j\}_{j\in \mathbb N}$ is considered instead of the variation. Suppose that $\{t_j\}_{j\in \mathbb N}$ is  a decreasing sequence in $(0,\infty)$ such that $t_j \rightarrow 0$, as $j \rightarrow \infty$. We define the oscillation $O(g,\{t_j\}_{j \in \mathbb N})$ associated with $\{t_j\}_{j \in \mathbb N}$ by
$$
O(g,\{t_j\}_{j \in \mathbb N})=\left( \sum_{j \in \mathbb N}\sup_{t_{j+1}\leq \varepsilon_{i+1} <\varepsilon_i\leq t_j}|g(\varepsilon_{i+1})-g(\varepsilon_i)|^2\right)^{\frac{1}{2}}.
$$
The variation and oscillation give information about the convergence properties of $g$.

Given a family of bounded operator $\mathcal T=\{T_t\}_{t>0}$ in $L^p(\Omega,\mu)$ with $1\leq p <\infty$ and for a certain measure space $(\Omega,\mu)$, an important problem in analysis is the existence of the limit $\lim_{t \rightarrow t_0}T_t(f)(x)$ with $t_0 \in [0,+\infty]$ and $x \in \Omega$ and its speed of convergence. In order to study those questions the following variation and oscillation operators can be considered
$$
V_\rho(\mathcal T)(f)(x)=V_\rho(t\rightarrow T_t(f)(x)),\quad x \in \Omega,
$$
where $\rho  > 2$, and
$$
O(\mathcal T,\{t_j\}_{j\in \mathbb N})(f)(x)= O(t \rightarrow T_t(f)(x),\{t_j\}_{j\in \mathbb N}),\quad x \in \Omega,
$$
being $\{t_j\}_{j\in \mathbb N} \subset (0,\infty)$ a sequence as above.

Variation inequalities have been investigated in probability,  ergodic theory and harmonic analysis. The first variation inequality is due to L\'epingle (\cite{Le}) for martingales. Later Bourgain (\cite{Bo}) proved a variation estimates for ergodic averages that inspired the interest of a number of authors in oscillation and variation inequalities for collections of operators in ergodic theory (\cite{JKRW} and \cite{JRW}) and harmonic analysis (\cite{CJRW1}, \cite{CJRw2}, \cite{Jo}, \cite{JSW}, \cite{LeMXu}, \cite{MTX1}, \cite{MTX2}, \cite{MTZ} and \cite{Mu1}).

The local-global strategy and the arguments developed in the proof of Theorem \ref{Th1.1} allow us to get the following variation and oscillation inequalities in our context.

\begin{thm}\label{Th1.3}
Let $\nu>0$, $k \in \mathbb N$ and $\alpha \in \mathbb N^d$. The operators
$$
V_\rho(\{t^{k+\widehat{\alpha}}\partial_t^k\partial_x^\alpha P_t^\nu\}_{t>0})
$$
and
$$
O(\{t^{k+\widehat{\alpha}}\partial_t^k\partial_x^\alpha P_t^\nu\}_{t>0}, \{t_j\}_{j\in \mathbb{N}}),
$$
where $\{t_j\}_{j\in \mathbb{N}}\subset (0,\infty)$ is decreasing and such that $\displaystyle \lim_{j\rightarrow \infty}t_j=0$, are bounded
\begin{enumerate}
\item from $L^p(\mathbb R^d,\gamma_\infty)$ into itself, for every $1<p<\infty$, when $\alpha =0$,
\item from $L^1(\mathbb R^d,\gamma_\infty)$ into $L^{1,\infty}(\mathbb R^d,\gamma_\infty)$ when $\widehat\alpha \leq 2$.
\end{enumerate}
\end{thm}

The result proved in \cite[Theorem 1.2]{HMMT} appears as a special case of the last theorem when the symmetric Ornstein-Uhlenbeck is considered with $k=0$ and $\alpha =0$.

The paper is organized as follow. Theorem \ref{Th1.1} is proved in Section \ref{section2} and Section \ref{section3}. In Section \ref{section2} we are concerned with $1<p<\infty$. The case $p=1$ is considered in Section \ref{section3}. Theorems \ref{Th1.2} and \ref{Th1.3} are proved in Section \ref{section4} and Section \ref{section5}, respectively.

Throughout this paper by $c$ and $C$ we denote always positive constants that can change in each occurrence.

\section{Proof of Theorem \ref{Th1.1} for $1<p<\infty$}\label{section2}
Assume that $1<p<\infty$. We consider firstly that $\alpha=0$ and $k \in \mathbb N$, $k \geq 1$. We have that
\begin{equation}\label{(2.0)}
\partial_t h_t(x,y)= \mathcal L_x h_t(x,y),\quad x,y\in \mathbb R^d,\;x \not= y.
\end{equation}

According to \cite[Lemma 4.1]{CCS3} we can write, for every $x,y \in \mathbb R^d$, $t>0$, and $i,j=1,\ldots,d$,
\begin{equation}\label{(2.1)}
\partial_{x_j} h_t(x,y)= h_t(x,y)P_j(t,x,y)
\end{equation}
where
$$
P_j(t,x,y)= \langle Q_\infty^{-1}x,e_j\rangle + \langle Q_t^{-1}e^{tB}e_j, y-D_tx\rangle,
$$
and
\begin{equation}\label{(2.2)}
\partial^2_{x_ix_j}h_t(x,y)=h_t(x,y)(P_i(t,x,y)P_j(t,x,y)+\Delta_{i,j}(t)),
\end{equation}
where
$$
\Delta_{i,j}(t)=-\langle e_j,e^{tB^*}Q_t^{-1}e^{tB}e_i\rangle.
$$
Here, for every $j=1,\ldots,d$, $e_j=(e_{j\ell})_{\ell=1}^d$ being $e_{j\ell}=0$, $\ell=1,\ldots,d$, $ \ell\not=j$, and $e_{jj}=1$.

By combining \cite[(2.10) and (4.5)]{CCS3} and \eqref{(2.1)} we get, for every $j=1,\ldots,d$,
\begin{align}\label{(2.3)}
\left|\partial_{x_j}h_t(x,y)\right| &\leq C\frac{e^{R(x)}}{t^\frac{d}{2}}e^{ -c\frac{|y-D_tx|^2}{t}} \left(|x|+\frac{|y-D_tx|}{t}\right)\nonumber\\ 
&\leq C\frac{e^{R(x)}}{t^{\frac{d}{2}}}e^{-c\frac{{|y-D_tx|^2}}{t}}\left(|x|+\frac{1}{\sqrt t}\right),\quad x,y \in \mathbb R^d\mbox{ and }0<t<1.
\end{align}
By \eqref{(2.1)} and \cite[(2.11), (4.5) and Lemma 2.1]{CCS3} it follows that, for every $j=1,\ldots,d$,
\begin{eqnarray}\label{(2.4)}
|\partial_{x_j}h_t(x,y)| &\leq& Ce^{R(x)}e^{-c|D_{-t}y-x|^2}(e^{-ct}|D_{-t}y -x|+|D_{-t}y|) \nonumber\\
& \leq & Ce^{R(x)}e^{-c|D_{-t}y-x|^2} e^{-ct}(1+|y|),\quad  x,y\in \mathbb R^d\mbox{ and } t\geq 1.
\end{eqnarray}

By using now (\ref{(2.2)}) and \cite[(2.10), (2.11), (4.5), (4.6) and Lemma 2.1]{CCS3} we obtain, for every $i,j=1,\ldots,d$,
\begin{equation}\label{(2.5)}
\left|\partial^2_{x_ix_j}h_t(x,y)\right| \leq C\frac{e^{R(x)}}{t^{\frac{d}{2}}}e^{-c\frac{|y-D_tx|^2}{t}}\left(|x|+\frac{1}{\sqrt{t}}\right)^2 , \quad x,y \in \mathbb R^d\mbox{ and }0<t<1,
\end{equation}
and
\begin{equation}\label{(2.6)}
\left|\partial^2_{x_ix_j}h_t(x,y)\right| \leq Ce^{R(x)}e^{-c|D_{-t}y-x|^2}e^{-ct}(1+|y|)^2, \quad x,y \in \mathbb R^d\mbox{ and }t\geq 1.
\end{equation}
By taking in account \eqref{(2.0)}, 
\begin{equation}\label{duxd2x}
    \begin{split}
        \partial_u h_u(x,y) & = 
        -\frac{1}{2} \textup{div}_x(Q\nabla_x h_u(x,y)) - \langle \nabla_x h_u(x,y), Bx\rangle
        \\ & =         \sum_{i,j=1}^{d} (c_{i,j}\partial^2_{x_ix_j}h_u(x,y) +
        d_{i,j} x_i\partial_{x_j}h_u(x,y))
        ,\quad 
        x,y\in\mathbb{R}^d
        \text{ and } u>0,
    \end{split}
\end{equation}
for certain $c_{i,j}$ and $d_{i,j}\in\mathbb{R}$, $i,j=1,\dots,d$. We deduce from \eqref{(2.3)} and \eqref{(2.5)} that
\begin{equation}\label{(2.7)}
|\partial_th_t(x,y)|\leq C\frac{e^{R(x)}}{t^{\frac{d}{2}+1}}e^{-c\frac{|y-D_tx|^2}{t}}(1+|x|)^2, \quad x,y\in\mathbb R^d\;\mbox{and}\;0<t<1,
\end{equation}
and from \eqref{(2.4)} and \eqref{(2.6)} that
\begin{equation}\label{(2.8)}
|\partial_th_t(x,y)|\leq C e^{R(x)}e^{-c|D_{-t}y-x|^2}e^{-ct}(1+|y|)^2, \quad  x,y\in\mathbb R^d\mbox{ and } t\geq 1.
\end{equation}

Let $f \in L^p(\mathbb R^d,\gamma_\infty)$. By \eqref{(2.7)} it follows that
\begin{equation}\label{(2.9)}
\int_{\mathbb R^d}|\partial_th_t(x,y)||f(y)|d\gamma_\infty(y) \leq C\frac{e^{R(x)}}{t^{\frac{d}{2}+1}}(1+|x|)^2\|f\|_{L^p(\mathbb R^d,\gamma_{\infty})},\quad x\in \mathbb R^d \mbox{ and } 0<t<1.
\end{equation}

From \eqref{(2.8)} we get
\begin{eqnarray}\label{(2.10)}
\int_{\mathbb R^d}|\partial_th_t(x,y)||f(y)|d\gamma_{\infty}(y) &\leq & Ce^{R(x)}e^{-ct}\int_{\mathbb R^d}|f(y)|(1+|y|)^2d\gamma_\infty(y)\\
&\leq & Ce^{R(x)}e^{-ct}\|f\|_{L^p(\mathbb R^d,\gamma_\infty)}, \quad x\in \mathbb R^d\mbox{ and } t>1. \nonumber
\end{eqnarray}

On the other hand, we have that
$$
P^\nu_t(f)(x)=\frac{1}{\Gamma(\nu)}\int_0^\infty e^{-s}s^{\nu-1}\mathcal H_{\frac{t^2}{4s}}(f)(x) ds, \quad x\in \mathbb R^d\mbox{ and }t>0.
$$

According to \eqref{(2.9)} and \eqref{(2.10)} we obtain
\begin{align*}
\int_0^\infty e^{-s}s^{\nu-1}\left|\partial_t\mathcal H_{\frac{t^2}{4s}}(f)(x)\right| ds &= \frac{t}{2}\int_0^\infty e^{-s}s^{\nu-2}|\partial_u\mathcal H_u(f)(x)_{|u=\frac{t^2}{4s}}|ds \\
& \hspace{-3cm}\leq Ct\left(\int_0^{\frac{t^2}{4}}e^{-s}s^{\nu-2}e^{-c\frac{t^2}{s}}ds+\int_{\frac{t^2}{4}}^\infty e^{-s}s^{\nu-2}\left(\frac{s}{t^2}\right)^{\frac{d}{2}+1}ds\right)e^{R(x)}(1+|x|)^2\|f\|_{L^p(\mathbb R^d,\gamma_\infty)}\\
&\hspace{-3cm} \leq Ct\left(\frac{1}{t^2}\int_0^{\frac{t^2}{4}} e^{-s}s^{\nu -1}ds+\frac{1}{t^{d+2}}\int_{\frac{t^2}{4}}^\infty e^{-s}s^{\nu +\frac{d}{2}-1}ds\right)e^{R(x)}(1+|x|)^2\|f\|_{L^p(\mathbb R^d,\gamma_\infty)}\\
&\hspace{-3cm}\leq Ce^{R(x)}\frac{(1+|x|)^2}{t}\Big(1+\frac{1}{t^d}\Big)\|f\|_{L^p(\mathbb R^d,\gamma_\infty)},\quad x \in \mathbb R^d\mbox{ and }t>0.
\end{align*}

Derivation under the integral sign is justified and we can write
$$
\partial_tP_t^\nu(f)(x) =\frac{1}{\Gamma(\nu)}\int_0^\infty e^{-s}s^{\nu-1}\partial_t\mathcal H_{\frac{t^2}{4s}}(f)(x)ds,\quad x\in \mathbb R^d\mbox{ and }t>0.
$$
It follows that, for every $x\in \mathbb{R}^d$ and $t>0$,
$$
\partial_tP_t^\nu(f)(x) = \frac{t}{2\Gamma(\nu)}\int_0^\infty e^{-s}s^{\nu-2}\partial_u\mathcal H_u(f)(x)_{|u=\frac{t^2}{4s}}ds 
= \frac{t^{2\nu-1}}{2^{2\nu-1}\Gamma(\nu)}\int_0^\infty \frac{e^{-\frac{t^2}{4u}}}{u^\nu}\partial_u\mathcal H_u(f)(x)du.
$$
As above we can justify the derivation under the integral sign and we get
\begin{eqnarray}\label{alpha0}
t^k\partial_t^k P_t^\nu(f)(x) &=& \frac{t^k}{2^{2\nu-1}\Gamma(\nu)}\int_0^\infty \frac{\partial_t^{k-1}\left[t^{2\nu-1}e^{-\frac{t^2}{4u}}\right]}{u^\nu}\partial_u\mathcal H_u(f)(x)du\nonumber\\
&=&\frac{2}{\Gamma(\nu)}\int_0^\infty [s^k\partial_s^{k-1}\mathfrak{h}_\nu(s)]_{|s=\frac{t}{2\sqrt u}}\partial_u\mathcal H_u(f)(x)du\nonumber\\
&=& \frac{4}{\Gamma(\nu)}\int_0^\infty s^{k-1}\partial_s^{k-1}\mathfrak{h}_\nu(s) \big[u\partial_u\mathcal H_u(f)(x)\big]_{|u=\frac{t^2}{4s^2}}ds, \quad x\in \mathbb R^d\mbox{ and }t>0,
\end{eqnarray}
where $\mathfrak{h}_\nu(s)=s^{2\nu-1}e^{-s^2}$, $s>0$.

By using Minkowski inequality we get
\begin{align*}
\left(\int_0^\infty |t^k\partial^k_tP_t^\nu(f)(x)|^2\frac{dt}{t}\right)^{\frac{1}{2}} & \leq C\int_0^\infty s^{k-1}|\partial_s^{k-1}\mathfrak{h}_\nu(s)|\Big\|u\partial_u\mathcal H_u(f)(x)_{|u=\frac{t^2}{4s^2}}\Big\|_{L^2((0,\infty ),\frac{dt}{t})}ds\\
&\leq CG^1(f)(x) \int_0 ^\infty  s^{k-1}|\partial^{k-1}_s\mathfrak{h}_\nu(s)|ds ,\quad x\in \mathbb R^d,
\end{align*}
where
$$
G^1(f)(x) = \left(\int_0^\infty |u\partial_u\mathcal H_u(f)(x)|^2\frac{du}{u}\right)^{\frac{1}{2}},\quad x\in \mathbb R^d.
$$
Since $\displaystyle \int_0^\infty s^{k-1}|\partial_s^{k-1}\mathfrak{h}_\nu(s)|ds<\infty$ (note that $s^\ell|\partial _s^\ell \mathfrak{h}_\nu (s)|\leq C\mathfrak{h}_\nu (s/2)$, $s>0$ and $\ell \in \mathbb{N}$) we conclude that
$$
g^\nu_{k,0}(f) \leq CG^1(f).
$$
As it was mentioned in Section \ref{S1} the square function $G^1$ is bounded from $L^p(\mathbb R^d,\gamma_\infty)$ into itself. Then, $g^\nu_{k,0}$ is bounded from $L^p(\mathbb R^d,\gamma_\infty)$ into itself.

Suppose now that $\alpha \in \mathbb N^d \setminus \{0\}$. By iterating the formulas in \cite[Lemma 4.1]{CCS3} and by using  \cite[(2.10), (2.11), (4.5), (4.6) and Lemma 2.1]{CCS3} we obtain
 \begin{equation}\label{(2.10.1)}
 |\partial_x^\alpha h_t(x,y)| \leq C\sum_{n=0}^{[\frac{\widehat\alpha}{2}]}\frac{1}{t^n}\left(|x|+\frac{|y-D_tx|}{t}\right)^{\widehat\alpha - 2n}\frac{e^{R(x)}}{t^{\frac{d}{2}}}e^{-c\frac{|y-D_tx|^2}{t}}, \quad x,y\in \mathbb R^d\mbox{ and }0<t< 1,
 \end{equation}
 and
 \begin{equation}\label{(2.11)}
  |\partial_x^\alpha h_t(x,y)| \leq C\sum_{n=0}^{[\frac{\widehat\alpha}{2}]}(|D_{-t}y-x|+|y|)^{\widehat\alpha - 2n} e^{-ct}e^{R(x)}e^{-c|x-D_{-t}y|^2}, \quad x,y\in \mathbb R^d\mbox{ and }t\geq 1.
 \end{equation}
 By proceeding as above we can see that the derivation under the integral sign is justified and we can write
 \begin{align}\label{tkalpha}
t^{k+\widehat{\alpha}}\partial_t^k\partial_x^\alpha P^\nu_t(f)(x) &= \frac{t^{k+\widehat{\alpha}}}{4^\nu\Gamma(\nu)}\int_0^\infty \partial_t^k \left(t^{2\nu}e^{-\frac{t^2}{4u}}\right) \partial_x^\alpha \mathcal H_u(f)(x)\frac{du}{u^{\nu+1}}\nonumber\\
 &=\frac{2^{\widehat{\alpha}}}{\Gamma(\nu)}\int_0^\infty [s^{k+\widehat{\alpha}}\partial_s^k{\mathfrak g}_\nu(s)]_{|s=\frac{t}{2\sqrt u}}\partial_x^\alpha \mathcal H_u(f)(x)u^{\frac{\widehat{\alpha}}{2}-1}du\\
 &=\frac{2^{\widehat{\alpha}+1}}{\Gamma(\nu)}\int_0^\infty s^{k+\widehat{\alpha}-1}\partial_s^k{\mathfrak g}_\nu(s)\big[u^{\frac{\widehat{\alpha}}{2}}\partial_x^\alpha \mathcal H_u(f)(x)\big]_{|u={\frac{t^2}{4s^2}}}ds,\quad x\in \mathbb R^d\mbox{ and }t>0,\nonumber
 \end{align}
 where ${\mathfrak g}_\nu(s)=s\mathfrak{h}_\nu (s)=s^{2\nu}e^{-s^2}$, $s>0$.

 Minkowski inequality leads to
 \begin{eqnarray*}
 \left(\int_0^\infty |t^{k+\widehat\alpha}\partial_t^k\partial^\alpha_xP_t^\nu(f)(x)|^2\frac{dt}{t}\right)^{\frac{1}{2}} &\leq & C\int_0^\infty s^{k+\widehat{\alpha}-1}|\partial_s^k\mathfrak{g}_\nu(s)|\left\|\big[u^{\frac{\widehat{\alpha}}{2}}\partial_x^\alpha \mathcal H_u(f)(x)\big]_{|u={\frac{t^2}{4s^2}}}\right\|_{L^2((0,\infty),\frac{dt}{t})}\\
 & \leq & CG_\alpha(f)(x)\int_0^\infty s^{k+\widehat{\alpha}-1}|\partial_s^k\mathfrak{g}_\nu(s)|ds ,\quad x \in \mathbb R^d,
 \end{eqnarray*}
 where
 $$
 G_\alpha(f)(x) = \left(\int_0^\infty |u^{\frac{\widehat\alpha}{2}}\partial_x^\alpha\mathcal H_u(f)(x)|^2\frac{du}{u}\right)^{\frac{1}{2}},\quad x \in \mathbb R^d.
 $$
Since $s^{\ell +\sigma}|\partial _s^\ell \mathfrak{g}_\nu (s)|\leq C\mathfrak{g}_\nu (s/2)$, $s>0$, $\ell \in \mathbb{N}$ and $\sigma \geq 0$, we have that $\displaystyle \int_0^\infty s^{k+\widehat\alpha-1} |\partial_s^k{\mathfrak g}_\nu(s)|ds < \infty$ and we conclude that
$$
g^\nu_{k,\alpha}(f) \leq C G_\alpha(f).
$$
Our next objective is to see that the square function $G_\alpha$ is bounded from $L^p(\mathbb R^d,\gamma_\infty)$ into itself.

According to \cite[Proposition 2.3]{MN2} the $\alpha$-Riesz transform $\mathcal R_\alpha$ is bounded from $L^p(\mathbb R^d,\gamma_\infty)$ into itself. Then, the set $\{t^{\frac{\widehat\alpha}{2}}\partial_x^\alpha \mathcal H_t\}_{t>0}$ of operators is $R$-bounded in $L^p(\mathbb R^d,\gamma_\infty)$. Indeed, let $t>0$. We consider the function $\Phi_t(z)=z^{\frac{\widehat\alpha}{2}}e^{-tz}$, $z\in \mathbb C\setminus (-\infty,0]$. Thus, $\Phi_t$ is analytic en $\mathbb C\setminus(-\infty,0]$ and the function $\Psi_t(z)=(1+z)^{\widehat\alpha}z^{-\frac{\widehat\alpha}{2}}\Phi_t(z)$, $z\in \mathbb C\setminus (-\infty,0]$, is analytic in $\mathbb C \setminus(-\infty,0]$ and bounded in $S_\mu=\{z\in \mathbb C:|{\rm Arg }\,z| <\mu\} \setminus \{0\}$, for every $\mu \in (0,\frac{\pi}{2})$. Since $\mathcal L_p$ is a sectorial and one to one operator in $L^p_0(\mathbb R^d,\gamma_\infty)$ with sharp angle $\theta_p$ we define the operator
$\Phi_t(\mathcal L_p)$ (see \cite[(2.1)]{CDMY}) as follows
\begin{equation}\label{(2.12)}
\Phi_t(\mathcal L_p) = \frac{1}{2\pi i}\int_{\Gamma_\sigma}(zI - \mathcal L_p)^{-1}\Phi_t(z)dz,
\end{equation}
where $\sigma\in (\theta_p,\frac{\pi}{2})$ and
$$
\Gamma_\sigma(s)=\left\{\begin{array}{ll}
-se^{-i\sigma}, & -\infty < s \leq 0,\\
se^{i\sigma}, & 0 \leq s < +\infty.
\end{array}\right.
$$
Note that integral in \eqref{(2.12)} converges in the Bochner sense with respect to the space of the bounded operators from $L^p_0(\mathbb R^d,\gamma_\infty)$ into itself.

Let $f \in L^p(\mathbb R^d,\gamma_\infty)$. Since $L^p(\mathbb R^d,\gamma_\infty)$ is the direct sum of $\mathcal N(\mathcal L_p)$, the kernel of $\mathcal L_p$, and $L^p_0(\mathbb R^d,\gamma_\infty)$ we can write $f=f_1+f_2$ where $f_1\in \mathcal N(\mathcal L_p)$ and $f_2 \in L^p_0(\mathbb R^d,\gamma_\infty)$.
Since $f_1$ is constant and $\{\mathcal{H}_t\}_{t>0}$ is conservative we have that
$$
\partial^\alpha_x\mathcal H_t(f) = \partial_x^\alpha\mathcal H_t(f_2).
$$
Furthermore, since $\gamma_\infty$ is a invariant measure for $\{\mathcal{H}_t\}_{t>0}$, $\mathcal H_t(f_2)\in L^p_0(\mathbb{R}^d, \gamma_\infty)$. By using \cite[(ii), p. 56]{CDMY} it follows that
$$
\partial_x^\alpha \mathcal H_t(f)=\partial^\alpha_x\mathcal L^{-\frac{\widehat\alpha}{2}}_p\Phi_t(\mathcal L_p)f_2.
$$
Then,
$$
t^{\frac{\widehat\alpha}{2}}\partial_x^\alpha \mathcal H_t(f)= \mathcal R_\alpha\Phi_1(t\mathcal L_p)\Pi_0(f).
$$
Suppose that $t_j>0$ and $f_j \in L^p(\mathbb R^d,\gamma_\infty)$, $j=1,\ldots, n$, $n\in \mathbb N$, $n\geq 1$. By \cite[Proposition 2.3]{MN2} we deduce inspired in \cite[Proposition 2.1]{CO} that
\begin{eqnarray*}
\mathbb E \left\| \sum_{j=1}^nr_jt^{\frac{\widehat\alpha}{2}}\partial_x^\alpha \mathcal H_{t_j}(f_j)\right\|_{L^p(\mathbb R^d,\gamma_\infty)}& = &\mathbb E\left\|\mathcal R_\alpha\sum_{j=1}^nr_{j_1}\Phi_1(t_j\mathcal L_p)\Pi_0(f_j)\right\|_{L^p(\mathbb R^d,\gamma_\infty)}\\
& \leq & C\mathbb E\left\|\sum_{j=1}^nr_j\Phi_1(t_j\mathcal L_p)\Pi_0(f_j)\right\|_{L^p(\mathbb R^d,\gamma_\infty)}.
\end{eqnarray*}
By $H^\infty(S_\mu)$ we denote the space of bounded and analytic functions in $S_\mu$ with $0<\mu<\pi$. On $H^\infty(S_\mu)$ we consider the norm $\|\cdot\|_{\infty,\mu}$ defined by
$$
\|f\|_{\infty,\mu}= \sup_{z\in S_\mu}|f(z)|,\quad f \in H^\infty(S_\mu).
$$
The set $\left\{t^{\frac{\widehat\alpha}{2}}\Phi_t(\mathcal L_p)\right\}_{t>0}$ is $R$-bounded in $H^\infty(S_\mu)$ for every $0<\mu < \frac{\pi}{2}$. According to \cite[Theorem 10]{CD} and \cite[Theorem 10.3.4]{HVNVW} the set $\left\{t^{\frac{\widehat\alpha}{2}}\Phi_t(\mathcal L_p)\right\}_{t>0}$ is $R$-bounded in $L^p(\mathbb R^d,\gamma_\infty)$. It follows that
\begin{eqnarray*}
\mathbb E\left\|\sum_{j=1}^nr_j\Phi_1(t_j\mathcal L_p)\Pi_0(f_j)\right\|_{L^p(\mathbb R^d,\gamma_\infty)} & \leq & C\mathbb E\left\|\sum_{j=1}^nr_j\Pi_0(f_j)\right\|_{L^p(\mathbb R^d,\gamma_\infty)}\\
&\leq& C \mathbb E\left\| \Pi_0(\sum_{j=1}^nr_jf_j)\right\|_{L^p(\mathbb R^d,\gamma_\infty)}\\
&\leq& C \mathbb E\left\|\sum_{j=1}^nr_jf_j\right\|_{L^p(\mathbb R^d,\gamma_\infty)}.
\end{eqnarray*}
We conclude that the set $\left\{t^{\frac{\widehat\alpha}{2}}\partial_x^\alpha \mathcal H_t\right\}_{t>0}$ is $R$-bounded in $L^p(\mathbb R^d,\gamma_\infty)$.

We continue the argument adapting a procedure developed in the proof of \cite[Proposition 5.1]{CO}. Let $f \in L^p(\mathbb R^d,\gamma_\infty)$. By partial integration we get
\begin{multline}
\int_0^\infty|t^{\frac{\widehat\alpha}{2}}\partial_x^\alpha \mathcal H_t(f)(x)|^2\frac{dt}{t} = \int_0^\infty t^{\widehat\alpha-1}|\partial_x^\alpha \mathcal H_t(f)(x)|^2dt \\
=\frac{1}{\widehat\alpha}\left(\lim_{t \rightarrow +\infty}t^{\widehat\alpha}|\partial_x^\alpha\mathcal H_t(f)(x)|^2 - 2 \mbox{\rm Re}\,\int_0^\infty t^{\widehat\alpha}\partial_x^\alpha \partial_t \mathcal H_t(f)(x)\overline{\partial_x^\alpha \mathcal H_t(f)(x)} dt\right), \quad x\in \mathbb R^d.
\end{multline}
 By \eqref{(2.11)} we deduce that
 $$
 t^{\widehat\alpha}|\partial_x^\alpha \mathcal H_t(f)(x)|^2  \leq Ct^{\widehat\alpha}e^{-ct}e^{R(x)}\|f\|_{L^p(\mathbb R^d,\gamma_\infty)},\quad x\in \mathbb R^d\mbox{ and }t>1.
 $$
 Then,
 \begin{align*}
 \int_0^\infty\left|t^{\frac{\widehat\alpha}{2}} \partial_x^\alpha \mathcal H_t(f)(x)\right|^2\frac{dt}{t} &\leq C \int_0^\infty t^{\widehat\alpha}\left|\partial_x^\alpha \partial_t\mathcal H_t(f)(x)\right|\left|\partial_x^\alpha\mathcal H_t(f)(x)\right|dt\\
 &\leq \left(\int_0^\infty \left|t^{\frac{\widehat\alpha}{2}} \partial_x^\alpha \mathcal H_t(f)(x)\right|^2\frac{dt}{t} \right)^{\frac{1}{2}}\left(\int_0^\infty \left|t^{\frac{\widehat\alpha}{2}+1}\partial_x^\alpha \partial_t\mathcal H_t(f)(x)\right|^2\frac{dt}{t}\right)^{\frac{1}{2}},\quad x\in \mathbb R^d.
 \end{align*}
 By using the semigroup property it follows that
 \begin{eqnarray*}
 G_\alpha(f)(x) &\leq & C \left(\int_0^\infty \left|t^{\frac{\widehat\alpha}{2}+1}\partial_x^\alpha \partial_t\mathcal H_t(f)(x)\right|^2\frac{dt}{t}\right)^{\frac{1}{2}}\\
 & \leq & C\left(\int_0^\infty\left|t^{\frac{\widehat\alpha}{2}}\partial_x^\alpha \mathcal H_{\frac{t}{2}}\left(t\partial_t\mathcal H_{\frac{t}{2}}(f)\right)(x)\right|^2\frac{dt}{t}\right)^{\frac{1}{2}},\;x\in \mathbb R^d.
 \end{eqnarray*}
 As it was seen in Section \ref{S1}, the Littlewood-Paley function $G^1$ is bounded from $L^p(\mathbb R^d,\gamma_\infty)$ into itself. Hence, for almost every $x\in \mathbb R^d$, the function $t  \rightarrow t\partial_t\mathcal H_{\frac{t}{2}}(f)(x)$ is in $L^2((0,\infty), \frac{dt}{t})$. Since the set $\{t^{\frac{\widehat\alpha}{2}}\partial_x^\alpha \mathcal H_t\}_{t>0}$ is $R$-bounded in $L^p(\mathbb R^d,\gamma_\infty)$, by using \cite[Lemma 2.3]{CO} we get
 $$
 \|G_\alpha(f)\|_{L^p(\mathbb R^d,\gamma_\infty)} \leq C\|G^1(f)\|_{L^p(\mathbb R^d,\gamma_\infty)} \leq C\|f\|_{L^p(\mathbb R^d,\gamma_\infty)}.
 $$
 We conclude that $g^\nu_{k,\alpha}$ is bounded from $L^p(\mathbb R^d,\gamma_\infty)$ into itself.

\section{Proof of Theorem \ref{Th1.1} for $p=1$.}\label{section3}
In order to prove that $g_{k,\alpha}^\nu$ is bounded from $L^1(\mathbb{R}^d,\gamma_\infty)$ into $L^{1,\infty}(\mathbb{R}^d,\gamma_\infty)$ we use the local-global technique.

Let $A>0$. We define the local region $L_A$ by
$$
L_A=\Big\{(x,y)\in \mathbb{R}^d\times \mathbb R^d: |x-y|\leq \frac{A}{1+|x|}\Big\},
$$
and the global region $G_A$ by $G_A=(\mathbb{R}^d\times \mathbb R^d)\setminus L_A$. The value of $A$ will be fixed later.

We choose a smooth function $\varphi_{\mbox{\tiny A}}$ on $\mathbb{R}^d\times \mathbb R^d$ such that $0\leq \varphi_{\mbox{\tiny A}}\leq 1$, $\varphi_{\mbox{\tiny A}} (x,y)=1$, $(x,y)\in L_A$, $\varphi_{\mbox{\tiny A}} (x,y)=0$, $(x,y)\in G_{2A}$, and there exists $C>0$ such that
$$
\sum_{i=1}^d(|\partial _{x_i}\varphi_{\mbox{\tiny A}} (x,y)|+|\partial _{y_i}\varphi_{\mbox{\tiny A}} (x,y)|)\leq \frac{C}{|x-y|},\quad x,y\in \mathbb{R} ^d,\;x\not=y.
$$

Given an operator $T$ defined on $L^p(\mathbb{R}^d,\gamma_\infty )$, $1\leq p<\infty$, we consider the local and global operators given by $T_{\rm loc}(f)(x)=T(\varphi_{\mbox{\tiny A}} (x,\cdot )f)(x)$, $x\in \mathbb{R}^d$, and $T_{\rm glob}=T-T_{\rm loc}$, respectively. 

In the study of the $L^p$-boundedness for our operators of local type we make use of the operator $\mathscr{S}_\eta$, with $\eta>0$, defined by  
$$
\mathscr{S}_\eta (g)(x)=\int_{(x,y)\in L_\eta}\frac{1+|x|}{|x-y|^{d-1}}g(y)dy,\quad x\in \mathbb R^d.
$$
Let $\eta>0$. We observe that $\mathscr{S}_\eta$ is a bounded operator from $L^p(\mathbb R^d,dx)$ into itself, for every $1\leq p\leq \infty$. Indeed, for every $x\in \mathbb R^d$, we can write
$$
\int_{(x,y)\in L_\eta}\frac{1+|x|}{|x-y|^{d-1}}dy\leq C\int_0^{\frac{\eta}{1+|x|}}(1+|x|)dr\leq C.
$$
Then
$$
\sup_{x\in \mathbb{R}^d}\int_{(x,y)\in L_\eta}\frac{1+|x|}{|x-y|^{d-1}}dy<\infty.
$$
On the other hand, when $y\in \mathbb R^d$, $|y|\leq 2\eta$, it follows that
$$
\int_{(x,y)\in L_\eta} \frac{1+|x|}{|x-y|^{d-1}}dx\leq \int_{|x-y|<\eta}\frac{1+|x-y|+|y|}{|x-y|^{d-1}}dx\leq (3\eta +1)\int_{|x-y|<\eta}\frac{dx}{|x-y|^{d-1}}\leq C,
$$
and if $|y|\geq 2\eta$ and $(x,y)\in L_\eta$, then, $|x|\geq |y|-|x-y|\geq |y|-\eta\geq |y|-\frac{|y|}{2}=\frac{|y|}{2}$. Hence, $|x-y|\leq \frac{2\eta}{1+|y|}$ when $(x,y)\in L_\eta$ and $|y|\geq 2\eta$, and we can write
$$
\int_{(x,y)\in L_\eta} \frac{1+|x|}{|x-y|^{d-1}}dx\leq \int_{|x-y|\leq \frac{2\eta}{2+|y|}}\frac{1+|x-y|+|y|}{|x-y|^{d-1}}dx\leq C\int_0^{\frac{2\eta}{2+|y|}}(1+\eta +|y|)dr\leq C,\quad y\in \mathbb R^d, \;|y|\geq 2\eta.
$$
Then,
$$
\sup_{y\in \mathbb{R}^d}\int_{(x,y)\in L_\eta}\frac{1+|x|}{|x-y|^{d-1}}dx<\infty.
$$
By using interpolation we conclude that $\mathscr S_\eta$ is bounded from $L^p(\mathbb R^d,dx)$ into itself, for every $1\leq p\leq \infty$. Since $\mathscr S_\eta$ is a local operator we can deduce that $\mathscr S_\eta$ is also bounded from $L^p(\mathbb{R}^d,\gamma_\infty )$ into itself, $1\leq p<\infty$.

To see this assertion let us consider the sequence $\{B_\ell=B(w_\ell,\frac{1}{20(1+|w_\ell|)})\}_{\ell\in \mathbb N}$ of balls given in \cite[Lemma 3.1]{GCMST2}, where $w_\ell \in \mathbb{R}^d$, $\ell \in \mathbb{N}$. We have that, for each $\sigma >0$ there exists $C>0$ such that 
$$
\frac{e^{R(w_\ell)}}{C}\leq e^{R(x)}\leq Ce^{R(w_\ell)},\quad x\in \sigma B_\ell\mbox{ and }\ell \in \mathbb{N}.
$$
 Here $\sigma B_\ell=B(w_\ell, \frac{\sigma}{20(1+|w_\ell|)})$, $\ell\in \mathbb{N}$ and $\sigma >0$. Indeed, for $\sigma >0$ we can write
\begin{align}\label{equiv}
    e^{|R(x)-R(w_\ell)|}&=e^{\frac{1}{2}\big||Q_\infty ^{-1/2}x|^2-|Q_\infty ^{-1/2}w_\ell|^2\big|} \leq e^{\frac{1}{2}|Q_\infty ^{-1/2}(x+w_\ell)||Q_\infty ^{-1/2}(x-w_\ell)|} \nonumber\\
   & \leq e^{c|x+w_\ell||x-w_\ell|} \leq e^{c(|x-w_\ell|^2+2|w_\ell||x-w_\ell|})\leq e^{c(\sigma ^2+2\sigma)}=C,\quad x\in \sigma B_\ell \mbox { and }\ell\in \mathbb N.
  \end{align}

We take $\sigma=\sigma (\eta)>0$ according to \cite[Lemma 3.1]{GCMST2} in such a way that for every $\ell \in \mathbb{N}$, if $x\in B_\ell$ and $(x,y)\in L_\eta$, then $y\in \sigma B_\ell$. Then, we have that $\mathscr{S}_\eta (f)(x)=\mathscr{S}_\eta (\mathcal{X}_{\sigma B_\ell}f)(x)$, $x\in B_\ell$ and $\ell\in \mathbb N$. 

Let $1\leq p<\infty$. According to the properties of the sequence $\{B_\ell\}_{\ell\in \mathbb N}$ (\cite[Lemma 3.1]{GCMST2}) we deduce that
\begin{align*}
    \int_{\mathbb R^d}|\mathscr{S}_\eta (f)(x)|^pd\gamma_\infty(x)&\leq C \sum_{\ell \in \mathbb{N}} \int_{B_\ell}|\mathscr{S}_\eta(\mathcal{X}_{\sigma B_\ell}f)(x)|^pe^{-R(x)}dx\leq C\sum_{\ell\in \mathbb N} e^{-R(w_\ell)}\int_{\mathbb R^d}|\mathscr{S}_\eta(\mathcal{X}_{\sigma B_\ell}f)(x)|^pdx\\
    &\leq C\sum_{\ell\in \mathbb N} e^{-R(w_\ell)}\int_{\sigma B_\ell}|f(y)|^pdy\leq  C\sum_{\ell\in \mathbb N}\int_{\sigma B_\ell}|f(y)|^pe^{-R(y)}dy\leq C\int_{\mathbb R^d}|f(y)|^pe^{-R(y)}dy.
\end{align*}

We note that the reasoning given above also works for any local operator $T_{\rm loc}$. Thus, given $1\leq p<\infty$, if $T_{\rm loc}$ is an operator bounded from $L^p(\mathbb{R}^d,dx)$ into itself, then is also bounded from $L^p(\mathbb{R}^d,\gamma_\infty )$ into itself. Moreover, in a similar way we can establish that if $T_{\rm loc}$ is bounded from $L^p(\mathbb{R}^d,dx)$ into $L^{p,\infty }(\mathbb{R}^d,dx)$, then it is also bounded from $L^p(\mathbb{R}^d,\gamma_\infty )$ into $L^{p,\infty }(\mathbb{R}^d,\gamma_\infty)$, $1\leq p<\infty$. We will turn to this argument in the proofs of our results.

Along this section we also consider the function $\mathfrak{m}(z)=\min\{1,|z|^{-2}\}$, $z\in \mathbb R^d\setminus\{0\}$, and $\mathfrak{m}(0)=1$. We observe that   $\mathfrak{m}(z)\sim (1+|z|)^{-2}$, $z\in \mathbb{R}^d$, an estimate that we will frequently use. 

\subsection{}
 Assume first that $0<\widehat{\alpha}\leq 2$ and $k\in \mathbb{N}$. According to \eqref{tkalpha}, for every $f\in L^p (\mathbb R^d,\gamma_\infty)$, $1\leq p<\infty$, we can write 
$$
t^{k+\widehat{\alpha}}\partial_t^k\partial_x^\alpha P_t^\nu (f)(x)=\int_{\mathbb{R}^d}P_{k,\alpha}^\nu (x,y,t)f(y)d\gamma_\infty (y),\quad x\in \mathbb R^d\mbox{ and }t>0,
$$
where
$$
P_{k,\alpha}^\nu (x,y,t)=\frac{2^{\widehat{\alpha}}}{\Gamma (\nu)}\int_0^\infty [s^{k+\widehat{\alpha}}\partial _s^k{\mathfrak g}_\nu(s)]_{\big|s=\frac{t}{2\sqrt{u}}}\partial _x^\alpha h_u(x,y)u^{\frac{\widehat{\alpha}}{2}-1}du,\quad x,y\in \mathbb{R}^d\mbox{ and }t>0.
$$
Here, as before, ${\mathfrak g}_\nu(s)=s^{2\nu }e^{-s^2}$, $s\in (0,\infty )$. 

We define the local part of $g_{k,\alpha}^\nu$ by
$$
g_{k,\alpha,{\rm loc}}^\nu (f)(x)=\left(\int_0^\infty \left|\int_{\mathbb R^d}P_{k,\alpha}^\nu (x,y,t)\varphi_{\mbox{\tiny A}} (x,y)f(y)d\gamma_\infty (y)\right|^2\frac{dt}{t}\right)^{1/2},\quad x\in \mathbb R^d,
$$
and the global part of $g_{k,\alpha}^\nu$ by
$$
g_{k,\alpha,{\rm glob}}^\nu (f)(x)=\left(\int_0^\infty \left|\int_{\mathbb R^d}P_{k,\alpha}^\nu (x,y,t)(1-\varphi_{\mbox{\tiny A}} (x,y))f(y)d\gamma_\infty (y)\right|^2\frac{dt}{t}\right)^{1/2},\quad x\in \mathbb R^d.
$$
It is clear that
$$
g_{k,\alpha}^\nu (f)\leq g_{k,\alpha,{\rm loc}}^\nu (f)+g_{k,\alpha,{\rm glob}}^\nu (f),\quad f\in L^p(\mathbb{R}^d,\gamma_\infty) \;\;\;(1\leq p<\infty).
$$

We study $g_{k,\alpha,{\rm loc}}^\nu$ and $g_{k,\alpha,{\rm glob}}^\nu$ separately. We take advantage of some estimates established in \cite{CCS3}. We also take into account that, for each $\ell \in \mathbb{N}$, $|s^\ell \partial _s^\ell \mathfrak{g}_\nu (s)|\leq C\mathfrak{g}_\nu (s/2)$, $s\in (0,\infty )$, from which it follows that $\|s^{\ell+\sigma}\partial _s^\ell \mathfrak{g}_\nu (s)\|_{L^2((0,\infty ),\frac{ds}{s})}\leq C$, $\sigma\geq 0$. 

\subsubsection{About the local part of $g_{k,\alpha}^\nu$ when $\widehat{\alpha}=1$}\label{S3.1.1}

Assume that $j\in \{1,\ldots,d\}$ and $\alpha =(\alpha_1,\ldots,\alpha _d)\in \mathbb{N}^d$ being $\alpha _j=1$ and $\alpha _i=0$, $i\in \{1,\ldots,d\}\setminus\{j\}$. We can write
$$
g_{k,\alpha,{\rm loc}}^\nu (f)(x)=\left\|\int_{\mathbb R^d}K_{k,j}^\nu(x,y,\cdot )\varphi_{\mbox{\tiny A}} (x,y)f(y)d\gamma_\infty (y)\right\|_{L^2((0,\infty ),\frac{dt}{t})},\quad x\in \mathbb{R}^d,
$$
where 
$$
K_{k,j}^\nu(x,y,t)=\frac{2}{\Gamma (\nu)}\int_0^\infty [s^{k+1}\partial _s^k{\mathfrak g}_\nu(s)]_{|s=\frac{t}{2\sqrt{u}}}\partial _{x_j}h_u(x,y)\frac{du}{\sqrt{u}},\quad x,y\in \mathbb{R}^d \mbox{ and }t>0.
$$
We consider the operator $S_{k,j}^\nu$ defined by
\begin{equation}\label{Snuj}
S_{k,j}^\nu(f)(x,t)=\int_{\mathbb R^d}S_{k,j}^\nu (x,y,t)f(y)d\gamma_\infty (y),\quad x\in \mathbb{R}^d\mbox{ and }t>0,
\end{equation}
where, for every $x,y\in \mathbb R^d$,
$$
S_{k,j}^\nu (x,y,t)=\frac{2}{\Gamma (\nu)}\int_0^\infty [s^{k+1}\partial _s^k{\mathfrak g}_\nu(s)]_{|s=\frac{t}{2\sqrt{u}}}\mathbb{S}_u^j(x,y)\frac{du}{\sqrt{u}},\quad t>0,
$$
and
\begin{align*}
\mathbb{S}_u^j(x,y)&=\left(\frac{{\rm det}\,Q_\infty}{u^d{\rm det}\,Q}\right)^{1/2}e^{R(x)}\partial _{x_j}e^{-\frac{1}{2u}\langle Q^{-1}(y-x),y-x\rangle }\\
&=\left(\frac{{\rm det}\,Q_\infty}{u^d{\rm det}\,Q}\right)^{1/2}e^{R(x)}\frac{1}{u}\langle Q^{-1}e_j,y-x\rangle e^{-\frac{1}{2u}\langle Q^{-1}(y-x),y-x\rangle },\quad u>0.
\end{align*}
Let us denote by $S_{k,j,{\rm loc}}^\nu$ the corresponding local operator, that is, $S_{k,j,{\rm loc}}^\nu (f)(x)=S_{k,j}^\nu (\varphi_{\mbox{\tiny A}} (x,\cdot)f)(x)$, $x\in \mathbb{R}^d$. In order to obtain the $L^p$-boundedness properties for $g_{k,\alpha,{\rm loc}}^\nu$ it is sufficient to establish the following two items:

\begin{itemize}
\item[(I1)] $S_{k,j,{\rm loc}}^\nu$ is bounded from $L^p(\mathbb R^d,\gamma _\infty)$ into $L^p_{L^2((0,\infty ),\frac{dt}{t})}(\mathbb R^d,\gamma_\infty)$, for every $1<p<\infty$, and from $L^1(\mathbb R^d,\gamma_\infty)$ into $L^{1,\infty}_{L^2((0,\infty ),\frac{dt}{t})}(\mathbb R^d,\gamma_\infty)$.
\item [(I2)] The operator given by $\mathcal{D}_{k,j}^\nu (f)(x)=g_{k,\alpha,{\rm loc}}^\nu (f)(x)-\|S_{k,j,{\rm loc}}^\nu(f)(x,\cdot)\|_{L^2((0,\infty ),\frac{dt}{t})}$, $x\in \mathbb{R}^d$, is bounded from $L^p(\mathbb{R}^d,\gamma_\infty)$ into itself, for each $1\leq p\leq \infty$.
\end{itemize}

First, we prove (I1). For that, we  decompose $S_{k,j,{\rm loc}}^\nu$ as follows 
\begin{align*}
S_{k,j,{\rm loc}}^\nu (f)(x,t)&=\int_{\mathbb R^d}S_{k,j}^\nu (x,y,t)[\varphi_{\mbox{\tiny A}} (x,y)-\varphi_{\mbox{\tiny A}} (Q^{-1/2}x,Q^{-1/2}y)]f(y)d\gamma_\infty (y)\\
&\quad +\int_{\mathbb R^d}S_{k,j}^\nu (x,y,t)\varphi_{\mbox{\tiny A}} (Q^{-1/2}x,Q^{-1/2}y)f(y)d\gamma_\infty (y)\\
&=S_{k,j,{\rm loc}}^{\nu,1}(f)(x,t) +S_{k,j,{\rm loc}}^{\nu,2}(f)(x,t) ,\quad x\in \mathbb R^d\mbox{ and }t>0.
\end{align*}

We can find $\eta >0$ such that $\frac{\eta ^{-1}}{1+|x|}\leq |x-y|\leq\frac{\eta }{1+|x|}$, provided that $\varphi_{\mbox{\tiny A}} (x,y)-\varphi_{\mbox{\tiny A}} (Q^{-1/2}x,Q^{-1/2}y)\not=0$. Then, we can write 
$$
\|S_{k,j,{\rm loc}}^{\nu,1}(f)(x,\cdot)\|_{L^2((0,\infty ),\frac{dt}{t})}\leq C\int_{\frac{\eta ^{-1}}{1+|x|}\leq |x-y|\leq\frac{\eta }{1+|x|}}\|S_{k,j}^\nu (x,y,\cdot )\|_{L^2((0,\infty ),\frac{dt}{t})}|f(y)|e^{-R(y)}dy. 
$$
We have that
\begin{align}\label{3.2'}
R(x)-R(y)&=\frac{1}{2}(|Q_\infty ^{-1/2}x|^2-|Q_\infty ^{-1/2}y|^2)\leq \frac{1}{2}|Q_\infty ^{-1/2}(x+y)||Q_\infty ^{-1/2}(x-y)|\nonumber\\
&\leq C|x+y||x-y|\leq C\frac{2|x|+|x-y|}{1+|x|}\leq C,\quad (x,y)\in L_\eta ,
\end{align}
and thus,  when $\frac{\eta ^{-1}}{1+|x|}\leq |x-y|\leq\frac{\eta }{1+|x|}$ we get
\begin{align*}
e^{-R(y)}\|S_{k,j}^\nu (x,y,\cdot )\|_{L^2((0,\infty ),\frac{dt}{t})}&\\
&\hspace{-3cm}\leq C\int_0^\infty \big\|[s^{k+1}\partial _s^k\mathfrak{g}_\nu(s)]_{|s=\frac{t}{2\sqrt{u}}}\big\|_{L^2((0,\infty ),\frac{dt}{t})}\Big|\frac{1}{u}\langle Q^{-1}e_j,y-x\rangle e^{-\frac{1}{2u}\langle Q^{-1}(y-x),y-x\rangle} \Big|\frac{du}{u^{\frac{d+1}{2}}}\\
&\hspace{-3cm}\leq C\int_0^\infty |Q^{-1/2}(y-x)|e^{-\frac{|Q^{-1/2}(y-x)|^2}{2u}}\frac{du}{u^{\frac{d+3}{2}}}\leq C\int_0^\infty \frac{e^{-c\frac{|y-x|^2}{u}}}{u^{\frac{d}{2}+1}}du\\
&\hspace{-3cm}=\frac{C}{|x-y|^d}\leq C\frac{1+|x|}{|x-y|^{d-1}}.
\end{align*}
Observe that $\|[s^{k+1}\partial _s^k\mathfrak{g}_\nu(s)]_{|s=\frac{t}{2\sqrt{u}}}\|_{L^2((0,\infty ),\frac{dt}{t})}=\|s^{k+1}\partial _s^k\mathfrak{g}_\nu(s)\|_{L^2((0,\infty ),\frac{ds}{s})}=C$, $u>0$.

Hence we conclude that
$$
\big\|S_{k,j,{\rm loc}}^{\nu,1}(f)(x,\cdot )\big\|_{L^2((0,\infty ),\frac{dt}{t})}\leq C\mathscr{S}_\eta(|f|)(x),\quad x\in \mathbb R^d,
$$
and then, $S_{k,j,{\rm loc}}^{\nu,1}$ is a bounded operator from $L^p(\mathbb{R}^d,dx)$ into $L^p_{L^2((0,\infty),\frac{dt}{t})}(\mathbb{R}^d,dx)$, for every $1\leq p\leq \infty$.

We now deal with $S_{k,j,{\rm loc}}^{\nu ,2}$. A straightforward change of variables allows us to write 
$$
S_{k,j,{\rm loc}}^{\nu ,2}(f)(x,t)=c\widetilde{S}_{k,j,{\rm loc}}^\nu(f(Q^{1/2}\cdot))(Q^{-1/2}x,t),\quad x\in \mathbb{R}^d\mbox{ and }t>0,
$$
where the operator $\widetilde{S}_{k,j,{\rm loc}}^\nu$ is defined by
$$
\widetilde{S}_{k,j,{\rm loc}}^\nu (g)(x,t)=\int_{\mathbb R^d}\widetilde{S}_{k,j}^\nu(x,y,t)\varphi_{\mbox{\tiny A}} (x,y)g(y)dy,\quad x\in \mathbb R^d\mbox{ and }t>0,
$$
and being $\widetilde{S}_{k,j}^\nu(x,y,t)=e^{-R(Q^{1/2}y)}S_{k,j}^\nu(Q^{1/2}x,Q^{1/2}y,t)$, $x,y\in \mathbb{R}^d$, $t>0$. Let us establish that the operator $\widetilde{S}_{k,j,{\rm loc}}^\nu$ is bounded from $L^p(\mathbb R^d,dx)$ into $L^p_{L^2((0,\infty ),\frac{dt}{t})}(\mathbb R^d,dx)$, for every $1<p<\infty$, and from $L^1(\mathbb R^d,dx)$ into $L^{1,\infty}_{L^2((0,\infty ),\frac{dt}{t})}(\mathbb R^d,dx)$. Consequently, we will have the same $L^p$-boundedness properties for ${S}_{k,j,{\rm loc}}^{\nu,2}$.

If $Q^{-1/2}=(c_{n\ell})_{n,\ell =1}^d$ with $c_{n\ell}\in \mathbb R$, $n,\ell=1,\ldots,d$, we have that
$$
\langle Q^{-1/2}e_j,z\rangle =\sum_{n=1}^dc_{nj}z_n, \quad z=(z_1,\ldots,z_d)\in \mathbb R^d.
$$
Then, we can write
\begin{align*}
\widetilde{ S}_{k,j}^\nu(x,y,t)&=\frac{2}{\Gamma (\nu)}\Big(\frac{{\rm det}\,Q_\infty}{{\rm det }\,Q}\Big)^{1/2}e^{R(Q^{1/2}x)-R(Q^{1/2}y)}\int_0^\infty [s^{k+1}\partial _s^k\mathfrak{g}_\nu(s)]_{|s=\frac{t}{2\sqrt{u}}}\langle Q^{-1/2}e_j,y-x\rangle e^{-\frac{|y-x|^2}{2u}}\frac{du}{u^{\frac{d+3}{2}}}\\
&=\frac{2}{\Gamma (\nu)}\Big(\frac{{\rm det}\,Q_\infty}{{\rm det}\,Q}\Big)^{1/2}e^{R(Q^{1/2}x)-R(Q^{1/2}y)}\sum_{n=1}^dc_{nj}(y_n-x_n)\int_0^\infty [s^{k+1}\partial _s^k\mathfrak{g}_\nu(s)]_{|s=\frac{t}{2\sqrt{u}}}e^{-\frac{|y-x|^2}{2u}}\frac{du}{u^{\frac{d+3}{2}}},
\end{align*}
and we get $\widetilde{S}_{k,j,{\rm loc}}^\nu (f)(x,t)=c\sum_{n=1}^dc_{nj}T_{k,n,{\rm loc}}^\nu(f)(x,t)$, $x\in \mathbb{R}^d$, $t>0$, where for every $n=1,\ldots,d$, the operator $T_{k,n}^\nu$ is given by
$$
T_{k,n}^\nu (f)(x,t)=\int_{\mathbb R^d}T_{k,n}^\nu (x,y,t)f(y)dy,\quad x\in \mathbb R^d\mbox{ and }t>0,
$$
and $T_{k,n,{\rm loc}}^\nu$ is defined in the usual way, and where, for every $x=(x_1,...,x_d)$, $y=(y_1,...,y_d)\in \mathbb R^d$ and $t>0$,
$$
T_{k,n}^\nu (x,y,t)=\frac{2}{\Gamma (\nu )}e^{R(Q^{1/2}x)-R(Q^{1/2}y)}\int_0^\infty [s^{k+1}\partial _s^k\mathfrak{g}_\nu(s)]_{|s=\frac{t}{2\sqrt{u}}}\partial _{x_n}\mathbb{W}_u(x-y)\frac{du}{\sqrt{u}}.
$$
Here and in the sequel $\{\mathbb W_u\}_{u>0}$ represents the Euclidean heat semigroup in $\mathbb R^d$, that is, $\mathbb{W}_u(z) = (2\pi u)^{-\frac{d}{2}} e^{-\frac{|z|^2}{2u}}$, $z\in\mathbb{R}^d$ and $u>0$.

Let $n=1,...,d$. To see that $T_{k,n,{\rm loc}}^\nu$ is bounded from $L^p(\mathbb R^d,dx)$ into $L^p_{L^2((0,\infty ),\frac{dt}{t})}(\mathbb R^d,dx)$, for every $1<p<\infty$, and from $L^1(\mathbb R^d,dx)$ into $L^{1,\infty}_{L^2((0,\infty ),\frac{dt}{t})}(\mathbb R^d,dx)$ we introduce the operator $\mathbb{T}_{k,n}^\nu$ as follows
$$
\mathbb T_{k,n}^\nu (f)(x,t)=\int_{\mathbb R^d}\mathbb{T}_{k,n}^\nu (x-y,t)f(y)dy, \quad x\in \mathbb R^d\mbox{ and }t>0,
$$
where
$$
\mathbb{T}_{k,n}^\nu (z,t)=\frac{2}{\Gamma (\nu )}\int_0^\infty [s^{k+1}\partial _s^k\mathfrak{g}_\nu(s)]_{\big|s=\frac{t}{2\sqrt{u}}}\partial _{z_n}\mathbb{W}_u(z)\frac{du}{\sqrt{u}},\quad z=(z_1,...,z_d)\in \mathbb{R}^d\mbox{ and }t>0.
$$
Note that $T_{k,n}^\nu(f)(x,t)=e^{R(Q^{1/2}x)}\mathbb{T}_{k,n}^\nu (e^{-R(Q^{1/2}\cdot )}f)(x,t)$, $x\in \mathbb{R}^d$ and $t>0$.

We will show below that the local operator $\mathbb{T}_{k,n,{\rm loc}}^\nu$ is bounded from $L^p(\mathbb R^d,dx)$ into $L^p_{L^2((0,\infty ),\frac{dt}{t})}(\mathbb R^d,dx)$, for every $1<p<\infty$, and from $L^1(\mathbb R^d,dx)$ into $L^{1,\infty}_{L^2((0,\infty ),\frac{dt}{t})}(\mathbb R^d,dx)$. From this fact, we can deduce the same $L^p$-properties for $T_{k,n,{\rm loc}}^\nu$. Effectively, let us consider the sequence $\{B_\ell\}_{\ell\in \mathbb N}$ of balls given in the proof of \eqref{equiv}. Following the same reasoning we obtain that
\begin{align*}
    e^{|R(Q^{1/2}x)-R(Q^{1/2}w_\ell)|}\leq C,\quad x\in \sigma B_\ell \mbox { and }\ell\in \mathbb N,
  \end{align*}
with $\sigma >0$. The constant $C>0$ depends on $\sigma $ but it does not depend on $\ell\in \mathbb N$.

We take $\sigma=\sigma (A)>0$ according to \cite[Lemma 3.1]{GCMST2} and for which we have that
$$
T_{k,n, {\rm loc}}^\nu (f)(x)=T_{k,n,{\rm loc}}^\nu(\mathcal{X}_{\sigma B_\ell}f)(x),\quad x\in B_\ell\mbox{ and }\ell\in \mathbb N. 
$$
Let $1<p<\infty$. According to the properties of the sequence $\{B_\ell\}_{\ell\in \mathbb N}$ (\cite[Lemma 3.1]{GCMST2}) we deduce that
\begin{align*}
    \int_{\mathbb R^d}\Big\|T_{k,n,{\rm loc }}^\nu (f)(x,\cdot)\Big\|_{L^2((0,\infty ),\frac{dt}{t})}^pdx&\leq C\sum_{\ell\in \mathbb N} \int_{B_\ell}\Big\|T_{k,n,{\rm loc}}^\nu(\mathcal{X}_{\sigma B_\ell}f)(x,\cdot)\Big\|_{L^2((0,\infty ),\frac{dt}{t})}^pdx\\
    &\leq C\sum_{\ell\in \mathbb N} e^{pR(Q^{1/2}w_\ell)}\int_{\mathbb R^d}\Big\|\mathbb T_{k,n,{\rm loc}}^\nu(e^{-R(Q^{1/2}\cdot)}\mathcal X_{\sigma B_\ell}f)(x,\cdot )\Big\|_{L^2((0,\infty ),\frac{dt}{t})}^pdx\\
    &\leq C\sum_{\ell\in \mathbb N} e^{pR(Q^{1/2}w_\ell)}\int_{\sigma B_\ell}|e^{-R(Q^{1/2}y)}f(y)|^pdy\\
    &\leq  C\sum_{\ell\in \mathbb N}\int_{\sigma B_\ell}|f(y)|^pdy\leq C\int_{\mathbb R^d}|f(y)|^pdy.
\end{align*}

In a similar way we can see that
$$
\big|\big\{x\in \mathbb R^d:\|T_{k,n,{\rm loc }}^\nu (f)(x,\cdot)\|_{L^2((0,\infty ),\frac{dt}{t})}>\lambda \big\}\big|\leq \frac{C}{\lambda}\|f\|_{L^1(\mathbb R^d,dx)},\quad \lambda >0.
$$

By all the estimations above, to obtain the $L^p$-boundedness properties for $S_{k,j,{\rm loc}}^{\nu ,2}$ it remains to establish that $\mathbb{T}_{k,n,{\rm loc}}^\nu$ is a bounded operator from $L^p(\mathbb R^d,dx)$ into $L^p_{L^2((0,\infty ),\frac{dt}{t})}(\mathbb R^d,dx)$, for every $1<p<\infty$, and from $L^1(\mathbb R^d,dx)$ into $L^{1,\infty}_{L^2((0,\infty ),\frac{dt}{t})}(\mathbb R^d,dx)$.

We observe that, for every $f\in L^p(\mathbb R ^d,dx)$, $1\leq p<\infty$,
$$
\mathbb{T}_{k,n}^\nu (f)(x,t)=t^{k+1}\partial _t^k\partial _{x_n}\mathbb P_t^\nu (f)(x),\quad x\in \mathbb R^d\mbox{ and }t>0,
$$
where, for every $t>0$,
$$
\mathbb{P}_t^\nu (f)(x)=\frac{t^{2\nu}}{4^\nu \Gamma (\nu )}\int _0^\infty e^{-\frac{t^2}{4u}}\mathbb W_u(f)(x)\frac{du}{u^{\nu +1}},\quad x\in \mathbb R^d.
$$
The derivations under the integral signs can be justified by proceeding as in Section \ref{section2}. 

We define the Littlewood-Paley function $\mathbb G_{k,n}^\nu$ as follows
$$
\mathbb G_{k,n}^\nu(f)(x)=\left(\int_0^\infty |t^{k+1}\partial _t^k\partial_{x_n}\mathbb{P}_t^\nu (f)(x)|^2\frac{dt}{t}\right)^{1/2},\quad x\in \mathbb R^d.
$$
Then, it is clear that 
$$
\mathbb G_{k,n}^\nu(f)(x)=\Big\|\mathbb T_{k,n}^\nu (f)(x,\cdot )\Big\|_{L^2((0,\infty ),\frac{dt}{t})},\quad x\in \mathbb R^d.
$$
We are going to see that the operator $\mathbb T_{k,n}^\nu$ can be seen as a $L^2((0,\infty ),\frac{dt}{t})$-valued Calder\'on-Zygmund operator. The integral kernel $\mathbb{T}_{k,n}^\nu (z,\cdot)$ is a standard $L^2((0,\infty ),\frac{dt}{t})$-valued Calder\'on-Zygmund kernel. Indeed, Minkowski inequality leads to 
\begin{align}\label{3.1}
    \Big\|\mathbb{T}_{k,n}^\nu (z,\cdot )\Big\|_{L^2((0,\infty ),\frac{dt}{t})}&\leq C\int_0^\infty \big\|[s^{k+1}\partial _s^k\mathfrak {g}_\nu (s)]_{\big|s=\frac{t}{2\sqrt{u}}}\big\|_{L^2((0,\infty ),\frac{dt}{t})}|z_n|e^{-\frac{|z|^2}{2u}}\frac{du}{u^{\frac{d+3}{2}}}\nonumber\\
    &\leq C|z|\int_0^\infty e^{-\frac{|z|^2}{2u}}\frac{du}{u^{\frac{d+3}{2}}}\leq \frac{C}{|z|^d},\quad z\not=0.
\end{align}
In a similar way we can see that
\begin{equation}\label{3.2}
    \sum_{\ell=1}^d\big\|\partial _{z_\ell}T_{k,n}^\nu (z,\cdot )\Big\|_{L^2((0,\infty ),\frac{dt}{t})}\leq \frac{C}{|z|^{d+1}},\quad z=(z_1,...,z_d)\in \mathbb{R}^d,\,z\not=0.
\end{equation}

The square function $\mathbb{G}_{k,n}^\nu$ is bounded from $L^2(\mathbb R^d,dx)$ into itself. In order to prove this property we use the Fourier transform $\mathcal{F}$ defined by 
\begin{equation}\label{Fouriertransform}
\mathcal{F}(f)(\xi)=\frac{1}{(2\pi )^{\frac{d}{2}}}\int_{\mathbb{R}^d}f(x)e^{-i\langle x,\xi\rangle}dx,\quad \xi\in \mathbb{R}^d.
\end{equation}
Let $f\in L^2(\mathbb R^d,dx)$. Plancherel equality allows us to get
\begin{align*}
  \big\|\mathbb{G}_{k,n}^\nu(f)\big\|_{L^2(\mathbb R^d,dx)}^2&=\int_0^\infty \int_{\mathbb R^d}|t^{k+1}\partial_t^k\partial _{x_n}\mathbb P_t^\nu (f)(x)|^2dx\frac{dt}{t}=\int_0^\infty \int_{\mathbb R^d}|t^{k+1}\partial _t^k\mathcal{F}(\partial _{x_n}\mathbb P_t^\nu (f))(y)|^2dy\frac{dt}{t}\\
  &\leq C\int_0^\infty \int_{\mathbb R^d}|y|^2|t^{k+1}\partial _t^k\mathcal{F}(\mathbb P_ t^\nu (f))(y)|^2dy\frac{dt}{t}.
\end{align*}
We can write
\begin{align*}
\mathcal{F}(\mathbb P_t^\nu (f))(y)&=\frac{t^{2\nu }}{4^\nu \Gamma (\nu )}\int_0^\infty e^{-\frac{t^2}{4u}}\mathcal{F}(\mathbb W_u(f))(y)\frac{du}{u^{\nu+1}}\\
&=\frac{\mathcal{F}(f)(y)}{4^\nu \Gamma (\nu )}\int_0^\infty t^{2\nu }e^{-\frac{t^2}{4u}}e^{-\frac{u|y|^2}{2}}\frac{du}{u^{\nu +1}}\\
&=\frac{\mathcal{F}(f)(y)}{\Gamma(\nu )} \int_0^\infty \mathfrak{g}_\nu (s)_{|s=\frac{t}{2\sqrt{u}}}e^{-\frac{u|y|^2}{2}}\frac{du}{u},\quad y\in \mathbb R^d\mbox{ and }t>0.
\end{align*}
Then,
$$
t^{k+1}\partial_t^k\mathcal{F}(\mathbb P_t^\nu (f))(y)=\frac{2}{\Gamma (\nu)}\mathcal{F}(f)(y)\int_0^\infty [s^{k+1}\partial _s^k\mathfrak{g}_\nu (s)]_{|s=\frac{t}{2\sqrt{u}}}e^{-\frac{u|y|^2}{2}}\frac{du}{\sqrt{u}}.
$$
By using Minkowski inequality we get
\begin{align}\label{3.3}
    \big\|\mathbb{G}_{k,n}^\nu(f)\big\|_{L^2(\mathbb R^d,dx)}^2&\leq C\int_{\mathbb R^d}|y|^2|\mathcal{F}(f)(y)|^2\Big\| \int_0^\infty [s^{k+1}\partial _s^k\mathfrak{g}_\nu (s)]_{|s=\frac{t}{2\sqrt{u}}}e^{-\frac{u|y|^2}{2}}\frac{du}{\sqrt{u}}\Big\|_{L^2(0,\infty ),\frac{dt}{t})}^2dy\nonumber\\
    &\leq C\int_{\mathbb R^d}|y|^2|\mathcal{F}(f)(y)|^2 \Big(\int_0^\infty \|[s^{k+1}\partial_s^k\mathfrak{g}_\nu (s)]_{|s=\frac{t}{2\sqrt{u}}}\|_{L^2(0,\infty ),\frac{dt}{t})}e^{-\frac{u|y|^2}{2}}\frac{du}{\sqrt{u}}\Big)^2dy\nonumber\\
    &\leq C\int_{\mathbb R^d}|y|^2|\mathcal{F}(f)(y)|^2\Big(\int_0^\infty e^{-\frac{u|y|^2}{2}}\frac{du}{\sqrt{u}}\Big)^2dy\nonumber\\
    &\leq C\int_{\mathbb R^d}|\mathcal{F}(f)(y)|^2dy= C\|f\|_{L^2(\mathbb R^d,dx)}^2.
\end{align}
Let $N\in \mathbb{N}$, $N\geq 2$. We define $E_N=L^2((\frac{1}{N},N),\frac{dt}{t})$. Assume that $f\in C_c^\infty (\mathbb R ^d)$. Let $x\in \mathbb R^d$. We consider the $E_N$-valued function $F_x$ defined on $\mathbb R ^d$ as follows
$$
[F_x(y)](t)=f(y)\mathbb{T}_{k,n}^\nu (x-y,t),\quad t\in \Big(\frac{1}{N},N\Big)\mbox{ and }y\in \mathbb R^d.
$$
The function $F_x$ is $E_N$-strongly measurable in $\mathbb R^d$ because $F_x$ is $E_N$-valued continuous in $\mathbb R ^d$ and $E_N$ is a separable Banach space. On the other hand, from \eqref{3.1} we deduce that
$$
\int_{\mathbb R^d}\|\mathbb{T}_{k,n}^\nu (x-y,\cdot )\|_{E_N}|f(y)|dy<\infty ,\quad x\not \in {\rm supp }f.
$$
We define
$$
\tau _{k,n}^\nu(f)(x)=\int_{\mathbb R^d}\mathbb{T}_{k,n}^\nu (x-y,\cdot )f(y)dy,\quad x\not\in {\rm supp }f,
$$
where the integral is understood in the $E_N$- Bochner sense.

By using the properties of the Bochner integral, for every $g\in E_N$, we have that 
$$
\int_{1/N}^N g(t)[\tau _{k,n}^\nu (f)(x)](t)\frac{dt}{t}=\int_{\mathbb R^d}f(y)\int_{1/N}^N\mathbb{T}_{k,n}^\nu (x-y,t)g(t)\frac{dt}{t}dy=\int_{1/N}^Ng(t)\mathbb T _{k,n}^\nu (f)(x,t)\frac{dt}{t},\quad x\not\in {\rm supp }f. 
$$
The interchange of the order of integration is justified because from \eqref{3.1} we deduce 
\begin{align*}
    \int_{\mathbb R^d}\int_{1/N}^N|f(y)||\mathbb{T}_{k,n}^\nu (x-y,t)||g(t)|\frac{dt}{t}dy&\leq C\int_{\mathbb R^d}|f(y)|\|g\|_{E_N}\|\mathbb{T}_{k,n}^\nu (x-y,\cdot )\|_{E_N}dy\\
    &\leq C\int_{{\rm supp }f}\frac{|f(y)|}{|x-y|^d}dy<\infty ,\quad x\not \in {\rm supp }f.
\end{align*}
It follows that, for every $x\not\in {\rm supp }f$,
$$
\tau _{k,n}^\nu (f)(x)=\mathbb T_{k,n}^\nu (f)(x,\cdot ),\quad \mbox{ in }E_N.
$$
According to the vector valued Calder\'on-Zygmund theory (\cite{RbFRT}), $\mathbb{T}_{k,n}^\nu $ can be extended from $L^p(\mathbb R^d,dx)\cap L^2(\mathbb R^d,dx)$ to $L^p(\mathbb R^d,dx)$ as a bounded operator $\mathscr{T}_{k,n}^\nu$ from $L^p(\mathbb R^d,dx)$ into $L^p_{E_N}(\mathbb R^d, dx)$, for every $1<p<\infty $, and from $L^1(\mathbb R^d,dx)$ into $L^{1,\infty }_{E_N}(\mathbb R^d,dx)$.

Furthermore, since \eqref{3.1}, \eqref{3.2} and \eqref{3.3} do not depend on $N$, we obtain, for every $1<p<\infty$,
$$
\sup_{N\in \mathbb{N}}\|\mathscr{T} _{k,n}^\nu \|_{L^p(\mathbb R ^d,dx)\rightarrow L^p_{E_N}(\mathbb R^d,dx)}<\infty ,
$$
and 
$$
\sup_{N\in \mathbb{N}}\|\mathscr{T} _{k,n}^\nu \|_{L^1(\mathbb R ^d,dx)\rightarrow L^{1,\infty}_{E_N}(\mathbb R^d,dx)}<\infty .
$$
Let $f\in L^p(\mathbb{R}^d,dx)$ with $1\leq p<\infty$. We have that
$$
\mathbb{T}_{k,n}^\nu (f)(x,t)\leq C\sup_{u>0}|\mathbb{W}_u(f)(x)|,\quad x\in \mathbb{R}^d \mbox{ and }t\in \Big(\frac{1}{N},N\Big).
$$
It follows that
$$
\|\mathbb{T}_{k,n}^\nu (f)(x,\cdot )\|_{E_N}\leq C\sup_{u>0}|\mathbb{W}_u(f)(x)|,\quad x\in \mathbb{R}^d.
$$
Here $C=C(N)$ is a positive constant depending on $N$.

Suppose that, for every $r\in \mathbb{N}$, $f_r\in L^p(\mathbb{R}^d,dx)\cap L^2(\mathbb{R}^d,dx)$ and that $f_r\longrightarrow f$, as $r\rightarrow \infty$, in $L^p(\mathbb{R}^d,dx)$. Let $1<p<\infty$. Since the maximal operator $\mathbb{W}_*$ associated with $\{\mathbb{W}_t\}_{t>0}$ is bounded from $L^p(\mathbb{R}^d,dx)$ into itself we deduce that
$$
\|\mathbb{T}_{k,n}^\nu (f)(x,\cdot )\|_{E_N}=\lim_{r\rightarrow \infty }\|\mathbb{T}_{k,n}^\nu(f_r)(x,\cdot )\|_{E_N},\quad \mbox{in }L^p(\mathbb{R}^d,dx).
$$
Hence $\|\mathbb{T}_{k,n}^\nu (f)(x,\cdot )\|_{E_N}=\|\mathscr{T}_{k,n}^\nu (f)(x,\cdot )\|_{E_N}$, a.e. $x\in \mathbb{R}^d$.

By taking into account that $\mathbb{W}_*$ is bounded from $L^1(\mathbb{R}^d,dx)$ into $L^{1,\infty }(\mathbb{R}^d,dx)$ as above we obtain that, for every $f\in L^1(\mathbb{R}^d, dx)$,
$$
\|\mathbb{T}_{k,n}^\nu (f)(x,\cdot )\|_{E_N}=\|\mathscr{T}_{k,n}^\nu (f)(x,\cdot )\|_{E_N},\quad \mbox{ a.e. }x\in \mathbb{R}^d.
$$
Then, we get that, for every $1<p<\infty$,
$$
\sup_{N\in \mathbb{N}}\|\mathbb{T}_{k,n}^\nu \|_{L^p(\mathbb{R}^d,dx)\rightarrow L^p_{E_N}(\mathbb{R}^d,dx)}<\infty,
$$
and
$$
\sup_{N\in \mathbb{N}}\|\mathbb{T}_{k,n}^\nu \|_{L^1(\mathbb{R}^d,dx)\rightarrow L_{E_N}^{1,\infty}(\mathbb{R}^d,dx)}<\infty.
$$

By using the monotone convergence theorem we conclude that the square function $\mathbb G_{k,n}^\nu $ is bounded from $L^p(\mathbb R^d,dx)$ into itself, for every $1<p<\infty $, and from $L^1(\mathbb R^d,dx)$ into $L^{1,\infty }(\mathbb R^d,dx)$, that is, $\mathbb{T}_{k,n}^\nu$ is bounded from $L^p(\mathbb R^d,dx)$ into $L^p_{L^2(0,\infty ),\frac{dt}{t})}(\mathbb R^d,dx)$, for every $1<p<\infty$, and from $L^1(\mathbb R^d,dx)$ into $L^{1,\infty}_{L^2((0,\infty ),\frac{dt}{t})}(\mathbb R^d,dx)$.

According to \cite[Proposition 3.4]{GCMST2} we deduce that the local operator $\mathbb T_{k,n,{\rm loc}}^\nu$ is also bounded from $L^p(\mathbb R^d,dx)$ into  $L^p_{L^2((0,\infty),\frac{dt}{t})}(\mathbb R^d,dx)$, for every $1<p<\infty$, and from $L^1(\mathbb R^d,dx)$ into $L^{1,\infty}_{L^2((0,\infty),\frac{dt}{t})}(\mathbb R^d,dx)$.

We conclude that $S_{k,j,{\rm loc}}^{\nu ,2}$ is a bounded operator from $L^p(\mathbb{R}^d,dx)$ into $L^p_{L^2((0,\infty),\frac{dt}{t})}(\mathbb{R}^d,dx)$, for every $1<p<\infty$, and from $L^1(\mathbb R^d,dx)$ into $L^{1,\infty}_{L^2((0,\infty),\frac{dt}{t})}(\mathbb R^d,dx)$. Then, the same $L^p$-boundedness properties are obtained for $S_{k,j,{\rm loc}}^\nu$ and by considering the comments after formula \eqref{equiv} the proof of (I1) can be finished. 

Let us show now property (I2). Minkowski inequality leads to
$$
   |\mathcal{D}_{k,j}^\nu (f)(x)|\leq  \int_{\mathbb{R}^d}\varphi_{\mbox{\tiny A}} (x,y)|f(y)|\left\|K_{k,j}^\nu (x,y,\cdot )-S_{k,j}^\nu (x,y,\cdot )\right\|_{L^2((0,\infty ),\frac{dt}{t})}e^{-R(y)}dy,\quad x\in \mathbb{R}.
$$
Our objective is to prove that
\begin{equation}\label{I2}
e^{-R(y)}\left\|K_{k,j}^\nu (x,y,\cdot )-S_{k,j}^\nu (x,y,\cdot )\right\|_{L^2((0,\infty ),\frac{dt}{t})}\leq C\frac{1+|x|}{|x-y|^{d-1}},\quad (x,y)\in L_{2A}.
\end{equation}
Thus, we obtain that $\mathcal{D}_{k,j}^\nu (f)\leq C\mathscr{S}_\eta(|f|)$ for $\eta =2A$, and, consequently, $\mathcal{D}_{k,j}^\nu$ is bounded from $L^p(\mathbb{R}^d,dx)$ into itself (and also from $L^p(\mathbb{R}^d,\gamma_\infty )$ into itself), for every $1\leq p\leq \infty$.

Again by Minkowski inequality, it follows that
 \begin{align*}
 e^{-R(y)}\left\|K_{k,j}^\nu (x,y,\cdot )-S_{k,j}^\nu (x,y,\cdot )\right\|_{L^2((0,\infty ),\frac{dt}{t})} &\\
&\hspace{-3cm}\leq Ce^{-R(y)}\int_0^\infty \|[s^{k+1}\partial _s^k{\mathfrak g}_\nu(s)]_{|s=\frac{t}{2\sqrt{u}}}\|_{L^2((0,\infty ),\frac{dt}{t})}|\partial _{x_j}(h_u(x,y))- \mathbb{S}_u^j(x,y)|\frac{du}{\sqrt{u}}\\
&\hspace{-3cm}\leq Ce^{-R(y)}\int_0^\infty |\partial _{x_j}(h_u(x,y))- \mathbb{S}_u^j(x,y)|\frac{du}{\sqrt{u}}\\
&\hspace{-3cm}=Ce^{-R(y)}\left(\int_0^{\mathfrak{m}(x)}+\int_{\mathfrak{m}(x)}^\infty \right)|\partial _{x_j}(h_u(x,y))- \mathbb{S}_u^j(x,y)|\frac{du}{\sqrt{u}}\\
&\hspace{-3cm}=I_j^0(x,y)+I_j^\infty (x,y),\quad x,y\in \mathbb{R}^d.
\end{align*}

We first estimate $I_j^\infty (x,y)$, $(x,y)\in L_{2A}$. Here the cancellation does not play now any role.

By using \eqref{(2.3)} and \eqref{(2.4)} we get
\begin{equation}\label{partialm1}
|\partial _{x_j}h_u(x,y)|\leq C\frac{e^{R(x)}}{u^{\frac{d}{2}}}(1+|x|),\quad x,y\in \mathbb{R}^d,\,\mathfrak{m}(x)\leq u\leq 1,
\end{equation}
and
\begin{equation}\label{partialinfty}
|\partial _{x_j}h_u(x,y)|\leq Ce^{R(x)}e^{-cu}(1+|y|),\quad x,y\in \mathbb{R}^d,\, u>1.
\end{equation}

On the other hand,
$$
|\mathbb{S}_u^j(x,y)|\leq C\frac{e^{R(x)}}{u^{\frac{d}{2}}}\frac{|y-x|}{u}e^{-c\frac{|y-x|^2}{u}}\leq C\frac{e^{R(x)}}{u^{\frac{d+1}{2}}}\leq C\frac{e^{R(x)}}{u^{\frac{d}{2}}}(1+|x|),\quad x,y\in \mathbb{R}^d \mbox{ and }\mathfrak{m}(x)\leq u.
$$ 

Thus, by using also \eqref{3.2'} with $\eta=2A$ and taking into account that if $(x,y)\in L_{2A}$, then $|x-y|\leq \frac{2A}{1+|x|}\leq 2A\sqrt{\mathfrak{m}(x)}\leq 2A$ and $|y|\leq C(1+|x|)$, we obtain
\begin{align}\label{Iinfty}
    I_j^\infty (x,y)&\leq Ce^{R(x)-R(y)}(1+|x|)\left(\int_{\mathfrak{m}(x)}^1\frac{du}{u^{\frac{d+1}{2}}}+\int_1^\infty \big(e^{-cu}+\frac{1}{u^{\frac{d}{2}}}\big)\frac{du}{\sqrt{u}}\right)\nonumber\\
    &\leq C(1+|x|)\left(\int_{\mathfrak{m}(x)}^\infty u^{-\frac{d+1}{2}}du+1\right)\leq C(1+|x|)\Big(\frac{1}{\mathfrak{m}(x)^{\frac{d-1}{2}}}+1\Big)\leq C\frac{1+|x|}{|x-y|^{d-1}},\quad (x,y)\in L_{2A}.
\end{align}

Next we analyze the term $I_j^0(x,y)$, $(x,y)\in L_{2A}$. Now the cancellation property of the difference is relevant to get the bound we need. We consider 
\begin{equation}\label{Htilde}
\widetilde{H}_u(x,y)=\Big(\frac{{\rm det }\,Q_\infty }{u^d{\rm det}\,Q}\Big)^{1/2}e^{R(x)}e^{ -\frac{1}{2}\langle(Q_t^{-1}-Q_\infty ^{-1})(y-D_ux),y-D_ux\rangle},\quad x,y\in \mathbb R^d\mbox{ and }u>0, 
\end{equation}
and
\begin{equation}\label{Wtilde}
\widetilde{W}_u(x,y)=\Big(\frac{{\rm det }\,Q_\infty }{u^d{\rm det}\,Q}\Big)^{1/2}e^{R(x)}e^{-\frac{1}{2u}\langle Q^{-1}(y-x),y-x\rangle},\quad x,y\in \mathbb R^d\mbox{ and }u>0. 
\end{equation}

We note that $h_u(x,y)=(u^d{\rm det}\,Q)^{1/2}({\rm det}\,Q_u)^{-1/2}\widetilde{H}_u(x,y)$, $x,y\in\mathbb{R}^d$ and $u>0$ and that 
\begin{equation}\label{Wtildederivada}
\partial _{x_j} \widetilde{W}_u(x,y)=\langle Q_\infty ^{-1}x,e_j\rangle \widetilde{W}_u(x,y)+\mathbb{S}_u^j(x,y),\quad x,y\in \mathbb R^d\mbox{ and }u>0.
\end{equation}
Hence, we can write
\begin{align}\label{I0}
I_j^0(x,y)&\leq Ce^{-R(y)}\int_0^{\mathfrak{m}(x)}\Big(|\partial _{x_j}[h_u(x,y)-\widetilde{W}_u(x,y)]|+|\langle Q_\infty^{-1}x,e_j\rangle| \widetilde{W}_u(x,y)\Big)\frac{du}{\sqrt{u}}\nonumber\\
&\leq Ce^{-R(y)}\left(\int_0^{\mathfrak{m}(x)}\Big(\big|(u^d{\rm det}\,Q)^{-1/2}({\rm det}\,Q_u)^{1/2}-1\big||\partial_{x_j}h_u(x,y)|\Big)\frac{du}{\sqrt{u}}\nonumber\right.\\
&\left.\quad +\int_0^{\mathfrak{m}(x)}\big|\partial_{x_j}\big[\widetilde{H}_u(x,y)-\widetilde{W}_u(x,y)\big]\big|\frac{du}{\sqrt{u}}+\int_0^{\mathfrak{m}(x)}|\langle Q_\infty^{-1}x,e_j\rangle| \widetilde{W}_u(x,y)\frac{du}{\sqrt{u}}\right)\nonumber\\
&=I_j^{0,1}(x,y)+I_j^{0,2}(x,y)+I_j^{0,3}(x,y),\quad x,y\in \mathbb{R}^d.
\end{align}

From \eqref{3.2'} and since $Q$ is symmetric and positive definite we deduce that
\begin{equation}\label{I03}
I_j^{0,3}(x,y)\leq C(1+|x|)\int_0^{\mathfrak{m}(x)}\frac{e^{-c\frac{|x-y|^2}{u}}}{u^{\frac{d+1}{2}}}du\leq C\frac{1+|x|}{|x-y|^{d-1}},\quad (x,y)\in L_{2A}.
\end{equation}

Next, we analyze $I_j^{0,1}(x,y)$, $(x,y)\in L_{2A}$. We recall that  $Q_t=\int_0^te^{Bs}Qe^{B^*s}ds$, $t>0$. We can write
$$
\frac{1}{t}Q_t-Q=\frac{1}{t}\int_0^t(e^{Bs}-I)Qe^{B^*s}ds+\frac{1}{t}\int_0^t Q(e^{B^*s}-I)ds,\quad t>0.
$$
Then, since $|e^{Bs}-I|\leq Cs$, $s\in (0,1)$, it follows that $\lim_{t\rightarrow 0^+}\frac{1}{t}Q_t=Q$. Hence, 
$$
\lim_{t\rightarrow 0^+}{\rm det }\big(\frac{1}{t}Q_t\big)=\lim_{t\rightarrow 0^+}\frac{1}{t^d}{\rm det }(Q_t)={\rm det }(Q).
$$
Furthermore, we have that $|\frac{1}{t}Q_t-Q|\leq Ct$, $t\in (0,1)$. We write $Q_t=(q_{ij}^{(t)})_{i,j=1}^d$ and $Q=(q_{ij})_{i,j=1}^d$. Then, $\frac{1}{t}q_{ij}^{(t)}=q_{ij}+\phi_{ij}(t)$, $t\in (0,1)$, where $|\phi _{ij}(t)|\leq Ct$, $t\in (0,1)$. We get $\frac{1}{t^d}{\rm det }(Q_t)={\rm det }(Q)+\phi(t)$, $t\in (0,1)$, being $|\phi (t)|\leq Ct$, $t\in (0,1)$. According to \cite[Lemma 2.2]{CCS3} we obtain
\begin{align}\label{3.0}
    \left|\frac{({\rm det }\,Q_u)^{1/2}}{(u^d{\rm det }\,Q)^{1/2}}-1\right|&=\frac{|({\rm det }\,Q_u)^{1/2}-u^{\frac{d}{2}}({\rm det }\,Q)^{1/2}|}{(u^d{\rm det }\,Q)^{1/2}}\nonumber\\
    &=\frac{|{\rm det }\,Q_u-u^d{\rm det }\,Q|}{(u^d{\rm det }\,Q)^{1/2}[({\rm det }\,Q_u)^{1/2}+u^{\frac{d}{2}}({\rm det }\,Q)^{1/2})}\leq C|\phi (u)|\leq Cu,\quad u\in (0,1).
\end{align}

By \cite[Lemma 2.3]{CCS3} we get
\begin{equation}\label{3.1'}
|y-D_tx|\geq |x-y|-|x-D_tx|\geq |y-x|-Ct|x|\geq |y-x|-C,\quad x,y\in \mathbb R^d\mbox{ and }0<t<\mathfrak{m}(x).
\end{equation}

Then, by taking into account also \eqref{(2.3)} and \eqref{3.2'} it follows that
\begin{align}\label{I01}
   I_j^{0,1}(x,y)&\leq Ce^{R(x)-R(y)}\int_0^{\mathfrak{m}(x)}e^{-c\frac{|y-x|^2}{u}}\Big(|x|+\frac{1}{\sqrt{u}}\Big)\frac{du}{u^{\frac{d-1}{2}}}\nonumber\\
   &\leq C(1+|x|)\int_0^{\mathfrak{m}(x)} e^{-c\frac{|y-x|^2}{u}}\frac{du}{u^{\frac{d}{2}}}\leq  C(1+|x|)\int_0^{\mathfrak{m}(x)} e^{-c\frac{|y-x|^2}{u}}\frac{du}{u^{\frac{d+1}{2}}}\nonumber\\
&\leq C(1+|x|)\int_0^\infty e^{-c\frac{|y-x|^2}{u}}\frac{du}{u^{\frac{d+1}{2}}}\leq C\frac{1+|x|}{|x-y|^{d-1}},\quad (x,y)\in L_{2A}.
\end{align}

To deal with $I_j^{0,2}(x,y)$, $(x,y)\in L_{2A}$, we write, according to \eqref{(2.1)}, 
$$
\partial _{x_j}\widetilde{H}_u(x,y)=\big(\langle Q_\infty ^{-1}x,e_j\rangle+\langle Q_u^{-1}e^{uB}e_j, y-D_ux\rangle\big)\widetilde{H}_u(x,y),\quad x,y\in \mathbb R^d\mbox{ and }u>0,
$$
and consider \eqref{Wtildederivada} to obtain
\begin{align}\label{diferenciaHW}
\big|\partial _{x_j}\big[\widetilde{H}_u(x,y)-\widetilde{W}_u(x,y)\big]\big|&\leq \Big|\langle Q_\infty ^{-1}x,e_j\rangle-\langle \frac{1}{u}Q^{-1}e_j, y-x\rangle\Big|\Big|\widetilde{H}_u(x,y)-\widetilde{W}_u(x,y)\Big|\nonumber\\
&\quad +\Big|\langle Q_u^{-1}e^{uB}e_j, y-D_ux\rangle-\langle \frac{1}{u}Q^{-1}e_j, y-x\rangle\Big|\widetilde{H}_u(x,y)\nonumber\\
&\leq\Big|\langle Q_\infty ^{-1}x,e_j\rangle-\langle \frac{1}{u}Q^{-1}e_j, y-x\rangle\Big|\Big|\widetilde{H}_u(x,y)-\widetilde{W}_u(x,y)\Big|\\
&\quad +\Big(|\langle Q_u^{-1}e^{uB}e_j-\frac{1}{u}Q^{-1}e_j, y-D_ux\rangle|\nonumber\\
&\quad +|\langle \frac{1}{u}Q^{-1}e_j, x-D_ux\rangle|\Big)\widetilde{H}_u(x,y),\quad x,y\in \mathbb R^d\mbox{ and }u>0.\nonumber
\end{align}

We now establish some useful estimates.

{\bf (i)} We can write
\begin{align*}
Q_t^{-1}e^{tB}e_j-\frac{1}{t}Q^{-1}e_j&=Q_t^{-1}e^{tB}e_j-\frac{1}{t}Q^{-1}e^{tB}e_j+\frac{1}{t}Q^{-1}(e^{tB}-I)e_j\\
&=Q_t^{-1}(tQ-Q_t)\frac{1}{t}Q^{-1}e^{tB}e_j+\frac{1}{t}Q^{-1}(e^{tB}-I)e_j,\quad t>0,
\end{align*}
and
$$
e^{tB^*}Q_t^{-1}e^{tB}e_j-\frac{1}{t}Q^{-1}e_j=(e^{tB^*}-I)Q_t^{-1}e^{tB}e_j-(Q_t^{-1}e^{tB}e_j-\frac{1}{t}Q^{-1}e_j),\quad t>0.
$$
By using \cite[Lemma 2.2]{CCS3}, since $|tQ-Q_t|\leq Ct^2$, $t\in (0, 1)$, we get
\begin{equation}\label{3.3'}
|Q_t^{-1}e^{tB}e_j-\frac{1}{t}Q^{-1}e_j|+|e^{tB^*}Q_t^{-1}e^{tB}e_j-\frac{1}{t}Q^{-1}e_j|\leq C,\quad t\in (0,1).
\end{equation}

{\bf (ii)} By \cite[Lemma 2.3]{CCS3} we obtain
\begin{equation}\label{3.4'}
    |y-D_tx|\leq |y-x|+Ct|x|,\quad x,y\in \mathbb R^d\mbox{ and }\;t\in (0,1).
\end{equation}

{\bf (iii)} For every $t>0$, since $Q_t$ is a symmetric positive definite matrix, $|Q_t^{1/2}|=|Q_t|^{1/2}$. We deduce that
$$
|z|=|Q_t^{1/2}Q_t^{-1/2}z|\leq |Q_t|^{1/2}|Q_t^{-1/2}z|,\quad z\in \mathbb R^d\mbox{ and }t>0.
$$
Furthermore,
$$
|Q_t|\leq |Q|\int_0^te^{s(|B|+|B^*|)}ds\leq Ct,\quad t\in (0,1).
$$
Then, we get 
\begin{equation}\label{3.4.1}
    |Q_t^{-1/2}z|\geq C\frac{|z|}{\sqrt{t}},\quad z\in  \mathbb R^d\mbox{ and }t\in (0,1).
\end{equation}
It follows that
\begin{equation}\label{3.5}
    e^{-\frac{1}{2}\langle Q_t^{-1}(y-D_tx),y-D_tx\rangle}=e^{-\frac{1}{2}|Q_t^{-1/2}(y-D_tx)|^2}\leq Ce^{-c\frac{|y-D_tx|^2}{t}},\quad x,y\in \mathbb R^d\mbox{ and }t\in (0,1).
\end{equation}
According to \cite[(2.10)]{CCS3} we obtain
\begin{equation}\label{3.6}
    e^{-\frac{1}{2}\langle(Q_t^{-1}-Q_\infty ^{-1})(y-D_tx),y-D_tx\rangle}\leq Ce^{-c\frac{|y-D_tx|^2}{t}},\quad x,y\in \mathbb R^d\mbox{ and }t\in (0,1).
\end{equation}
We can write
$$
e^{-\frac{1}{2}\langle Q_t^{-1}(y-x),y-D_tx\rangle}=e^{-\frac{1}{2}\langle Q_t^{-1}(y-x),y-x\rangle}e^{-\frac{1}{2}\langle Q_t^{-1}(y-x),x-D_tx\rangle},\quad x,y\in \mathbb R^d\mbox{ and }t>0.
$$
By \cite[Lemmas 2.2, (ii), and 2.3]{CCS3} we have that
\begin{align*}
    -\langle Q_t^{-1}(y-x),x-D_tx\rangle&\leq |\langle Q_t^{-1}(y-x),x-D_tx\rangle|\leq |Q_t^{-1}(y-x)||x-D_tx|\\
    &\leq C|y-x||x|\leq C\frac{|x|}{1+|x|}\leq C,\quad (x,y)\in L_{2A}\mbox{ and }t\in (0,1).
\end{align*}
Then, by \eqref{3.4.1}
\begin{align}\label{3.7}
e^{-\frac{1}{2}\langle(Q_t^{-1}(y-x),y-D_tx\rangle}&\leq Ce^{-\frac{1}{2}\langle(Q_t^{-1}(y-x),y-x\rangle}\leq Ce^{-\frac{1}{2}|Q_t^{-1/2}(y-x)|^2}\nonumber\\
&\leq Ce^{-c\frac{|x-y|^2}{t}},\quad (x,y)\in L_{2A}\mbox{ and }t\in (0,1).
\end{align}

{\bf (iv)} In the following three estimates we use again \cite[Lemmas 2.2 and 2.3]{CCS3}. It follows that
\begin{align}\label{3.8}
    |\langle(Q_t^{-1}-Q_\infty ^{-1})(y-D_tx),y-D_tx\rangle-\langle Q_t^{-1}(y-D_tx),y-D_tx\rangle|&\nonumber\\
    &\hspace{-4cm}\leq C|y-D_tx|^2\leq C(|y-x|^2+t^2|x|^2),\quad x,y\in \mathbb R^d\mbox{ and }t\in (0,1);
\end{align}
\begin{align}\label{3.9}
    |\langle Q_t^{-1}(y-D_tx),y-D_tx\rangle-\langle Q_t^{-1}(y-x),y-D_tx\rangle&\leq |\langle Q_t^{-1}(x-D_tx),y-D_tx\rangle|\nonumber \\
    &\hspace{-4cm}\leq  |Q_t^{-1}(x-D_tx)||y-D_tx|\leq C|x|(|y-x|+t|x|),\quad x,y\in \mathbb R^d\mbox{ and }t\in (0,1);
\end{align}
\begin{align}\label{3.10}
    |\langle Q_t^{-1}(y-x),y-D_tx\rangle-\langle Q_t^{-1}(y-x),y-x\rangle |&\leq |Q_t^{-1}(y-x)||x-D_tx|\nonumber \\
    &\leq  C|y-x||x|,\quad x,y\in \mathbb R^d\mbox{ and }t\in (0,1);
\end{align}
and
\begin{align}\label{3.11}
    |\langle Q_t^{-1}(y-x),y-x\rangle-\frac{1}{t}\langle Q^{-1}(y-x),y-x\rangle |&\leq | (Q_t^{-1}-\frac{1}{t}Q^{-1})(y-x)||y-x|\nonumber \\
    &\leq C|y-x|^2,\quad x,y\in \mathbb R^d\mbox{ and }t\in (0,1).
\end{align}
{\bf (v)} We have that, for every $x,y\in \mathbb R^d$ and $t>0$,
\begin{align*}
    e^{-\frac{1}{2}\langle(Q_t^{-1}-Q_\infty ^{-1})(y-D_tx),y-D_tx\rangle}-e^{-\frac{1}{2t}\langle Q^{-1}(y-x),y-x\rangle}&\\
    &\hspace{-8cm}=\Big(e^{-\frac{1}{2}\langle(Q_t^{-1}-Q_\infty ^{-1})(y-D_tx),y-D_tx\rangle}-e^{-\frac{1}{2}\langle Q_t^{-1}(y-D_tx),y-D_tx\rangle}\Big) +\Big(e^{-\frac{1}{2}\langle Q_t^{-1}(y-D_tx),y-D_tx\rangle}-e^{-\frac{1}{2}\langle Q_t^{-1}(y-x),y-D_tx\rangle}\Big)\\
    &\hspace{-8cm}\quad +\Big(e^{-\frac{1}{2}\langle Q_t^{-1}(y-x),y-D_tx\rangle}-e^{-\frac{1}{2}\langle Q_t^{-1}(y-x),y-x\rangle}\Big) +\Big(e^{-\frac{1}{2}\langle Q_t^{-1}(y-x),y-x\rangle}-e^{-\frac{1}{2t}\langle Q^{-1}(y-x),y-x\rangle}\Big).
\end{align*}
We now use the elementary inequality
$$
|e^{-a}-e^{-b}|\leq |b-a|e^{-\min\{a,b\}},\quad a,b\in \mathbb R.
$$
By estimations \eqref{3.5}-\eqref{3.11} we get
\begin{align*}
    \Big| e^{-\frac{1}{2}\langle(Q_t^{-1}-Q_\infty ^{-1})(y-D_tx),y-D_tx\rangle}-e^{-\frac{1}{2t}\langle Q^{-1}(y-x),y-x\rangle}\Big|&\\
    &\hspace{-6cm}\leq C\Big[(|y-x|^2+t^2|x|^2)e^{-c\frac{|y-D_tx|^2}{t}} +|x|(|y-x|+t|x|)\Big(e^{-c\frac{|y-D_tx|^2}{t}}+e^{-c\frac{|y-x|^2}{t}}\Big)\\
    &\hspace{-6cm} \quad +(|y-x||x|+|y-x|^2)e^{-c\frac{|y-x|^2}{t}}\Big],\quad (x,y)\in L_{2A}\mbox{ and }t\in (0,1).
\end{align*}
By \eqref{3.1'} we deduce that
\begin{align}\label{3.7.1}
   \Big| e^{-\frac{1}{2}\langle(Q_t^{-1}-Q_\infty ^{-1})(y-D_tx),y-D_tx\rangle}-e^{-\frac{1}{2t}\langle Q^{-1}(y-x),y-x\rangle}\Big|&\leq Ce^{-c\frac{|y-x|^2}{t}}(|y-x|^2+|y-x||x|+t|x|^2)  \nonumber\\
   &\hspace{-6cm}\leq Ce^{-c\frac{|y-x|^2}{t}}(|y-x|^2+|y-x||x|+t|x|^2)\leq Ce^{-c\frac{|y-x|^2}{t}}(t+\sqrt{t}|x|+t|x|^2)\nonumber\\
&\hspace{-6cm}\leq Ce^{-c\frac{|y-x|^2}{t}}\sqrt{t}(1+|x|),\quad (x,y)\in L_{2A}\mbox{ and }0<t<\mathfrak{m}(x).
\end{align}
In the last inequality we have used that $\mathfrak{m}(x)\sim (1+|x|^2)^{-1}$, $x\in \mathbb{R}^d$. 

Hence, by using \eqref{3.2'} and \eqref{3.7.1} and that $\sqrt{u}|x|\leq C$, when $x\in \mathbb{R}^d$ and $0<u<\mathfrak{m}(x)$, we obtain that
\begin{align*}
e^{-R(y)}\Big|\langle Q_\infty ^{-1}x,e_j\rangle-\langle \frac{1}{u}Q^{-1}e_j, y-x\rangle\Big|\Big|\widetilde{H}_u(x,y)-\widetilde{W}_u(x,y)\Big|&\leq C\Big(|x|+\frac{|y-x|}{u}\Big)\frac{e^{R(x)-R(y)}}{u^{\frac{d}{2}}}e^{-c\frac{|y-x|^2}{u}}\sqrt{u}(1+|x|)\\
&\hspace{-6cm}\leq C(\sqrt{u}|x|+1)\frac{e^{-c\frac{|y-x|^2}{u}}}{u^{\frac{d}{2}}}(1+|x|)\leq C\frac{e^{-c\frac{|y-x|^2}{u}}}{u^{\frac{d}{2}}}(1+|x|),\quad (x,y)\in L_{2A}\mbox{ and }0<u<\mathfrak{m}(x).
\end{align*}

On the other hand, from \eqref{3.3'}, \eqref{3.6}, \cite[Lemma 2.3]{CCS3}, and again \eqref{3.2'} and \eqref{3.1'}, it follows that
\begin{align*}
e^{-R(y)}\Big(|\langle Q_u^{-1}e^{uB}e_j-\frac{1}{u}Q^{-1}e_j, y-D_ux\rangle|+|\langle \frac{1}{u}Q^{-1}e_j, x-D_ux\rangle|\Big)\widetilde{H}_u(x,y)&\leq Ce^{-R(y)}(|y-D_ux|+|x|)\widetilde{H}_u(x,y)\\
&\hspace{-10cm}\leq C(|y-D_ux|+|x|)\frac{e^{R(x)-R(y)}}{u^{\frac{d}{2}}}e^{-c\frac{|y-D_ux|^2}{u}}\leq C(1+|x|)\frac{e^{-c\frac{|y-x|^2}{u}}}{u^{\frac{d}{2}}},\quad (x,y)\in L_{2A}\mbox{ and }0<u<\mathfrak{m}(x).
\end{align*}

These estimations, jointly \eqref{diferenciaHW}, lead to
$$
e^{-R(y)}\big|\partial _{x_j}\big[\widetilde{H}_u(x,y)-\widetilde{W}_u(x,y)\big]\big|\leq C(1+|x|)\frac{e^{-c\frac{|y-x|^2}{u}}}{u^{\frac{d}{2}}},\quad (x,y)\in L_{2A}\mbox{ and }0<u<\mathfrak{m}(x),
$$
and we conclude that
\begin{equation}\label{I02}
I_j^{0,2}(x,y)\leq C(1+|x|)\int_0^{\mathfrak{m}(x)}\frac{e^{-c\frac{|x-y|^2}{u}}}{u^{\frac{d+1}{2}}}du\leq C\frac{1+|x|}{|x-y|^{d-1}},\quad (x,y)\in L_{2A}.
\end{equation}

By considering \eqref{Iinfty}, \eqref{I0}, \eqref{I03}, \eqref{I01} and\eqref{I02} we get the estimation \eqref{I2} and property (I2) is established.

\subsubsection{About the local part of $g_{k,\alpha}^\nu$: $\widehat{\alpha}=2$}\label{S3.1.2}

Assume that $i,j\in\{1,\ldots,d\}$ and $\alpha=(\alpha_1,\ldots,\alpha_d)$ being $\alpha_i=\alpha_j=1$ and $\alpha_\ell=0$, $\ell\in\{1,\ldots,d\}\setminus\{i,j\}$, when $i\neq j$, and $\alpha_i=2$ and $\alpha_\ell=0$, $\ell\in\{1,\ldots,d\}\setminus\{i\}$, when $i=j$.

In this case we have 
$$
g_{k,\alpha,{\rm loc}}^\nu (f)(x)=\left\|\int_{\mathbb R^d}K_{k,i,j}^\nu(x,y,\cdot )\varphi_{\mbox{\tiny A}} (x,y)f(y)d\gamma_\infty (y)\right\|_{L^2((0,\infty ),\frac{dt}{t})},\quad x\in \mathbb{R}^d,
$$
where 
$$
K_{k,i,j}^\nu(x,y,t)=\frac{4}{\Gamma (\nu)}\int_0^\infty [s^{k+2}\partial _s^k{\mathfrak g}_\nu(s)]_{|s=\frac{t}{2\sqrt{u}}}\partial^2 _{x_ix_j}h_u(x,y)du,\quad x,y\in \mathbb{R}^d \mbox{ and }t>0.
$$

We consider the operator $S_{k,i,j}^\nu$ defined by
$$
S_{k,i,j}^\nu(f)(x,t)=\int_{\mathbb R^d}S_{k,i,j}^\nu(x,y,t)f(y)d\gamma_\infty (y)\quad x\in\mathbb R^d\mbox{ and }t>0,
$$
where
$$
S_{k,i,j}^\nu(x,y,t)=\frac{4}{\Gamma (\nu)}\int_0^\infty [s^{k+2}\partial _s^k{\mathfrak g}_\nu(s)]_{|s=\frac{t}{2\sqrt{u}}}\mathbb{S}_u^{i,j}(x,y)du,\quad x,y\in\mathbb R^d\mbox{ and }t>0,
$$
and
\begin{equation}\label{Su}
\mathbb{S}_u^{i,j}(x,y)=\left(\frac{{\rm det}\,Q_\infty}{u^d{\rm det}\,Q}\right)^{1/2}e^{R(x)}\partial^2_{x_ix_j}\big[ e^{-\frac{1}{2u}\langle Q^{-1}(y-x),y-x\rangle }\big],\quad x,y\in \mathbb{R}^d\mbox{ and }u>0.
\end{equation}

We also introduce the operator $\widetilde S_{k,i,j}^\nu$ by
$$
\widetilde S_{k,i,j}^\nu(f)(x,t)=\int_{\mathbb R^d}\widetilde S_{k,i,j}^\nu(x,y,t)f(y)dy,\quad x\in\mathbb R^d\mbox{ and }t>0,
$$
being
$$
\widetilde S_{k,i,j}^\nu(x,y,t)=e^{-R(Q^{1/2}y)}S_{k,i,j}^\nu(Q^{1/2}x,Q^{1/2}y,t),\quad x,y\in\mathbb R^d\;\mbox{and}\;t>0.$$
We define, for every $\ell,m=1,\ldots,d$, the operators
$$
T_{k,\ell,m}^\nu(f)(x,t)=\int_{\mathbb R^d}e^{R(Q^{1/2}x)-R(Q^{1/2}y)}\mathbb{T}_{k,\ell,m}^\nu(x-y,t)f(y) dy,\quad x\in\mathbb R^d\mbox{ and }t>0,
$$
and
$$
\mathbb T_{k,\ell,m}^\nu(f)(x,t)=\int_{\mathbb R^d}\mathbb{T}_{k,\ell,m}^\nu(x-y,t)f(y)dy,\quad x\in\mathbb R^d\mbox{ and }t>0,
$$
where
$$
\mathbb{T}_{k,\ell,m}^\nu(z,t)=\frac{4}{\Gamma (\nu)}\int_0^\infty  [s^{k+2}\partial_s^k(\mathfrak{g}_\nu(s)]_{|s=\frac{t}{2\sqrt{u}}}\partial^2_{z_\ell z_m}\mathbb{W}_u(z)du,\quad z=(z_1,\ldots,z_d)\in\mathbb R^d\mbox{ and }t>0.
$$
Let $\ell,m\in\{1,\ldots,d\}$. For every $f\in C_c^\infty(\mathbb R^d)$ we have that
$$
\mathbb T_{k,\ell,m}^\nu(f)(x,t)=t^{k+2}\partial^k_t\partial^2_{x_\ell x_m}\mathbb P_t^\nu(f)(x),\quad x\in\mathbb R^d\mbox{ and }t>0.
$$
We consider the square function $\mathbb G_{k,\ell,m}^\nu$ defined by
$$
\mathbb G_{k,\ell,m}^\nu(f)(x)=\left(\int_0^\infty|\mathbb T_{k,\ell,m}^\nu(f)(x,t)|^2\frac{dt}{t}\right)^{1/2},\quad x\in\mathbb R^d.
$$
We are going to see the $L^p$-boundedness properties of $\mathbb G_{k,\ell,m}^\nu$ by using vector valued Calder\'on-Zygmund theory. The arguments are similar to the ones used for the studying of $\mathbb{G}_{k,n}^\nu$ in the Section \ref{S3.1.1}. We sketch the proof in this case.

The operator $\mathbb T_{k,\ell,m}^\nu$ can be seen as a $L^2((0,\infty),\frac{dt}{t})$-singular integral having as integral kernel $t\rightarrow \mathbb{T}_{k,\ell,m}^\nu(x-y,t)$, $t\in (0,\infty)$, for every $x,y\in\mathbb R^d$.  Minkowski inequality leads to
\begin{align*}
\big\|\mathbb{T}_{k,\ell,m}^\nu(z,\cdot)\big\|_{L^2((0,\infty),\frac{dt}{t})} & \leq C \int_0^\infty \big\|[s^{k+2}\partial_s^k(\mathfrak{g}_\nu(s)]_{s=\frac{t}{2\sqrt{u}}}\big\|_{L^2((0,\infty),\frac{dt}{t})}\Big|\partial^2_{z_\ell z_m}\big(e^{-\frac{|z|^2}{2u}}\big)\Big|\frac{du}{u^{\frac{d}{2}}} \\
& \leq C \int_0^\infty \frac{|z|^2}{u^{\frac{d}{2}+2}}e^{-\frac{|z|^2}{4u}}du\leq C\int_0^\infty \frac{e^{-c\frac{|z|^2}{4u}}}{u^{\frac{d}{2}+1}}du\leq \frac{C}{|z|^d},\quad z\in\mathbb R^d,\,z\neq 0.
\end{align*}
In a similar way we can see that,
$$
\sum_{i=1}^d\|\partial_{z_i}\mathbb{T}_{k,\ell,m}^\nu(z,\cdot)\|_{L^2((0,\infty),\frac{dt}{t})}\leq \frac{C}{|z|^{d+1}},\quad z\in\mathbb R^d,\,z\neq 0.
$$
By using Fourier transform we can see that the square function $\mathbb G_{k,\ell,m}^\nu$ is bounded from $L^2(\mathbb R^d,dx)$ into itself.

Arguing as in the Section \ref{S3.1.1} for $\mathbb{G}_{k,n}^\nu$ and by using Calder\'on-Zygmund theory for vector valued singular integrals, we can conclude that $\mathbb G_{k,\ell,m}^\nu$ is a bounded operator from $L^p(\mathbb R^d,dx)$ into itself, for every $1<p<\infty$, and from $L^1(\mathbb R^d,dx)$ into $L^{1,\infty}(\mathbb R^d,dx)$. According to \cite[Proposition 3.4]{GCMST2} the local operator $\mathbb T^\nu_{k,\ell,m,{\rm loc}}$ defined in the usual way is a bounded operator from $L^p(\mathbb R^d,dx)$ into  $L^p_{L^2((0,\infty),\frac{dt}{t})}(\mathbb R^d,dx)$, when $1<p<\infty$, and from $L^1(\mathbb R^d,dx)$ into $L^{1,\infty}_{L^2((0,\infty),\frac{dt}{t})}(\mathbb R^d,dx)$. 

We can now follow the proof of $(I1)$ in Section \ref{S3.1.1} to establish that $S_{k,i,j,{\rm loc}}^\nu$ is a bounded operator $L^p(\mathbb R^d,\gamma_\infty)$ into  $L^p_{L^2((0,\infty),\frac{dt}{t})}(\mathbb R^d,\gamma_\infty)$, when $1<p<\infty$, and from $L^1(\mathbb R^d,\gamma_\infty)$ into $L^{1,\infty}_{L^2((0,\infty),\frac{dt}{t})}(\mathbb R^d,\gamma_\infty)$.
For that, we have to take into account that $T_{k,\ell ,m}^\nu(f)(x,t)=e^{R(Q^{1/2}x)}\mathbb{T}_{k,\ell,m}^\nu (f)(x,t)$, $x\in \mathbb{R}^d$ and $t>0$, and also that $\widetilde S_{k,i,j,{\rm loc}}^\nu(f)=\sum_{\ell,m=1}^da_{\ell,m}^ {i,j}T^\nu_{k,\ell,m,{\rm loc}}$, for certain constants $a_{\ell,m}^{i,j}\in \mathbb{R}$. To see this last property it is sufficient to note that
\begin{align}\label{partial2}
 \partial_{x_ix_j}^2  e^{-\frac{1}{2t}\langle Q^{-1}(y-x),y-x\rangle}&=\left(\langle\frac{1}{t}Q^{-1}e_j,y-x\rangle\langle\frac{1}{t}Q^{-1}e_i,y-x\rangle-\langle e_j,\frac{1}{t}Q^{-1}e_i\rangle\right)e^{-\frac{1}{2t}\langle Q^{-1}(y-x),y-x\rangle} \nonumber\\
& = \left(\langle\frac{1}{t}Q^{-1/2}e_j,Q^{-1/2}(y-x)\rangle\langle\frac{1}{t}Q^{-1/2}e_i,Q^{-1/2}(y-x)\rangle-\langle Q^{-1/2}e_j,\frac{1}{t}Q^{-1/2}e_i\rangle\right) \\
& \quad \times e^{-\frac{1}{2t}\langle Q^{-1/2}(y-x),Q^{-1/2}(y-x)\rangle},\quad x,y\in\mathbb R^d\mbox{ and }t>0.\nonumber
\end{align}

Suppose that $Q^{-1/2}=(c_{\ell m})_{\ell,m=1}^d$ with $c_{\ell m}\in\mathbb R$, $\ell,m=1,\ldots,d$. It follows that
\begin{align*}
\langle Q^{-1/2}e_j,y-x\rangle  \langle Q^{-1/2}e_i,y-x\rangle&=\langle(c_{\ell j})_{\ell=1}^d,y-x\rangle\langle(c_{\ell i})_{\ell=1}^d,y-x\rangle=\left(\sum_{\ell=1}^dc_{\ell j}(y_\ell-x_\ell)\right)\left(\sum_{\ell=1}^dc_{\ell i}(y_\ell-x_\ell)\right) \\
& =\sum_{\ell=1}^dc_{\ell j}c_{\ell i}(y_\ell-x_\ell)^2+\sum_{\substack{\ell,m=1\\
\ell\not=m}}^dc_{\ell j}c_{m i}(y_\ell-x_\ell)(y_m-x_m),
\end{align*}
for every $x=(x_1,\ldots,x_d)\in\mathbb R^d$, $y=(y_1,\ldots,y_d)$, and 
$$
\langle Q^{-1/2}e_j,Q^{-1/2}e_i\rangle=\langle(c_{\ell j})_{\ell=1}^d,(c_{\ell i})_{\ell=1}^d\rangle=\sum_{\ell=1}^d c_{\ell j}c_{\ell i}.
$$
We get
\begin{align*}
\left(\langle\frac{1}{t}Q^{-1/2}e_j,y-x\rangle\langle\frac{1}{t}Q^{-1/2}e_i,y-x\rangle-\langle Q^{-1/2}e_j,\frac{1}{t}Q^{-1/2}e_i\rangle\right)e^{-\frac{|y-x|^2}{2t}} \\
&\hspace{-8cm} =\Big(\frac{1}{t^2}\sum_{\substack{\ell,m=1\\l\not =m}}^dc_{\ell j}c_{m i}(y_\ell-x_\ell)(y_m-x_m)+\sum_{\ell=1}^dc_{\ell j}c_{\ell i}\Big(\frac{(y_\ell-x_\ell)^2}{t^2}-1\Big)\Big)e^{-\frac{|y-x|^2}{2t}} \\
&\hspace{-8cm} =\sum_{\substack{\ell,m=1\\\ell \not =m}}^dc_{\ell j}c_{mi}\partial_{x_\ell x_m}^2e^{-\frac{|y-x|^2}{2t}}+\sum_{\ell=1}^dc_{\ell j}c_{\ell i}\partial_{x_\ell x_\ell}^2e^{-\frac{|y-x|^2}{2t}} ,\quad x,y\in\mathbb R^d\mbox{ and }t>0.
\end{align*}

As in the Section \ref{S3.1.1}, in order to finish the proof it remains to establish that the operator $\mathcal{D}_{k,i,j}^\nu (f)(x)=g_{k,\alpha,{\rm loc}}^\nu (f)(x)-\|S_{k,i,j,{\rm loc}}^\nu(f)(x,\cdot)\|_{L^2((0,\infty ),\frac{dt}{t})}$, $x\in \mathbb{R}^d$, is bounded from $L^p(\mathbb{R}^d,\gamma_\infty)$ into itself, for each $1\leq p\leq \infty$. 

Minkowski inequality leads to
$$
e^{-R(y)}|\mathcal{D}_{k,i,j}^\nu (f)(x)|\leq C\int_{\mathbb R^d}\varphi_{\mbox{\tiny A}}(x,y)|f(y)|\big\|K_{k,i,j}^\nu (x,y,\cdot )-S_{k,i,j}^\nu (x,y,\cdot )\big\|_{L^2((0,\infty ),\frac{dt}{t})}e^{-R(y)}dy,\quad x\in \mathbb{R}^d.
$$
Thus, it is sufficient to prove that
\begin{equation}\label{diferenciakij}
e^{-R(y)}\big\|K_{k,i,j}^\nu (x,y,\cdot )-S_{k,i,j}^\nu (x,y,\cdot )\big\|_{L^2((0,\infty ),\frac{dt}{t})}\leq C\frac{1+|x|}{|x-y|^{d-1}},\quad (x,y)\in L_{2A}.
\end{equation}
Again, by Minkowski inequality we can write
\begin{align*}
e^{-R(y)}\big\|K_{k,i,j}^\nu (x,y,\cdot )-S_{k,i,j}^\nu (x,y,\cdot )\big\|_{L^2((0,\infty ),\frac{dt}{t})}&\\
&\hspace{-6cm}\leq Ce^{-R(y)}\int_0^\infty \|[s^{k+2}\partial _s^k\mathfrak{g}_\nu (s)]_{|s=\frac{t}{2\sqrt{u}}}\|_{L^2((0,\infty ),\frac{dt}{t})}|\partial ^2_{x_ix_j}h_u(x,y)-\mathbb{S}_u^{i,j}(x,y)|du\\
&\hspace{-6cm}\leq Ce^{-R(y)}\left(\int_0^{\mathfrak{m}(x)}+\int_{\mathfrak{m}(x)}^\infty\right) |\partial ^2_{x_ix_j}h_u(x,y)-\mathbb{S}_u^{i,j}(x,y)|du\\
&\hspace{-6cm}=I_{i,j}^0(x,y)+I_{i,j}^\infty (x,y),\quad x,y\in \mathbb{R}^d.
\end{align*}

First, we observe that from \eqref{3.2'} and (\ref{partial2}) we obtain that
\begin{align*}
e^{-R(y)}|\mathbb{S}_u^{i,j}(x,y)|&= C\frac{e^{R(x)-R(y)}}{u^{\frac{d}{2}}}\left|\partial^2_{x_ix_j}(e^{-\frac{1}{2u}\langle Q^{-1}(y-x),y-x\rangle})\right|\leq C\frac{e^{-c\frac{|y-x|^2}{u}}}{u^{\frac{d}{2}}}\Big(\frac{|y-x|^2}{u^2}+\frac{1}{u}\Big)\\
&\leq \frac{C}{u^{\frac{d}{2}+1}},\quad (x,y)\in L_{2A}\mbox{ and }u>0.
\end{align*}
Then, since $\mathfrak{m}(x)\sim (1+|x|)^{-2}$, $x\in \mathbb{R}^d$,
\begin{equation}\label{3.24}
   \int_{\mathfrak{m}(x)}^\infty e^{-R(y)}|\mathbb{S}_u^{i,j}(x,y)|du \leq C \int_{\mathfrak{m}(x)}^\infty\frac{du}{u^{\frac{d}{2}+1}} \leq \frac{C}{\mathfrak{m}(x)^{\frac{d}{2}}}\leq C(1+|x|)^d\leq C\frac{1+|x|}{|x-y|^{d-1}},\quad (x,y)\in L_{2A}. 
\end{equation}

On the other hand, by \eqref{(2.5)}, \eqref{(2.6)} and \eqref{3.2'} we get, for every $(x,y)\in L_{2A}$, 
$$
e^{-R(y)}|\partial^2_{x_ix_j}h_u(x,y)|  \leq 
 C\frac{e^{R(x)-R(y)}}{u^{\frac{d}{2}}}(1+|x|)^2\leq C\frac{(1+|x|)^2}{u^{\frac{d}{2}}},\quad \mbox{ when }\mathfrak{m}(x)<u<1,
$$
and
$$
e^{-R(y)}|\partial^2_{x_ix_j}h_u(x,y)|  \leq C e^{R(x)-R(y)}e^{-cu}(1+|y|)^2\leq C e^{-cu}(1+|x|)^2,\quad u\geq 1.
$$

By assuming that $d>2$ and using that $\mathfrak{m}(x)\sim (1+|x|)^{-2}$, $x\in \mathbb{R}^d$, it follows that 
\begin{align}\label{3.25}
\int_{\mathfrak{m}(x)}^\infty e^{-R(y)}  |\partial^2_{x_ix_j}h_u(x,y)|du&\leq  C(1+|x|)^2\left(\int_{\mathfrak{m}(x)}^\infty\frac{du}{u^{\frac{d}{2}}}+\int_1^\infty e^{-cu}du\right)\leq C(1+|x|)^2\Big(\frac{1}{\mathfrak{m}(x)^{\frac{d}{2}-1}}+1\Big)\nonumber\\
&\leq C(1+|x|)^d\leq C\frac{1+|x|}{|x-y|^{d-1}},\quad (x,y)\in L_{2A}.
\end{align}

By using together (\ref{3.24}),  and (\ref{3.25}) we obtain that
\begin{equation}\label{3.27}
    I_{i,j}^\infty (x,y)\leq C\frac{1+|x|}{|x-y|^{d-1}},\quad (x,y)\in L_{2A}. 
\end{equation}

We now deal with $I_{i,j}^0(x,y)$, $(x,y)\in L_{2A}$. Let us consider again $\widetilde{H}_u(x,y)$ and $\widetilde{W}_u(x,y)$, $x,y\in \mathbb{R}^d$ and $u>0$, as in \eqref{Htilde} and \eqref{Wtilde}, respectively. We write, for each $x,y\in \mathbb{R}^d$ and $u>0$, 
\begin{align}\label{Dif}
\partial ^2_{x_ix_j}h_u(x,y)-\mathbb{S}_u^{i,j}(x,y)&=\partial ^2_{x_ix_j}[h_u(x,y)-\widetilde{H}_u(x,y)]+\partial ^2_{x_ix_j}\widetilde{H}_u(x,y)-\mathbb{S}_u^{i,j}(x,y)\nonumber\\
&=\Big(1-\frac{({\rm det}\,Q_u)^{1/2}}{(u^d{\rm det }\,Q)^{1/2}}\Big)\partial ^2_{x_ix_j}h_u(x,y)+\partial ^2_{x_ix_j}\widetilde{H}_u(x,y)-\mathbb{S}_u^{i,j}(x,y),
\end{align}
where again we use that $h_u(x,y)=(u^d{\rm det}\,Q)^{1/2}({\rm det}\,Q_u)^{-1/2}\widetilde{H}_u(x,y)$, $x,y\in\mathbb{R}^d$ and $u>0$.

By considering \eqref{(2.5)}, \eqref{3.2'}, \eqref{3.0} and \eqref{3.1'} it follows that
\begin{align}\label{dif1}
e^{-R(y)}\left|\Big(1-\frac{({\rm det}\,Q_u)^{1/2}}{(u^d{\rm det }\,Q)^{1/2}}\Big)\partial ^2_{x_ix_j}h_u(x,y)\right|&\leq Cue^{R(x)-R(y)}\frac{e^{-c\frac{|y-D_ux|^2}{u}}}{u^{\frac{d}{2}}}\Big(|x|+\frac{1}{\sqrt{u}}\Big)^2\nonumber\\
&\hspace{-3cm}\leq C\frac{e^{-c\frac{|y-x|^2}{u}}}{u^{\frac{d}{2}}}(1+|x|)^2\leq C(1+|x|)\frac{e^{-c\frac{|y-x|^2}{u}}}{u^{\frac{d+1}{2}}},\quad (x,y)\in L_{2A} \mbox{ and }0<u<\mathfrak{m}(x).
\end{align}

On the other hand, according to \eqref{(2.2)}, we can write
\begin{align*}
 \partial_{x_ix_j}^2  \widetilde{H}_u(x,y)&=\widetilde{H}_u(x,y)(P_i(u,x,y)P_j(u,x,y)+\Delta_{i,j}(u)) \\
& = \widetilde{H}_u(x,y)\Big(\langle Q^{-1}_\infty x,e_i\rangle P_j(u,x,y) +\langle Q^{-1}_u e^{uB}e_i,y-D_u x\rangle\langle Q^{-1}_\infty x,e_j\rangle\Big)\\
& \quad +\widetilde{H}_u(x,y)\Big(\langle Q^{-1}_u e^{uB}e_j,y-D_u x\rangle\langle Q^{-1}_u e^{uB}e_i,y-D_u x\rangle +\Delta_{ij}\Big) \\
& =Z_{i,j}^1(x,y,u)+Z_{i,j}^2(x,y,u),\quad x,y\in\mathbb R^d\mbox{ and }u>0.
\end{align*}
Also, by considering \eqref{partial2} we can see that
$$
\mathbb{S}_u^{i,j}(x,y)=\widetilde{W}_u(x,y)\left(\langle \frac{1}{u}Q^{-1}e_j,y-x\rangle\langle \frac{1}{u}Q^{-1}e_i,y-x\rangle -\langle e_j,\frac{1}{u}Q^{-1}e_i\rangle\right),\quad x,y\in\mathbb R^d\mbox{ and }u>0.
$$
Then,
\begin{align*}
\Big|\partial ^2_{x_ix_j}\widetilde{H}_u(x,y)-\mathbb{S}_u^{i,j}(x,y)\Big|&\leq |Z_{i,j}^1(x,y,u)|+|Z_{i,j}^2(x,y,u)-\mathbb{S}_u^{i,j}(x,y)|\\
&\hspace{-3cm}\leq |Z_{i,j}^1(x,y,u)|+|\widetilde{H}_u(x,y)-\widetilde{W}_u(x,y)||\langle Q^{-1}_u e^{uB}e_j,y-D_u x\rangle\langle Q^{-1}_u e^{uB}e_i,y-D_u x\rangle +\Delta_{ij}|\\
&\hspace*{-3cm}\quad +\widetilde{W}_u(x,y)\Big|\Delta_{ij}+\langle e_j,\frac{1}{u}Q^{-1}e_i\rangle\Big|\\
&\hspace{-3cm}\quad +\widetilde{W}_u(x,y)\Big|\langle Q^{-1}_u e^{uB}e_j,y-D_u x\rangle\langle Q^{-1}_u e^{uB}e_i,y-D_u x\rangle -\langle \frac{1}{u}Q^{-1}e_j,y-x\rangle\langle \frac{1}{u}Q^{-1}e_i,y-x\rangle\Big|\\
&\hspace{-3cm}=\sum_{\ell =1}^4D_{i,j}^\ell (x,y,u),\quad x,y\in \mathbb{R}^d\mbox{ and }u>0.
\end{align*}

We take into account \eqref{3.2'}, \eqref{3.1'}, \eqref{3.6} and \cite[(4.5) and Lemma 2.2]{CCS3} to see that, for every $(x,y)\in L_{2A}$ and $0<u<\mathfrak{m}(x)$,
\begin{equation}\label{D1}
e^{-R(y)}D_{i,j}^1(x,y,u)\leq C|x|\frac{e^{-c\frac{|y-x|^2}{u}}}{u^{\frac{d}{2}}}\Big(|x|+\frac{1}{\sqrt{u}}\Big)\leq C|x|\frac{e^{-c\frac{|y-x|^2}{t}}}{u^{\frac{d+1}{2}}}\leq C(1+|x|)\frac{e^{-c\frac{|y-x|^2}{u}}}{u^{\frac{d+1}{2}}}.
\end{equation}
Now, by using \eqref{3.2'}, \eqref{3.4'}, \eqref{3.7.1} and \cite[(4.6) and Lemma 2.2]{CCS3} we get
\begin{align}\label{D2}
e^{-R(y)}D_{i,j}^2(x,y,u)&\leq C(1+|x|)\frac{e^{-c\frac{|y-x|^2}{u}}}{u^{\frac{d-1}{2}}}\Big(\frac{|y-D_ux|^2}{u^2}+\frac{1}{u}\Big)\leq
C(1+|x|)\frac{e^{-c\frac{|y-x|^2}{t}}}{u^{\frac{d-1}{2}}}\Big(|x|^2+\frac{1}{u}\Big)\nonumber\\
&\leq C(1+|x|)\frac{e^{-c\frac{|y-x|^2}{u}}}{u^{\frac{d+1}{2}}},\quad (x,y)\in L_{2A}\mbox{ and }0<u<\mathfrak{m}(x).
\end{align}
On the other hand, \eqref{3.2'}, \eqref{3.3'} lead to
\begin{align}\label{D3}
e^{-R(y)}D_{i,j}^3(x,y,u)&=e^{-R(y)}\widetilde{W}_u(x,y)\big|\langle e_j, \frac{1}{u}Q^{-1}e_i-e^{uB^*}Q_u^{-1}e^{uB}e_i\rangle \big|\leq C\frac{e^{-\frac{|Q^{-1/2}(y-x)|^2}{2u}}}{u^{\frac{d}{2}}}\nonumber\\
&\leq C(1+|x|)\frac{e^{-c\frac{|y-x|^2}{u}}}{u^{\frac{d+1}{2}}},\quad (x,y)\in L_{2A} \mbox{ and }0<u<1.
\end{align}
Finally, we estimate the term $e^{-R(y)}D_{i,j}^4(x,y,u)$, for $(x,y)\in L_{2A}$ and $0<u<\mathfrak{m}(x)$. 
We write
\begin{align*}
D_{i,j}^4(x,y,u)&=\widetilde{W}_u(x,y)\Big|\Big(\langle Q^{-1}_u e^{uB}e_j,y-D_u x\rangle-\langle \frac{1}{u}Q^{-1}e_j,y-x\rangle\Big)\langle Q^{-1}_u e^{uB}e_i,y-D_u x\rangle\\
&\quad+\langle \frac{1}{u}Q^{-1}e_j,y-x\rangle\Big(\langle Q^{-1}_u e^{uB}e_i,y-D_u x\rangle-\langle \frac{1}{u}Q^{-1}e_i,y-x\rangle\Big)\Big|,\quad x,y\in \mathbb{R}^d\mbox{ and }u>0.
\end{align*}
From (\ref{3.3'}), (\ref{3.4'}) and \cite[Lemmas 2.2 and 2.3]{CCS3} it follows that, for every $x,y\in\mathbb R^d$ and $u\in (0,1)$,
\begin{align*}
 \Big|\Big(\langle Q^{-1}_u e^{uB}e_j,y-D_u x\rangle-\langle \frac{1}{u}Q^{-1}e_j,y-x\rangle\Big)\langle Q^{-1}_u e^{uB}e_i,y-D_u x\rangle\Big| &\\
& \hspace{-8cm}\leq  \Big(\Big|\langle Q^{-1}_u e^{uB}e_j-\frac{1}{u}Q^{-1}e_j,y-D_u x\rangle\Big|+\Big|\langle \frac{1}{u}Q^{-1}e_j,x-D_u x\rangle\Big|\Big)|\langle Q^{-1}_u e^{uB}e_i,y-D_u x\rangle |\\
& \hspace{-8cm}\leq C \big(|y-D_u x|+|x|\big)\frac{|y-D_u x|}{u}\leq C\big(|y-x|+|x|\big)\Big(\frac{|y-x|}{u}+|x|\Big), 
\end{align*}
and, in the same way,
$$
\Big|\langle \frac{1}{u}Q^{-1}e_j,y-x\rangle\Big(\langle Q^{-1}_u e^{uB}e_i,y-D_u x\rangle-\langle\frac{1}{u} Q^{-1}e_i,y-x\rangle\Big)\Big| \leq  C\frac{|y-x|}{u}(|y-x|+|x|).
$$
Thus, by taking into account again \eqref{3.2'}, we obtain that
\begin{align}\label{D4}
e^{-R(y)}D_{i,j}^4(x,y,u)&\leq C\frac{e^{-c\frac{|y-x|^2}{u}}}{u^{\frac{d}{2}}}\big(|y-x|+|x|\big)\Big(\frac{|y-x|}{u}+|x|\Big)\leq C\frac{e^{-c\frac{|y-x|^2}{u}}}{u^{\frac{d}{2}}}(1+|x|)\Big(\frac{1}{\sqrt{u}}+|x|\Big)\nonumber\\
&\leq C(1+|x|)\frac{e^{-c\frac{|y-x|^2}{u}}}{u^{\frac{d+1}{2}}},\quad (x,y)\in L_{2A}\mbox{ and } 0<u<\mathfrak{m}(x).
\end{align}

From \eqref{D1}, \eqref{D2}, \eqref{D3} and \eqref{D4} we deduce that 
$$
e^{-R(y)}\Big|\partial ^2_{x_ix_j}\widetilde{H}_u(x,y)-\mathbb{S}_u^{i,j}(x,y)\Big|\leq C(1+|x|)\frac{e^{-c\frac{|y-x|^2}{u}}}{u^{\frac{d+1}{2}}},\quad (x,y)\in L_{2A}\mbox{ and } 0<u<\mathfrak{m}(x),
$$
that, jointly \eqref{Dif} and \eqref{dif1}, allows us to conclude that
\begin{equation}\label{H0}
I_{i,j}^0(x,y)\leq C(1+|x|)\int_0^{\mathfrak{m}(x)}\frac{e^{-c\frac{|y-x|^2}{u}}}{u^{\frac{d+1}{2}}}du\leq C\frac{1+|x|}{|x-y|^{d-1}},\quad (x,y)\in L_{2A}.
\end{equation}

By (\ref{3.27}) and (\ref{H0}) we get \eqref{diferenciakij} and then, we obtain that $\mathcal{D}_{k,i,j}^\nu(f)\leq C\mathscr{S}_{2A}(f)$. As consequence we obtain that $\mathcal{D}_{k,i,j}^\nu$ is bounded from $L^p(\mathbb{R}^d,\gamma_\infty)$ into itself, for each $1\leq p\leq \infty$. 

We conclude that the local square function $g^\nu_{k,\alpha,{\rm loc}}$ is bounded from $L^p(\mathbb R^d,\gamma_\infty)$ into itself, for every $1<p<\infty$, and from $L^1(\mathbb R^d,\gamma_\infty)$ into $L^{1,\infty}(\mathbb R^d,\gamma_\infty)$.

\subsubsection{About the global part of $g_{k,\alpha}^\nu$: $0<\widehat{\alpha}\leq 2$}\label{S3.1.3}
In this subsection we choose $A$ large enough as in \cite[Section 9]{CCS3}. By using Minkowski inequality we get
\begin{align}\label{3.28}
g_{k,\alpha,{\rm glob}}^\nu (f)(x)&\leq C\int_{\mathbb R^d}(1-\varphi_{\mbox{\tiny A}} (x,y))|f(y)|\|P_{k,\alpha}^\nu (x,y,\cdot )\|_{L^2((0,\infty ),\frac{dt}{t})}d\gamma_\infty(y)\nonumber\\
&\leq C\int_{\mathbb R^d}(1-\varphi_{\mbox{\tiny A}} (x,y))|f(y)|\int_0^\infty \big\|[s^{k+\widehat{\alpha}}\partial _s^k{\mathfrak g}_\nu(s)]_{\big|s=\frac{t}{2\sqrt{u}}}\big\|_{L^2((0,\infty ),\frac{dt}{t})}|\partial _x^\alpha h_u(x,y)|u^{\frac{\widehat{\alpha}}{2}-1}dud\gamma_\infty(y)\nonumber\\
&\leq C\int_{\mathbb R^d}(1-\varphi_{\mbox{\tiny A}} (x,y))|f(y)|\int_0^\infty|\partial _x^\alpha h_u(x,y)|u^{\frac{\widehat{\alpha}}{2}-1}dud\gamma_\infty(y),\quad x\in \mathbb R^d.
\end{align}
We recall that ${\mathfrak g}_\nu(s)=s^{2\nu}e^{-s^2}$, $s>0$. By using \cite[Lemma 4.1, Corollary 5.3 and  Propositions 7.1 and 9.1]{CCS3} and by taking into account that in the global operators the cancellation property of the integral kernels does not play any role, we can deduce from \eqref{3.28} that the global square functions $g_{k,\alpha,{\rm glob}}^\nu$ are bounded from $L^1(\mathbb R^d,\gamma_\infty)$ into $L^{1,\infty }(\mathbb R^d,\gamma_\infty )$.
\subsection{}\label{S3.2}
We consider $\alpha = 0$ and $k\geq 1$. Let $f\in C^\infty_c(\mathbb{R}^d)$. According to \eqref{alpha0} we can write
$$
t^k\partial^k_t P^\nu_t (f)( x) = 
        \int_{\mathbb{R}^d}P^{\nu}_k (x,y,t)f(y)d\gamma_\infty(y),
        \quad  x\in\mathbb{R}^d
        \text{ and }
        t>0,
$$
where
$$
P^\nu_k(x,y,t)=\frac{2}{\Gamma(\nu)}\int_0^\infty [s^k\partial^{k-1}_s\mathfrak{h}_\nu(s)]_{|s=\frac{t}{2\sqrt{u}}}\partial_u h_u(x,y) du,\quad  x,y\in\mathbb{R}^d\text{ and } t>0,
$$
with $\mathfrak{h}_\nu(s) = s^{2\nu-1}e^{-s^2}$, $s>0$. We observe that $\mathfrak{h}_\nu(s)=s^{-1}\mathfrak{g}_\nu(s)$, $s>0$. A simple computation shows that, for every $\ell \in \mathbb{N}$, we have that $|s^{\ell +1}\partial_s^\ell \mathfrak{h}_\nu (s)|\leq Cg_\nu (s/2)$, $s>0$. Then, $\|s^{\ell +1}\partial _s^\ell \mathfrak{h}_\nu(s)\|_{L^2((0,\infty ),\frac{ds}{s})}\leq C$, for each $\ell \in \mathbb{N}$.

According to \eqref{duxd2x}
$$
\partial_u h_u(x,y) \sum_{i,j=1}^{d} (c_{i,j}\partial^2_{x_ix_j}h_u(x,y) +
        d_{i,j} x_i\partial_{x_j}h_u(x,y))
        ,\quad 
        x,y\in\mathbb{R}^d
        \text{ and } u>0,
$$
for certain $c_{i,j}$ and $d_{i,j}\in\mathbb{R}$, $i,j=1,\dots,d$.

We define, for every $i,j=1,\dots,d$,
\begin{equation*}
    \mathcal{K}^\nu_{k,i,j}(x,y,t) =\int_0^\infty [s^k\partial^{k-1}_s\mathfrak{h}_\nu(s)]_{|s=\frac{t}{2\sqrt{u}}}\partial^2_{x_ix_j}
        h_u(x,y) du,\quad 
        x,y\in\mathbb{R}^d
        \text{ and } t>0,
\end{equation*}
and
\begin{equation*}
    \mathcal{R}^\nu_{k,i,j}(x,y,t) = 
    \int_0^\infty[s^k\partial^{k-1}_s\mathfrak{h}_\nu(s)]_{|s=\frac{t}{2\sqrt{u}}}
        x_i \partial_{x_j}
        h_u(x,y) du,\quad 
        x,y\in\mathbb{R}^d
        \text{ and } t>0.
\end{equation*}
and consider the functions given by
\begin{equation*}
    \mathscr{K}^{\nu}_{k,i,j} (f) (x,t) 
    =    \int_{\mathbb{R}^d}
    \mathcal{K}^\nu_{k,i,j}(x,y,t) f(y)
    d\gamma_\infty(y),
    \quad  x\in \mathbb{R}^d \mbox{ and }t>0,
\end{equation*}
and
\begin{equation*}
    \mathscr{R}^{\nu}_{k,i,j} (f) (x,t) 
    = \int_{\mathbb{R}^d}
    \mathcal{R}^\nu_{k,i,j}(x,y,t) f(y)
    d\gamma_\infty(y),
    \quad x\in \mathbb{R}^d\mbox{ and }t>0.
\end{equation*}

It is clear that 
\begin{equation*}
    g^\nu_{k,0}(f)(x) \leq C \sum_{i,j=1}^d
    (\|\mathscr{K}^\nu_{k,i,j}(f)(x,\cdot )\|_{L^2((0,\infty ),\frac{dt}{t})} + \|\mathscr{R}^\nu_{k,i,j}(f)(x,\cdot)\|_{L^2((0,\infty ),\frac{dt}{t})}),\quad x\in \mathbb{R}^d.
\end{equation*}
We define the local and the global parts of $g^\nu_{k,0}$, $\mathscr{K}^\nu_{k,i,j}$ and $\mathscr{R}^\nu_{k,i,j}$ in the usual way.
\subsubsection{About the local part of $g_{k,0}^\nu$}\label{S3.2.1}
Let $i, j=1,\dots,d$. We firstly study $\mathscr{R}^\nu_{k,i,j,\loc}$. Minkowski inequality leads to
$$
    \|\mathscr{R}^\nu_{k,i,j,\loc}(f)(x,\cdot )\|_{L^2((0,\infty ),\frac{dt}{t})} \leq C\int_{\mathbb{R}^d} \varphi_{\mbox{\tiny A}}(x,y)|f(y)| \|\mathcal{R}_{k,i,j}^\nu (x,y,\cdot)\|_{L^2((0,\infty ),\frac{dt}{t})}e^{-R(y)}dy.
$$
We assert that, when $d>2$,
\begin{equation}\label{Rkij}
    e^{-R(y)}\|\mathcal{R}_{k,i,j}^\nu (x,y,\cdot)\|_{L^2((0,\infty ),\frac{dt}{t})}\leq C \frac{1+|x|}{|x-y|^{d-1}},
        \quad  (x,y)\in L_{2A}.
\end{equation}
In this way, we obtain that $\|\mathscr{R}^\nu_{k,i,j,\loc}(f)(x,\cdot )\|_{L^2((0,\infty ),\frac{dt}{t})}\leq C\mathscr{S}_\eta(|f|)(x)$, $x\in \mathbb{R}^d$, with $\eta =2A$, and we conclude that $\mathscr{R}^\nu_{k,i,j,\loc}$ is bounded from $L^p(\mathbb{R}^d, \gamma_\infty)$ into $L^p_{L^2((0,\infty ),\frac{dt}{t})}(\mathbb{R}^d, \gamma_\infty)$, for every $1\leq p\leq \infty$.

Let us establish \eqref{Rkij}. Since
$$
\|s^{k}\partial _s^{k-1} \mathfrak{h}_\nu(s)_{|s=\frac{t}{2\sqrt{u}}}\|_{L^2((0,\infty ),\frac{ds}{s})}=\|s^k\partial _s^{k-1}\mathfrak{h}_\nu(s)\|_{L^2((0,\infty ),\frac{ds}{s})},\quad t,u>0,
$$
by using again Minkowski inequality and \eqref{(2.3)}, \eqref{(2.4)}, \eqref{3.2'} and \eqref{3.1'} we can write, when $d>2$,
\begin{align*}
     e^{-R(y)}\|\mathcal{R}_{k,i,j}^\nu (x,y,\cdot)\|_{L^2((0,\infty ),\frac{dt}{t})} &\leq C|x|e^{-R(y)}\int_0^\infty
    |\partial_{x_j}h_u(x,y)|du\\
&
    \hspace*{-3cm}
    \leq C|x|e^{R(x)-R(y)}
        \left(
        \int_0^{\mathfrak{m}(x)} \frac{e^{-c\frac{|y-x|^2}{u}}}{u^{\frac{d}{2}}}\Big(|x|+\frac{1}{\sqrt{u}}\Big)du
        +
       \int_{\mathfrak{m}(x)}^1\Big(|x|+\frac{1}{\sqrt{u}}\Big)
        \frac{du}{u^{\frac{d}{2}}}
                +
        (1+|y|)\int_1^\infty e^{-cu} du 
        \right)
        \\
&
    \hspace*{-3cm}
    \leq C|x|
        \left(
        \int_0^\infty \frac{e^{-\frac{c|y-x|^2}{u}}}{u^{\frac{d+1}{2}}}du
        +
       (1+|x|)\int_{\mathfrak{m}(x)}^\infty
        \frac{du}{u^{\frac{d}{2}}}
                + 1+|x|
        \right)
        \\
&\hspace*{-3cm} \leq C|x|
        \left( \frac{1}{|x-y|^{d-1}}
        +
        \frac{1+|x|}{\mathfrak{m}(x)^{\frac{d}{2}-1}}
        +
        1+|x|
        \right)\leq C|x|\left(\frac{1}{|x-y|^{d-1}}+(1+|x|)^{d-1}\right)
        \\ 
&\hspace*{-3cm} \leq C \frac{1+|x|}{|x-y|^{d-1}},
        \quad  (x,y)\in L_{2A}.
\end{align*}
In the third inequality we have taken into account that $|y|\leq C(1+|x|)$, when $(x,y)\in L_\eta$, $\eta >0$, and that $\mathfrak{m}(x)\sim (1+|x|)^{-2}$ which leads to the following estimates
$$
|x|+\frac{1}{\sqrt{u}}\leq C
\left\{\begin{array}{ll}
\displaystyle \frac{1}{\sqrt{u}},&0<u\leq \mathfrak{m}(x),\\[0.5cm]
1+|x|, &u\geq \mathfrak{m}(x),
\end{array}
\right., \quad x\in \mathbb{R}^d.
$$

In order to study $\mathscr{K}^\nu_{k,i,j,{\rm loc}}$ we proceed as in Section \ref{S3.1.2}. We define the operator \begin{equation*}
    \mathscr{S}^\nu_{k,i,j}(f)(x,t) = 
\int_{\mathbb{R}^d} \mathscr{S}^\nu_{k,i,j}(x,y,t)f(y) d\gamma_\infty (y),
\quad  x\in \mathbb{R}^d \text{ and } t>0.
\end{equation*}
where
\begin{equation*}
    \mathscr{S}^\nu_{k,i,j}(x,y,t) =
    \int_0^\infty [s^k\partial^{k-1}_s\mathfrak{h}_\nu(s)]_{|s=\frac{t}{2\sqrt{u}}}\mathbb{S}_u^{i,j}(x,y)du, \quad  x, y \in \mathbb{R}^d \text{ and } t>0.
\end{equation*}
Here $\mathbb{S}_u^{i,j}(x,y)$, $x,y\in \mathbb{R}^d$, $u>0$, is given by \eqref{Su}.

By using Minkowski inequality it follows that
$$
\big\|\mathcal{K}_{k,i,j}^\nu (x,y,\cdot)-\mathscr{S}^\nu_{k,i,j}(x,y,\cdot)\big\|_{L^2((0,\infty),\frac{dt}{t})}\leq C\int_0^\infty |\partial_{x_ix_j}^2h_u(x,y)-\mathbb{S}_u^{i,j}(x,y)|du,\quad x,y\in \mathbb{R}^d.
$$
Then, defining the local operator $\mathscr{S}^\nu_{k,i,j,\loc}$ in the usual way, from the proof of estimate \eqref{diferenciakij} we get that 
$$
\big\|\mathscr{K}^\nu_{k,i,j,\loc}(f)(x,\cdot ) -\mathscr{S}^\nu_{k,i,j,\loc}(f)(x,\cdot )\big\|_{L^2((0,\infty),\frac{dt}{t})}\leq C\mathscr{S}_{2A}(f)(x),\quad  x\in \mathbb{R}^d.
$$
On the other hand, consider the square function $\mathscr{G}^\nu_{k,i,j}$ defined by 
$$
\mathscr{G}_{k,i,j}^\nu (f)(x)=\big\|\mathscr{T}_{k,i,j}^\nu (f)(x,\cdot )\big\|_{L^2((0,\infty),\frac{dt}{t})},\quad x\in \mathbb{R}^d,
$$
where
\begin{equation*}
    \mathscr{T}^\nu_{k,i,j} (f)(x,t)
    = 
    \int_{\mathbb{R}^d}
    \mathcal{T}^\nu_{k,i,j} (x-y,t)f(y)dy,
        \quad  x\in\mathbb{R}^d\mbox{ and }t>0,
\end{equation*}
and
\begin{equation*}
\mathcal{T}^\nu_{k,i,j}(z,t) = \int_0^{\infty}
 [s^k\partial^{k-1}_s\mathfrak{h}_\nu(s)]_{|s=\frac{t}{2\sqrt{u}}}
\partial^2_{x_ix_j}\mathbb{W}_u(z)du,\quad z\in\mathbb{R}^d \mbox{ and }t>0.
\end{equation*}
By arguing as in the analysis of $\mathbb{G}_{k,\ell,m}^ \nu$ (see Section \ref{S3.1.2}) we can prove that the operator $\mathscr{G}^\nu_{k,i,j}$ is bounded from $L^p(\mathbb{R}^d,dx)$ into itself, for every $1<p<\infty$ and from $L^1(\mathbb{R}^d,dx)$ into $L^{1,\infty}(\mathbb{R}^d,dx)$ and that $\mathscr{S}_{k,i,j,{\rm loc}}^\nu$ is a bounded operator $L^p(\mathbb R^d,\gamma_\infty)$ into  $L^p_{L^2((0,\infty),\frac{dt}{t})}(\mathbb R^d,\gamma_\infty)$, when $1<p<\infty$, and from $L^1(\mathbb R^d,\gamma_\infty)$ into $L^{1,\infty}_{L^2((0,\infty),\frac{dt}{t})}(\mathbb R^d,\gamma_\infty)$.

As in previous sections we can conclude that $g^\nu_{k,0,\loc}$ is bounded from $L^p(\mathbb{R}^d,\gamma_\infty)$ into itself, when $1<p<\infty$, and from $L^1(\mathbb{R}^d,\gamma_\infty)$ into $L^{1,\infty}(\mathbb{R}^d,\gamma_\infty)$. 

\subsubsection{About the global part of $g_{k,0}^\nu$}\label{S3.2.2}
We assume that $A$ is large enough as in~\cite[Section 9]{CCS3}. Let $i,j=1,\dots,d$. We can write
\begin{equation*}
   \big\|\mathcal{K}^{\nu}_{k,i,j}(x,y,\cdot)\big\|_{L^2((0,\infty ),\frac{dt}{t})}
    \leq C\int_0^\infty  |\partial^2_{x_ix_j} h_u(x,y)| du,
    \quad  x,y\in\mathbb{R}^d.
\end{equation*}
and
\begin{equation*}
    \big\|\mathcal{R}^{\nu}_{k,i,j}(x,y,\cdot )\big\|_{L^2((0,\infty ),\frac{dt}{t})}
    \leq C \int_0^\infty |x_i| |\partial_{x_j} h_u(x,y)| du,
    \quad  x,y\in\mathbb{R}^d.
\end{equation*}

According to~\cite[Lemma 4.1 and Corollary 5.3]{CCS3} we get
\begin{equation}\label{3.30}
    |x_i|\int_1^\infty 
    |\partial_{x_j} h_u(x,y)|du
    +
    \int_1^\infty |\partial^2_{x_i,x_j} h_u(x,y)|du\leq 
    C (1+|x|) e^{R(x)}, \quad  x,y\in\mathbb{R}^d.
\end{equation}
On the other hand, according to \eqref{(2.3)} and \eqref{(2.5)} we obtain
\begin{equation}\label{3.31}
    |x_i| 
    |\partial_{x_j} h_u(x,y)|
    +
     |\partial^2_{x_i,x_j} h_u(x,y)| \leq 
    C  \frac{e^{R(x)}}{
    u^{\frac{d}{2}}} 
    e^{-c\frac{|y-D_u x|^2}{u}}\Big(|x| +
    \frac{1}{\sqrt{u}}
    \Big)^2,\quad x,y\in\mathbb{R}^d \mbox{ and }u\in (0,1).
\end{equation}

The estimates in~\eqref{3.30} and~\eqref{3.31} allow us to proceed as in~\cite[Propositions 7.1 and 9.1]{CCS3} to establish that the global operator $g^{\nu}_{k,0,\glob}$ is bounded from $L^1(\mathbb{R}^d,\gamma_\infty)$ into $L^{1,\infty}(\mathbb{R}^d,\gamma_\infty)$.
 
\section{Proof of Theorem \ref{Th1.2}}\label{section4}
Let $f\in C_c^\infty (\mathbb R^d)$. According to \eqref{tkalpha} we have that
$$
t^{k+\widehat{\alpha}}\partial _t^k\partial _x^\alpha P_t^\nu (f)(x)=\frac{2^{\widehat{\alpha}}}{\Gamma (\nu)}\int_0^\infty [s^{k+\widehat{\alpha}}\partial _s^k\mathfrak{g}_\nu (s)]_{|s=\frac{t}{2\sqrt{u}}}\int_{\mathbb R^d}f(y)\partial _x^\alpha h_u(x,y)d\gamma_\infty (y)u^{\frac{\widehat{\alpha}}{2}-1}du,\quad x\in \mathbb R^d\mbox{ and }t>0.
$$

Suppose that $\alpha \not=0$. By taking into account that $|s^k\partial _s^k\mathfrak{g}_\nu(s)|\leq C\mathfrak{g}_\nu (s/2)$, $s>0$, and using estimates \eqref{(2.10.1)} and \eqref{(2.11)} we obtain that
\begin{align*}
    |t^{k+\widehat{\alpha}}\partial _t^k\partial _x^\alpha P_t^\nu (f)(x)|&\leq C\int_0^\infty [s^{\widehat{\alpha}}\mathfrak{g}_\nu \big(\frac{s}{2}\big)]_{|s=\frac{t}{2\sqrt{u}}}\int_{{\rm supp }f}|\partial _x^\alpha h_u(x,y)|dyu^{\frac{\widehat{\alpha}}{2}-1}du\\
    &\hspace{-2cm}\leq Ct^{\widehat{\alpha}}e^{R(x)}\sum_{n=0}^{[\frac{\widehat{\alpha}}{2}]}\left(\int_0^1e^{-c\frac{t^2}{u}}\Big(|x|+\frac{1}{\sqrt{u}}\Big)^{\widehat{\alpha}-2n}u^{-n-\frac{d}{2}-1}du+\int_1^\infty e^{-c\frac{t^2}{u}}e^{-cu}\int_{{\rm supp}f}(1+|y|)^{\widehat{\alpha}-2n}dy\frac{du}{u}\right)\\
    &\hspace{-2cm}\leq C(x)t^{\widehat{\alpha}}\left(\int_0^1e^{-c\frac{t^2}{u}}u^{-\frac{\widehat{\alpha}}{2}-\frac{d}{2}-1}du+t^{-\widehat{\alpha}-d}\int_1^\infty u^{\frac{\widehat{\alpha}}{2}+\frac{d}{2}} e^{-cu}du\right)\leq \frac{C(x)}{t^d},\quad x\in \mathbb{R}^d\mbox{ and }t>0.
\end{align*}
Thus, 
$$
\lim_{t\rightarrow +\infty}\partial_t^k\partial_x^\alpha P_t^\nu (f)(x)=0,\quad x\in \mathbb R^d.
$$
We can write
$$
\partial_t^k\partial_x^\alpha P_t^\nu (f)(x)=-\int_t^\infty \partial _u^{k+1}\partial _x^\alpha P_u^\nu (f)(x)du,\quad x\in \mathbb R^d.
$$
Then, by proceeding as in \cite[p. 478]{LS} we obtain that
$$
P_{*,k,\alpha}^\nu (f)(x)\leq Cg_{k+1,\alpha}^\nu (f)(x),\quad x\in \mathbb R^d,
$$
and Theorem \ref{Th1.2} follows from Theorem \ref{Th1.1}.

On the other hand, we can see that
$$
P_{*,k,0}^\nu (f)(x)\leq C\sup_{t>0}|\mathcal H_t(f)(x)|,\quad x\in \mathbb R^d.
$$
By using \cite[Theorem 1.1 and Corollary 1.2]{CCS2} we can establish Theorem \ref{Th1.2} for $\alpha =0$.

\section{Proof of Theorem \ref{Th1.3}}\label{section5}
We show Theorem~\ref{Th1.3} for the variation operator $V_\rho(\{t^{k+\widehat{\alpha}} \partial^k_t \partial^{\alpha}_x P^\nu_t\}_{t>0})$.
The result for the oscillation operator can be proved in a similar way.

We consider the operator $\mathbb{T}^\nu_{k,\alpha}$ defined by
\begin{equation*}
    \mathbb{T}^\nu_{k,\alpha} (f)(x,t) = t^{k + \widehat{\alpha}} \partial^k_t \partial^\alpha_x \mathbb{P}^\nu_t (f)(x),
    \quad x\in \mathbb{R}^d \text{ and } t>0.
\end{equation*}
We recall that
\begin{equation*}
    \mathbb{P}^\nu_t (f)(x) = \frac{t^{2\nu}}{4^\nu\Gamma(\nu)}
    \int_{\mathbb{R}^d}e^{-\frac{t^2}{4u}}\mathbb{W}_u(f)(x)\frac{du}{u^{1+\nu}},
    \quad  x\in \mathbb{R}^d \text{ and } t>0,
\end{equation*}
where
\begin{equation*}
    \mathbb{W}_u(f)(x)=\int_{\mathbb{R}^d}\mathbb{W}_u(x-y)f(y)dy,
     \quad x\in \mathbb{R}^d 
     \text{ and } u>0,
\end{equation*}
and $\mathbb{W}_u(z) = (2\pi u)^{-\frac{d}{2}} e^{-\frac{|z|^2}{2u}} $, $z\in\mathbb{R}^d$ and $u>0$.

Let $f\in C^{\infty}_c(\mathbb{R}^d)$. As in \eqref{tkalpha} we can write
\begin{align*}
  \mathbb{T}^\nu_{k,\alpha} (f)(x,t) &= \frac{2^{\widehat{\alpha}}}{\Gamma(\nu)}\int_0^\infty [s^{k+\widehat{\alpha}}\partial_s^k{\mathfrak g}_\nu(s)]_{|s=\frac{t}{2\sqrt u}}\partial_x^\alpha \mathbb{W}_u(f)(x)u^{\frac{\widehat{\alpha}}{2}-1}du\\
 &=\frac{2^{\widehat{\alpha}+1}}{\Gamma(\nu)}\int_0^\infty s^{k+\widehat{\alpha}-1}\partial_s^k{\mathfrak g}_\nu(s)\big[u^{\frac{\widehat{\alpha}}{2}}\partial_x^\alpha \mathbb{W}_u(f)(x)\big]_{|u={\frac{t^2}{4s^2}}}ds,\quad x\in \mathbb R^d\mbox{ and }t>0,\nonumber
 \end{align*}
 
It follows that
\begin{equation}\label{4.0}
    \begin{split}
        V_\rho (t\rightarrow \mathbb{T}^\nu_{k,\alpha}(f)(x,t)) & 
        \leq C V_\rho (t\rightarrow t^{\frac{\widehat{\alpha}}{2}} \partial^\alpha_x \mathbb{W}_t(f)(x)) \int_0^\infty s^{k + \widehat{\alpha}-1}\partial^k_s[{\mathfrak g}_\nu(s)] ds 
               \\ & \leq C 
        V_\rho (t\rightarrow t^{\frac{\widehat{\alpha}}{2}} \partial^\alpha_x \mathbb{W}_t(f)(x)),
        \quad x\in\mathbb{R}^d.
    \end{split}
\end{equation}

Firstly we deal with the case $\widehat{\alpha}=1$. Let $j=1,\dots,d$. By considering the Fourier transform $\mathcal{F}$ in $\mathbb{R}^d$ given in \eqref{Fouriertransform} we have that
\begin{equation*}
    \mathcal{F}(\partial_{z_j}\mathbb{W}_t(f))(y) =i y_je^{-\frac{1}{2}t|y|^2}\mathcal{F}(f)(y)=-|y|e^{-\frac{1}{2}t|y|^2}\mathcal{F}(\mathcal{R}_jf)(y), \quad y\in \mathbb{R}^d
    \text{ and } t>0,
\end{equation*}
where $\mathcal R_j$ represents the $j$-th Euclidean Riesz transform. 

We define $\phi = -\mathcal{F}^{-1}\{(2\pi )^{-\frac{d}{2}}|y| e^{-\frac{|y|^2}{2}}\}$. Then $\phi$ is a radial and smooth function in $\mathbb{R}^d$ and satisfies
\begin{equation*}
    \frac{1}{t^{\frac{d}{2}}}\phi\left(
    \frac{x}{\sqrt{t}}
    \right) =-
    \mathcal{F}^{-1}\{(2\pi )^{-\frac{d}{2}}
    \sqrt{t} |y| e^{-\frac{1}{2}t|y|^2}
    \}(x), \quad
    x\in \mathbb{R}^d \text{ and }
    t>0.
\end{equation*}

If we denote by $\phi_s(x) =s^{-d}\phi(x/s)$, $x\in \mathbb R^d$ and $s>0$, then we have that
\begin{equation*}
    \mathcal{F}(\sqrt{t}\partial_{z_j}\mathbb{W}_t(f))(y) =(2\pi)^{\frac{d}{2}} \mathcal{F}(\phi_{\sqrt{t}})(y)\mathcal{F}(\mathcal{R}_jf)(y), \quad y\in \mathbb{R}^d
    \text{ and } t>0,
\end{equation*}
that is,
\begin{equation*}
    \sqrt{t} \partial_{z_j} \mathbb{W}_t(f)(z) = 
    (\phi_{\sqrt{t}} *(\mathcal R_j f))(z),
    \quad z\in \mathbb{R}^d
    \text{ and } t>0.
\end{equation*}
Here $*$ denotes the usual convolution in $\mathbb R^d$.

We represent by $h_\mu$ the $\mu$-Hankel transform defined by
\begin{equation*}
    h_\mu (\Psi)(u) = 
    \int_0^\infty (uv)^{-\mu}J_\mu
    (uv) \Psi(v) v^{2\mu+1} dv,
    \quad u>0.
\end{equation*}
We have that
\begin{equation*}
    \phi(x) =-(2\pi )^{-\frac{d}{2}}
    h_{\frac{d-2}{2}} (v e^{-\frac{v^2}{2}})(|x|),
    \quad x\in\mathbb{R}^d.
\end{equation*}

According to~\cite[p. 30, (14)]{EMOT}, we get

\begin{equation*}
    \begin{split}
        h_{\frac{d-2}{2}}(ve^{-\frac{v^2}{2}})(u) & 
        = \int_0^\infty (uv)^{\frac{2-d}{2}} 
        J_{\frac{d-2}{2}}(uv) ve^{-\frac{v^2}{2}} v^{d-1} dv=
        u^{\frac{1-d}{2}} \int_0^\infty \sqrt{uv}  J_{\frac{d-2}{2}}(uv)
        e^{-\frac{v^2}{2}} v^{\frac{d+1}{2}} dv
        \\ & 
        = 
        \frac{\sqrt{2}\Gamma(\frac{d+1}{2})}{\Gamma(\frac{d}{2})}
        \,\prescript{}{1}F_1
        \left(\frac{d+1}{2}, \frac{d}{2}; -\frac{u^2}{2}\right), 
        \quad u>0,
    \end{split}
\end{equation*}
where $\prescript{}{1}F_1$ denotes the Kummer's confluent hypergeometric function.

According to~\cite[(9.9.4), (9.11.2) and (9.12.8)]{Leb} (note that in \cite[(9.11.2)]{Leb} there must be a minus sign in the argument), we obtain
\begin{equation*}
    \lim_{u\rightarrow \infty} \prescript{}{1}F_1 
    \left(\frac{d+1}{2}, \frac{d}{2}; -\frac{u^2}{2}\right) = 0,
\end{equation*}
and
\begin{equation*}
    \begin{split}
        \int_0^\infty \left|\partial_u \prescript{}{1}F_1 
    \left(\frac{d+1}{2}, \frac{d}{2}; -\frac{u^2}{2}\right)
    \right| u^d du
    & = 
     \frac{d+1}{d} \int_0^\infty \left| \prescript{}{1}F_1 
    \left(\frac{d+3}{2}, \frac{d+2}{2}; -\frac{u^2}{2}\right)
    \right| u^{d+1} du
   \\ &  \leq C \left( \int_0^1 du
    + \int_1^\infty \frac{du}{u^2} \right)
    < \infty.
    \end{split}
\end{equation*}

By using~\cite[Lemma 2.4]{CJRW1}, the variation operator defined by the family
$\{T_t\}_{t>0}$, where for every $t>0$, $T_t(f) = \phi_{\sqrt{t}}*f$, 
is bounded from $L^{p}(\mathbb{R}^d,dx)$ into itself, for every $1<p<\infty$, 
and from $L^{1}(\mathbb{R}^d,dx)$ into $L^{1,\infty}(\mathbb{R}^d,dx)$.

Since $\mathcal{R}_j$ is bounded from $L^{p}(\mathbb{R}^d,dx)$ into itself, for every $1<p<\infty$,  we conclude that $V_\rho(\{t^{1/2} \partial_x \mathbb{W}_t\}_{t>0})$ is bounded from $L^{p}(\mathbb{R}^d,dx)$ into itself, for every $1<p<\infty$.

In order to see that $V_\rho(\{t^{1/2} \partial_x \mathbb{W}_t\}_{t>0})$ is bounded from $L^{1}(\mathbb{R}^d,dx)$ into $L^{1,\infty}(\mathbb{R}^d,dx)$ we use the Calder\'on-Zygmund theory for vector valued singular integrals.

We consider the space $E_\rho$ that consists of complex functions $g$ defined on $(0,\infty)$ such that $V_\rho(g)<\infty$. It is clear that $V_\rho(g) = 0$ if and only if $g$ is constant in $(0,\infty)$. $V_\rho$ defines a norm on the corresponding quotient space. Thus $(E_\rho,V_\rho)$ can be seen as a Banach space.

Suppose that $g$ is a derivable complex function on $(0,\infty)$. If $0<t_n < t_{n-1}< \dots < t_1$ we get
\begin{equation*}
    \left(\sum_{i=1}^{n-1} |g(t_{i+1}) - g(t_i)|^\rho
    \right)^{1/\rho} 
    =
    \sum_{i=1}^{n-1} \left|
    \int_{t_{i+1}}^{t_{i}} g'(s) ds
    \right|
    \leq 
    \int_0^\infty |g'(t)|dt.
\end{equation*}
Then
\begin{equation}\label{4.1}
    V_\rho(g) \leq \int_0^\infty |g'(t)| dt.
\end{equation}

We consider the operator $\tau_j$ defined by
\begin{equation*}
    \tau_j(f) (x,t) = \sqrt{t} \partial_{x_j} \mathbb{W}_t(f)(x)
    =
    \int_{\mathbb{R}^d} \sqrt{t} \partial_{x_j} \mathbb{W}_t
    (x-y) f(y)dy, \quad x\in\mathbb{R}^d 
    \text{ and } t>0.
\end{equation*}

Let $f\in C^\infty_c(\mathbb{R}^d)$. For every $x\in\mathbb{R}^d$, 
the function $t\rightarrow \tau_j(f)(x,t)$ is continuous in $(0,\infty)$.
Then,
\begin{equation*}
    V_{\rho} ( \{\tau_j (f)(x,t) \}_{t>0})
    = \sup_{\substack{ 0<t_n<\dots<t_1 \\ t_i\in\mathbb{Q}, \, i=1,\dots,n}}
    \left(
    \sum_{i=1}^{n-1} |\tau_j(f)(x, t_{i+1}) - 
    \tau_j(f)(x,t_i)|^\rho
    \right)^{1/\rho},\quad x\in \mathbb{R}^d.
\end{equation*}
Hence, $V_\rho(\{ \tau_j(f)(x,t)\}_{t>0})$ is a Lebesgue measurable
function on $\mathbb{R}^d$.

By~\eqref{4.1} we obtain
\begin{equation}\label{4.2}
    V_{\rho}\left(
    \{\sqrt{t} \partial_{z_j} \mathbb{W}_t(z)\}_{t>0}
    \right) 
       \leq
    \int_0^\infty |\partial_t [\sqrt{t}\partial_{z_j} \mathbb{W}_t(z)]|dt
       \leq C 
    \int_0^\infty \frac{e^{-c\frac{|z|^2}{t}}}{t^{\frac{d}{2}+1}} dt
    \leq \frac{C}{|z|^d},
    \quad z\in\mathbb{R}^d, z\neq 0,
\end{equation}
and, for every $i=1,\dots,d$,
\begin{equation}\label{4.3}
    V_{\rho}\left(
    \{\sqrt{t} \partial^2_{z_j,z_i} \mathbb{W}_t(z)\}_{t>0}
    \right)
    \leq \frac{C}{|z|^{d+1}},
    \quad z\in\mathbb{R}^d,\,z\neq 0.
\end{equation}

Let $N\in\mathbb{N}$, $N\geq 2$. We define $E_{\rho,N}$ the space consisting of all those complex functions $g$ defined in $[1/N,N]$ such that $V_{\rho}(\{g(t)\}_{t\in[1/N,N]})<\infty$, where this variation is defined in the natural way. $(E_{\rho,N},V_\rho)$ can be seen as a Banach space by identifying the functions differing by a constant. 

Let $a\in (1/N,N)\setminus\{1\}$. We define $L_a g = g(a) - g(1)$, $g\in E_{\rho,N}$. It is clear that $L_a\in(E_{\rho})'$, the dual space of $E_\rho$.

Let $x\in \mathbb{R}^d$. We define $F_x: \mathbb{R}^d\rightarrow E_{\rho,N}$ such that, for $y\in\mathbb{R}^d$,
\begin{align*}
    F_x(y):\;\big[\frac{1}{N},N\big] \quad & \longrightarrow \quad \mathbb{C}
    \\
    t\quad \quad & \longmapsto [F_x(y)](t) = \sqrt{t} \partial_{x_j} \mathbb{W}_t(x-y)f(y).
\end{align*}
The function $F_x$ is continuous. Indeed, let $y_0\in \mathbb{R}^d$. By~\eqref{4.1} we have that
\begin{equation*}
    \begin{split}
        V_{\rho}(F_x(y) - F_x(y_0)) 
        & \leq 
        \int_{1/N}^N |\partial_t\left(
        [F_x(y)](t) - [F_x(y_0)](t)
        \right)| dt
        \\ & =
        \int_{1/N}^N \big|\partial_t\left(
        \sqrt{t}\partial_{x_j}\mathbb{W}_t(x-y_0) - 
        \sqrt{t}\partial_{x_j}\mathbb{W}_t(x-y)
        \right)\big| dt.
    \end{split}
\end{equation*}
Since the function $\partial_t[\sqrt{t} \partial_{x_j} \mathbb{W}_t(x-y)]$ is uniformly continuous in $(t,y)\in [1/N,N]\times B(y_0,1)$, we conclude that
\begin{equation*}
    \lim_{y\rightarrow y_0}V_\rho(F_x(y) - F_x(y_0))= 0.
\end{equation*}
Since $F_x$ is continuous, $F_x$ is $E_{\rho,N}$-strongly measurable.

By using~\eqref{4.2} we get
\begin{equation*}
    \int_{\mathbb{R}^d} V_\rho(F_x(y)) dy<\infty,
    \quad x\notin \supp(f).
\end{equation*}

We define
\begin{equation*}
    \mathbb{T}_j(f)(x) (t) = \int_{\mathbb{R}^d}
    \sqrt{t} \partial_{x_j}\mathbb{W}_t(x-y) f(y) dy,
    \quad x\in\mathbb{R}^d,
\end{equation*}
where the integral is understood in the $E_{\rho,N}$-Bochner sense.

Let $a\in [1/N,N]\setminus \{1\}$. Well known properties of the Bochner integral allow us to obtain
\begin{equation*}
    \begin{split}
        L_a(\mathbb{T} _j(f)(x)) 
        & = 
        \int_{\mathbb{R}^d} \sqrt{t}
        \partial_{x_j} \mathbb{W}_t (x-y)\rvert_{t=a}
        f(y) dy -
        \int_{\mathbb{R}^d} \sqrt{t}
        \partial_{x_j} \mathbb{W}_t (x-y)\rvert_{t=1}
        f(y) dy
        \\ & = 
        \left(\int_{\mathbb{R}^d} \sqrt{t}
        \partial_{x_j} \mathbb{W}_t (x-y)
        f(y) dy 
        \right)(a)
        -
        \left(\int_{\mathbb{R}^d} \sqrt{t}
        \partial_{x_j} \mathbb{W}_t (x-y)
        f(y) dy 
        \right)(1).
    \end{split}
\end{equation*}

Then, $\mathbb{T}_j(f)(x) = \tau_j(f)(x,\cdot)$ in $E_{\rho,N}$.

By using the vector valued Calder\'on-Zygmund theory we conclude that $\mathbb{T}_j$ is bounded from $L^1(\mathbb{R}^d,dx)$ into $L^{1,\infty}_{E_{\rho,N}}(\mathbb{R}^d,dx)$. Furthermore,
\begin{equation*}
    \sup_{N\in\mathbb{N},\, N>2} \|\mathbb{T}_j\|_{L^1(\mathbb R^d,dx)\longrightarrow L^{1,\infty}_{E_{\rho,N}}(\mathbb{R}^d,dx)}< \infty.
\end{equation*}
The monotone convergence theorem leads to see that the operator $V_\rho(\{\sqrt{t}\partial_{x_j}\mathbb{W}_t\}_{t>0})$ is bounded from
$L^1(\mathbb{R}^d,dx)$ into $L^{1,\infty}(\mathbb{R}^d,dx)$.

By proceeding as in Section 3 we can prove that the local variation operator $V_\rho(\{\mathbb{T}_{j,\loc}(f)(x,\cdot)\}_{t>0})$
is bounded from
$L^1(\mathbb{R}^d,\gamma_\infty)$ into $L^{1,\infty}(\mathbb{R}^d,\gamma_\infty)$. Here, the local operator $\mathbb{T}_{j,\loc}$ is defined by
\begin{equation*}
    \mathbb{T}_{j,\loc}(f)(x,t) = 
    \int_{\mathbb{R}^d}\sqrt{t}\partial_{x_j} \mathbb{W}_t(x-y)
    \varphi_{\mbox{\tiny A}}(x,y) f(y) dy, \quad x\in\mathbb{R}^d
    \text{ and } t>0.
\end{equation*}

The estimate \eqref{4.0} allows us to conclude that the local variation operator
$V_{\rho}(\{\mathbb{T}^\nu_{k,\alpha,\loc}(f)(x,t)\}_{t>0})$
is bounded from 
$L^1(\mathbb{R}^d,\gamma_\infty)$ into $L^{1,\infty}(\mathbb{R}^d,\gamma_\infty)$,
where $\alpha=(\alpha_1,\dots,\alpha_d)$ with $\alpha_i=0$ for $i=1,\dots,d$, $i\neq j$, $\alpha_j=1$ and
\begin{equation*}
    \mathbb{T}^\nu_{k,\alpha,\loc} (f)(x,t) = 
    \int_{\mathbb{R}^ d} t^{k + \widehat{\alpha}} 
    \mathbb{P}^\nu_{k,\alpha} (x-y,t) \varphi_{\mbox{\tiny A}}(x,y) f(y) dy,
    \quad x\in\mathbb{R}^d
    \text{ and } t>0,
\end{equation*}
being
\begin{equation*}
    \mathbb{P}^\nu_{k,\alpha} (z,t) = 
    \frac{1}{4^\nu \Gamma(\nu)} 
    \int_0^\infty \partial^k_t [t^{2\nu}e^{-\frac{t^2}{4u}}]
    \partial^\alpha_x \mathbb{W}_u(z) \frac{du}{u^{\nu+1}},
    \quad z\in \mathbb{R}^d
    \text{ and } t>0.
\end{equation*}

We consider the operator $S^\nu_{k,j}$ and the kernel $K_{k,j}^\nu$ as in Section \ref{S3.1.1}. We also define
\begin{equation*}
    V_{\rho,\loc}(\{t^{k+1}\partial^k_t \partial_{x_j} P^\nu_t(f)(x)\}_{t>0})
    = V_{\rho}\Big(
    \Big\{\int_{\mathbb{R}^d} K_{k,j}^\nu(x,y,t)\varphi_{\mbox{\tiny A}}(x,y)f(y)
    d\gamma_\infty (y)\Big\}_{t>0}
    \Big), \quad x\in\mathbb{R}^d.
\end{equation*}

We have that
\begin{equation*}
    \begin{split}
        \Big|
        V_{\rho,\loc} & (\{t^{k+1}\partial^k_t \partial_{x_j} P^\nu_t(f)(x)\}_{t>0})
        - V_{\rho}(\{S_{k,j,\loc}^\nu (f)(x,t)\}_{t>0})
        \Big|
        \\ & 
        \leq 
        V_{\rho}\left(
        \left\{\int_{\mathbb{R}^d} K_{k,j}^\nu(x,y,t)\varphi_{\mbox{\tiny A}}(x,y)f(y)
    dy
        -
        S_{k,j,\loc}^\nu (f)(x,t)\right\}_{t>0}
        \right)
        \\ & \leq 
        C\int_{\mathbb{R}^d} \varphi_{\mbox{\tiny A}}(x,y)|f(y)|e^{-R(y)}\int_0^\infty
        \int_0^\infty \Big|\partial_t\big([s^{k+1}\partial^k_s\mathfrak{g}_\nu (s)]_{|s=\frac{t}{2\sqrt{u}}}\big)\Big|dt  |\partial_{x_j} h_u(x,y)
        - \mathbb{S}_u^j(x,y)|\frac{du}{\sqrt{u}}dy\\
        &\leq C\int_{\mathbb{R}^d} \varphi_{\mbox{\tiny A}}(x,y)|f(y)|e^{-R(y)}\int_0^\infty |\partial_{x_j} h_u(x,y)
        - \mathbb{S}_u^j(x,y)|\frac{du}{\sqrt{u}}dy,\quad x\in\mathbb{R}^d.
    \end{split}
\end{equation*}
In the last inequality we have used that
$$
\Big|\partial_t\big([s^{k+1}\partial^k_s\mathfrak{g}_\nu (s)]_{|s=\frac{t}{2\sqrt{u}}}\big)\Big|\leq \frac{C}{\sqrt{u}}\big(|s^k\partial _s^k\mathfrak{g}_\nu (s)|+|s^{k+1}\partial_s^{k+1}\mathfrak{g}_\nu (s)|\big)_{|s=\frac{t}{2\sqrt{u}}}\leq \frac{C}{\sqrt{u}}\mathfrak{g}_\nu \big(\frac{s}{2}\big)_{|s=\frac{t}{2\sqrt{u}}},\quad t,u>0.
$$

The arguments in Section~\ref{S3.1.1} allow us to establish that $V_{\rho,\loc}  (\{t^{k+1}\partial^k_t \partial_{x_j} P^\nu_t(f)(x)\}_{t>0})$ is a bounded operator from $L^1(\mathbb{R}^d,\gamma_\infty)$ into $L^{1,\infty}(\mathbb{R}^d,\gamma_\infty)$.

We consider the global variation operator defined by
\begin{equation*}
    V_{\rho,\glob}(\{t^{k+1}\partial^k_t \partial_{x_j} P^\nu_t(f)(x)\}_{t>0})
    = V_{\rho}\Big(
    \Big\{\int_{\mathbb{R}^d} P^{\nu}_{k,\alpha}(x,y,t)(1-\varphi_{\mbox{\tiny A}}(x,y))f(y)
    d\gamma_\infty(y)\Big\}_{t>0}
    \Big), \quad x\in\mathbb{R}^d.
\end{equation*}

Since
\begin{equation*}
    \int_0^\infty \Big|\partial_t\big([s^{k+1} \partial^k_s\mathfrak{g}_\nu (s)]_{|s=\frac{t}{2\sqrt{u}}}\big)\Big|dt \leq C,
    \quad u>0,
\end{equation*}
we get
\begin{equation*}
    V_{\rho,\glob}(\{t^{k+1}\partial^k_t \partial_{x_j} P^\nu_t(f)(x)\}_{t>0})
    \leq C \int_{\mathbb{R}^d} (1-\varphi_{\mbox{\tiny A}}(x,y))  |f(y)|
   e^{-R(y)}
    \int_0^\infty |\partial_{x_j} h_u(x,y)| \frac{du}{\sqrt{u}}dy, \quad x\in\mathbb{R}^d.
\end{equation*}

As in Section~\ref{S3.1.3} we conclude that the global variation operator 
$V_{\rho,\glob}  (\{t^{k+1}\partial^k_t \partial_{x_j} P^\nu_t(f)(x)\}_{t>0})$ is bounded from $L^1(\mathbb{R}^d,\gamma_\infty)$ into $L^{1,\infty}(\mathbb{R}^d,\gamma_\infty)$. Thus, we establish Theorem \ref{Th1.3} when $\widehat{\alpha}=1$.

We consider the case $\widehat{\alpha}= 2$. We can proceed as above by taking into account the properties that we are going to mark.

We denote by $\phi = \mathcal{F}^{-1}\{|y|^2 e^{-|y|^2}\}$. We have that
\begin{equation*}
    \phi(x) = h_{\frac{d-2}{2}} (v^2 e^{-v^2})(|x|), \quad
    x\in \mathbb{R}^d.
\end{equation*}

According to~\cite[p. 30, (13)]{EMOT} we get
\begin{equation*}
    \psi(u) = h_{\frac{d-2}{2}}(v^{2}e^{-v^2})(u) = 2e^{-\frac{u^2}{2}} L_1^{\frac{d-2}{2}}\big(\frac{u^2}{2}\big), \quad u>0,
\end{equation*}
where $L^\sigma_n$ denotes the $n$-th $\sigma$-Laguerre polynomial. Note that 
    $\lim_{u\rightarrow \infty}\psi(u) = 0$ and 
    $\int_0^{\infty} |\psi'(u)| u^d du<\infty$.
    
    On the other hand, we have that
    \begin{equation*}
        \int_0^\infty |\partial_t [s^{k+2} \partial^k_s \mathfrak{g}_\nu (s)]_{|s=\frac{t}{2\sqrt{u}}}|
    dt \leq C, \quad u>0.
    \end{equation*}
    Our above arguments can be used now to prove Theorem~\ref{Th1.3} when $\widehat{\alpha} = 2$.
    
    To finish we study the case $\alpha = 0$. As in~\eqref{4.0} we can obtain
    \begin{equation*}
        V_{\rho}(\{t^k \partial^k_t P^\nu_t (f)(x)\}_{t>0}) 
        \leq C 
        V_{\rho}(\{\mathcal{H}_t(f)(x)\}_{t>0}),\quad x\in\mathbb{R}^d.
    \end{equation*}
The semigroup $\{\mathcal{H}_t\}_{t>0}$ is contractive and analytic in 
$L^p(\mathbb{R}^d,\gamma_\infty)$, for $1<p<\infty$ (see~\cite[Theorem 2]{CFMP}).
Furthermore, $\{\mathcal{H}_t\}_{t>0}$ is contractive in 
$L^1(\mathbb{R}^d,\gamma_\infty)$ and $L^\infty(\mathbb{R}^d,\gamma_\infty)$.
Then,  $\{\mathcal{H}_t\}_{t>0}$  is contractively regular and, according to~\cite[Corollary 4.5]{LeMXu}, the variation operator  $V_\rho(\{\mathcal{H}_t\}_{t>0})$ is bounded from  $L^p(\mathbb{R}^d,\gamma_\infty)$ into itself, for $1<p<\infty$. We conclude that the operator
$V_{\rho}( \{t^k \partial^k_t P^\nu_t \}_{t>0})$ is bounded from $L^p(\mathbb{R}^d,\gamma_\infty)$ into itself, for $1<p<\infty$.

As far as we know, it has not been proved that the variation operator 
$V_\rho(\{\mathcal{H}_t\}_{t>0})$  is bounded from 
$L^1(\mathbb{R}^d,\gamma_\infty)$ into $L^{1,\infty}(\mathbb{R}^d,\gamma_\infty)$. Then, we need to proceed in a different way to prove that the operator 
$V_{\rho}( \{t^k \partial^k_t P^\nu_t \}_{t>0}) $ 
is bounded from 
$L^1(\mathbb{R}^d,\gamma_\infty)$ into $L^{1,\infty}(\mathbb{R}^d,\gamma_\infty)$.

Suppose that $k\in\mathbb{N}$ and $k\geq 1$ and consider the elements given in Section \ref{S3.2}. Let $f\in C^\infty_c(\mathbb{R}^d)$. We can write
$$
V_\rho (\{t^k\partial_t^kP_t^\nu (f)(x)\}_{t>0})\leq C\sum_{i,j=1}^dV_\rho (\{\mathscr{K}^{\nu}_{k,i,j} (f) (x,t)\}_{t>0})+V_\rho (\{ \mathscr{R}^{\nu}_{k,i,j} (f) (x,t) \}_{t>0}),\quad x\in \mathbb{R}^d.
$$
To proceed as in Section \ref{S3.2} we need the following estimation 
\begin{equation*}
    \int_0^\infty 
    \big|\partial_t ([s^k\partial_s^{k-1}\mathfrak{h}_\nu (s)]_{|s=\frac{t}{2\sqrt{u}}})\big|dt=\int_0^\infty |\partial_s(s^k\partial _s^{k-1}\mathfrak{h}_\nu (s))|ds\leq C\int_0^\infty \mathfrak{g}_\nu \big(\frac{s}{2}\big)\frac{ds}{s} 
    \leq C
    , \quad  u>0.
\end{equation*}
Here, we have taken into account that $|s^\ell\partial _s^\ell \mathfrak{h}_\nu (s)|\leq Cs^{-1}\mathfrak{g}_\nu (\frac{s}{2})$, $s>0$, $\ell \in \mathbb{N}$. 

On the other hand, using the above arguments for $\widehat{\alpha}=2$ and proceeding as in the case of $\widehat{\alpha}=1$ we get that $V_{\rho,{\rm loc}}(\{\mathscr{T}_{k,i,j}^\nu(f)(x)\}_{t>0})$
 is bounded from $L^1(\mathbb{R}^d,\gamma_\infty)$ into $L^{1,\infty}(\mathbb{R}^d,\gamma_\infty )$.

We can apply the reasoning in Section \ref{S3.2} to establish Theorem \ref{Th1.3} for $\alpha =0$ and $k\geq 1$.

Suppose now that $k=0$. We are going to see that the operator $V_\rho(\{P^\nu_t\}_{t>0})$ is bounded from
$L^1(\mathbb{R}^d,dx)$ into $L^{1,\infty}(\mathbb{R}^d,dx)$. 
We first study the global variation operator 
$$
V_{\rho,\glob}(\{P^\nu_t\}_{t>0})(f)(x)=V_\rho(\{P_{t,\glob}^\nu(f)(x)\}_{t>0}),\quad x\in \mathbb{R}^d.
$$
We write
\begin{equation*}
    P^\nu_t(x,y) = \frac{1}{\Gamma(\nu)}
    \int_0^\infty e^{-u} u^{\nu-1} h_{\frac{t^2}{4u}} (x,y) du,\quad x,y\in \mathbb{R}^d\mbox{ and }t>0.
\end{equation*}
It follows that
\begin{align*}
  V_{\rho,\glob}(\{P^\nu_t\}_{t>0})(f)(x)   & \leq 
        C \int_{\mathbb{R}^d}(1-\varphi_{\mbox{\tiny A}}(x,y))|f(y)| 
        \int_0^\infty e^{-u} u^{\nu-1} \int_0^\infty 
        |\partial_t  h_{\frac{t^2}{4u}}(x,y)| dt du
        d\gamma_\infty (y)\\
& \leq C \int_{\mathbb{R}^d}(1-\varphi_{\mbox{\tiny A}}(x,y))|f(y)|
    \int_0^\infty e^{-u} u^{\nu-1} \int_0^\infty 
    |\partial_s  h_s(x,y)| ds du   d\gamma_\infty (y)
    \\ 
& \leq C
    \int_{\mathbb{R}^d}(1-\varphi_{\mbox{\tiny A}}(x,y))|f(y)| 
    \int_0^\infty 
    |\partial_s  h_{s}(x,y)| ds   d\gamma_\infty (y),
    \quad x\in\mathbb{R}^d.
\end{align*}
By taking into account \eqref{duxd2x} and proceeding as in Section \ref{S3.2.2} we obtain that $V_{\rho,\glob}(\{P^\nu_t\}_{t>0})$ 
 is bounded from $L^1(\mathbb{R}^d,\gamma_\infty)$ into 
$L^{1,\infty}(\mathbb{R}^d,\gamma_\infty)$, provided that $A$ is large enough. 

To deal with the the local operator 
$$
V_{\rho,{\rm loc}}(\{P_t^\nu \}_{t>0})(f)(x)=V_\rho (\{P_{t,{\rm loc}}^\nu (f)(x)\}_{t>0}),\quad x\in \mathbb{R}^d,
$$
we consider the operator defined by
\begin{equation*}
    S^\nu(f)(x,t) = 
    \int_{\mathbb{R}^d} S^{\nu}(x,y,t)f(y) d\gamma_\infty (y), \quad
    x\in\mathbb{R}^d \text{ and } t>0,
\end{equation*}
where
\begin{equation*}
    S^\nu(x,y,t) = \frac{1}{\Gamma(\nu)}\int_0^\infty \mathfrak{g}_\nu \Big(\frac{t}{2\sqrt{u}}\Big) \widetilde{W}_u(x,y)\frac{du}{u},\quad x,y \in\mathbb{R}^d \text{ and } t>0,
\end{equation*}
where $\widetilde{W}_u(x,y)$, $x,y\in \mathbb{R}^d$, $u>0$ is the kernel given by \eqref{Wtilde}.

Our next arguments are inspired by some ideas developed in~\cite[Section 4]{HMMT} where the symmetric case is studied.
We consider the local operator $S_{\loc}^\nu$ in the usual way and write $D^\nu _{\rm loc}(f)(x,t)=P_{t,{\rm loc}}^\nu (f)(x)-S_{\rm loc}^\nu (f)(x,t)$, $x\in \mathbb{R}^d$, $t>0$. We have that
\begin{align*}
    \mathcal{D}^\nu _{\rm loc}(f)(x,t)&=\frac{1}{\Gamma (\nu )}\int_{\mathbb{R}^d}\varphi_{\mbox{\tiny A}} (x,y)f(y)\int_0^\infty \mathfrak{g}_\nu \Big(\frac{t}{2\sqrt{u}}\Big)\big[h_u(x,y)-\widetilde{W}_u(x,y)\big]\frac{du}{u}d\gamma_\infty (y)\\
    &=\frac{1}{\Gamma (\nu )}\int_{\mathbb{R}^d}\varphi_{\mbox{\tiny A}} (x,y)f(y)\int_0^\infty \mathfrak{g}_\nu \Big(\frac{t}{2\sqrt{u}}\Big)[h_u(x,y)-\mathcal{X}_{(1,\infty )}(u)-\widetilde{W}_u(x,y)]\frac{du}{u}d\gamma_\infty (y)\\
    &\quad +\frac{1}{\Gamma (\nu )}\int_{\mathbb{R}^d}\varphi_{\mbox{\tiny A}} (x,y)f(y)\int_1^\infty \mathfrak{g}_\nu \Big(\frac{t}{2\sqrt{u}}\Big)\frac{du}{u}d\gamma _\infty (y)\\
    &=\mathcal{D}_{\rm loc}^{\nu ,1}(f)(x,t)+\mathcal{D}_{\rm loc}^{\nu ,2}(f)(x,t),\quad x\in \mathbb{R}^d\mbox{ and }t>0.
\end{align*}

Since $\int_0^\infty \mathfrak{g}_\nu \big(\frac{t}{2\sqrt{u}}\big)\frac{du}{u}= \Gamma(\nu)$, $t>0$, we have that 
 $\frac{d}{dt}\left[\int_0^\infty \mathfrak{g}_\nu \big(\frac{t}{2\sqrt{u}}\big)\frac{du}{u}\right]=0$, and then
\begin{equation*}
    \frac{d}{dt}\int_1^{\infty}\mathfrak{g}_\nu \big(\frac{t}{2\sqrt{u}}\big)\frac{du}{u}
     = -\frac{d}{dt}\int_0^1 \mathfrak{g}_\nu \big(\frac{t}{2\sqrt{u}}\big)\frac{du}{u}, \quad t>0.
\end{equation*}
Thus, by using also that $s\mathfrak{g}_\nu '(s)\leq C\mathfrak{g}_\nu (\frac{s}{2})$, $s>0$, we can write,
\begin{align*}
    V_\rho (\{\mathcal{D}_{\rm loc}^{\nu ,2}(f)(x,t)\}_{t>0})&\leq C\int_{\mathbb{R}^d}\varphi_{\mbox{\tiny A}} (x,y)|f(y)|\int_0^\infty \left|\frac{d}{dt}\int_1^\infty \mathfrak{g}_\nu \Big(\frac{t}{2\sqrt{u}}\Big)\frac{du}{u}\right|dtd\gamma_\infty (y)\\
    &\hspace{-3cm}\leq C\int_{\mathbb{R}^d}\varphi_{\mbox{\tiny A}} (x,y)|f(y)|\left(\int_0^1 \int_1^\infty \left|\partial _t\Big[\mathfrak{g}_\nu \Big(\frac{t}{2\sqrt{u}}\Big)\Big]\right|\frac{du}{u}dt \right. +\left.\int_1^\infty \int_0^1\left|\partial _t\Big[\mathfrak{g}_\nu \Big(\frac{t}{2\sqrt{u}}\Big)\Big]\right|\frac{du}{u}dt\right)d\gamma_\infty(y)\\
    &\hspace{-3cm}\leq C\int_{\mathbb{R}^d}\varphi_{\mbox{\tiny A}} (x,y)|f(y)|\left(\int_0^1 \int_1^\infty +\int_1^\infty \int_0^1\right)t^{2\nu -1}e^{-c\frac{t^2}{u}}\frac{du}{u^{\nu +1}}dtd\gamma_\infty (y)\\
    &\hspace{-3cm}\leq C\int_{\mathbb{R}^d}\varphi_{\mbox{\tiny A}} (x,y)|f(y)|\left(\int_0^1 \int_1^\infty\frac{t^{2\nu -1}}{u^{\nu +1}}dudt +\int_1^\infty \int_0^1\frac{u^{\delta -1}}{t^{\delta +1}}dudt\right)d\gamma_\infty (y)\\
    &\hspace{-3cm}\leq C\int_{\mathbb{R}^d}\varphi_{\mbox{\tiny A}} (x,y)|f(y)|d\gamma_\infty (y)\leq C\int_{\mathbb{R}^d}\varphi_{\mbox{\tiny A}} (x,y)|f(y)|\frac{1+|x|}{|x-y|^{d-1}}dy=C\mathscr{S}_{2A}(|f|)(x),\quad x\in \mathbb{R}^d.
\end{align*}
Here $\delta$ can be any positive real number. In the last inequality we have taken into account that $|x-y|^{d-1}\leq 1\leq 1+|x|$, when $(x,y)\in L_{2A}$.

On the other hand, by considering that
$\int_0^\infty |\partial _t[\mathfrak{g}_\nu \big(\frac{t}{2\sqrt{u}}\big)]|dt\leq C$, $u>0$, it follows that
\begin{align*}
    V_\rho (\{\mathcal{D}_{\rm loc}^{\nu ,1}(f)(x,t)\}_{t>0})&\\
    &\hspace{-3cm}\leq C\int_{\mathbb{R}^d}\varphi_{\mbox{\tiny A}} (x,y)|f(y)| \int_0^\infty \int_0^\infty \Big|\partial _t\Big[\mathfrak{g}_\nu \big(\frac{t}{2\sqrt{u}}\big)\Big]\Big||h_u(x,y)-\mathcal{X}_{(1,\infty )}(u)-\widetilde{W}_u(x,y)|dt\frac{du}{u}d\gamma_\infty (y)\\
    &\hspace{-3cm}\leq C\int_{\mathbb{R}^d}\varphi_{\mbox{\tiny A}} (x,y)|f(y)|\int_0^\infty|h_u(x,y)-\mathcal{X}_{(1,\infty )}(u)-\widetilde{W}_u(x,y)|\frac{du}{u}d\gamma_\infty (y)
    ,\quad x\in \mathbb{R}^d.
\end{align*}

We now consider
$$
R(x,y)= e^{-R(y)}\int_0^\infty
    |h_u(x,y) - \mathcal{X}_{(1,\infty )}(u) -\widetilde{W}_u(x,y)| 
     \frac{du}{u}, \quad x,y\in\mathbb{R}^d.
$$
Our objective is to establish that
\begin{equation}\label{R}
|R(x,y)|\leq C\frac{1+|x|}{|x-y|^{d-1}},\quad (x,y)\in L_{2A}.
\end{equation}
We decompose $R(x,y)=R_0(x,y)+R_\infty (x,y)$, $x,y\in \mathbb{R}^d$, where
$$
    R_0(x,y) = e^{-R(y)}
    \int_0^1
    |h_u(x,y)-\mathcal{X}_{(1,\infty )}(u) -\widetilde{W}_u(x,y)| \frac{du}{u},\quad x,y\in \mathbb{R}^d.
$$

According to \eqref{3.2'}, \eqref{3.0}, \eqref{3.1'}, \eqref{3.6}, \eqref{3.7.1} and \cite[Lemma 2.2]{CCS3} we get, for $(x,y)\in L_{2A}$ and $0<u<\mathfrak{m}(x)$, 
\begin{align*}
e^{-R(y)}|h_u(x,y)-\widetilde{W}_u(x,y)|
    &\leq C
    e^{R(x)-R(y)}\left(\Big|\frac{1}{({\rm det}\,Q_u)^{1/2}}-\frac{1}{(u^d{\rm det}\,Q)^{1/2}}\Big|e^{ -\frac{1}{2}\langle  (Q_u^{-1}-Q^{-1}_{\infty})(y-D_u x),y- D_u x\rangle}\right.\\
    &\left.\quad +\frac{1}{u^{\frac{d}{2}}}\Big|e^{-\frac{1}{2}\langle  (Q_u^{-1}-Q^{-1}_{\infty})(y-D_u x),y- D_u x\rangle}-e^{-\frac{1}{2u}\langle  Q^{-1}(y- x),y- x\rangle}\Big|\right) \\
    &\leq C\frac{e^{R(x)-R(y)}}{u^{\frac{d}{2}}}e^{-c\frac{|y-x|^2}{u}}\Big(u+(1+|x|)\sqrt{u}\Big)\leq C(1+|x|)\frac{e^{-c\frac{|y-x|^2}{u}}}{u^{\frac{d-1}{2}}}.
\end{align*}
Then,
\begin{equation}\label{5.34}
        e^{-R(y)} \int_0^{\mathfrak{m}(x)} |h_u(x,y)-    \widetilde{W}_u(x,y)| \frac{du}{u}
        \leq C(1+|x|)\int_0^{\mathfrak{m}(x)} 
        \frac{e^{-c\frac{|y-x|^2}{u}}}{u^{\frac{d+1}{2}}} du
        \leq C \frac{1+|x|}{|x-y|^{d-1}}, \quad (x,y)\in L_{2A}.
\end{equation}

By using~\cite[(2.10)]{CCS3} and~\eqref{3.2'} we get
\begin{align}\label{5.37}
    e^{-R(y)} \int_{\mathfrak{m}(x)}^1  | h_u(x,y) - \widetilde{W}_u(x,y)| \frac{du}{u} 
    &     \leq Ce^{R(x)-R(y)} \int_{\mathfrak{m}(x)}^\infty \frac{du}{u^{\frac{d}{2}+1}}\leq  \frac{C}{\mathfrak{m}(x)^{\frac{d}{2}}}
  \nonumber\\
   & \leq C(1+|x|)^d\leq C \frac{1+|x|}{|x-y|^{d-1}},
    \quad (x,y)\in L_{2A}.
\end{align}

From~\eqref{5.34} and~\eqref{5.37} we get that
\begin{equation}\label{5.38}
    0\leq R_0(x,y) \leq C \frac{1+|x|}{|x-y|^{d-1}},
    \quad (x,y)\in L_{2A}.
\end{equation}

We now write
\begin{align*}
R_\infty(x,y)&\leq Ce^{-R(y)}\left(\int_1^\infty |h_u(x,y)-1|\frac{du}{u}+\int_1^\infty \widetilde{W}_u(x,y)\frac{du}{u}\right),\quad x,y\in \mathbb{R}^d.
\end{align*}
We observe that
\begin{equation}\label{S}
e^{-R(y)}\int_1^\infty \widetilde{W}_u(x,y)\frac{du}{u}\leq Ce^{R(x)-R(y)} \int_1^{\infty}  \frac{du}{u^{\frac{d}{2}+1}} \leq C \leq C \frac{1+|x|}{|x-y|^{d-1}},
\quad (x,y)\in L_{2A}.
\end{equation}
On the other hand, we have that $\lim_{u\rightarrow \infty} h_u(x,y) = 1$, 
$x$, $y\in\mathbb{R}^d$. Indeed, since $D_u^{-1}=D_{-u}$, $u>0$, as in the proof of~\cite[Proposition 3.3]{CCS2}, we get
\begin{align*}
    \langle 
    (Q_u^{-1}-Q_{\infty}^{-1})(y-D_u x),
    y-D_u x\rangle &
    = \langle 
    (Q_u^{-1}-Q_{\infty}^{-1}) D_u(D_{-u}y- x),
    D_u(D_{-u}y- x)\rangle
    \nonumber\\ 
    &\hspace{-1cm} = 
    \langle Q^{-1}_\infty (D_{-u}y- x), D_{-u}y -x 
    \rangle
    + |Q^{-1/2}_u e^{uB}(D_{-u}y-x)|^2, 
    \quad x,y\in\mathbb{R}^d \text{ and } u>0.
\end{align*}
Then, by taking into account that $|D_{-u}y-x|\simeq |D_{-u}y|$, $x,y\in \mathbb{R}^d$, $u>0$, (see the proof of Lemma 5.1 in \cite{CCS1}) and considering \cite[Lemmas 2.1 and 2.2]{CCS3} we obtain
\begin{equation*}
    \lim_{t\rightarrow\infty} \langle(Q^{-1}_t - Q^{-1}_\infty)(y-D_t x),
    y-D_tx\rangle = \langle Q^{-1}_\infty x, x\rangle, 
    \; x,y\in\mathbb{R}^d. 
\end{equation*}
Then, we conclude that $\lim_{u\rightarrow \infty} h_u(x,y) = 1$, $x, y\in\mathbb{R}^d$. 

Thus, by considering also \eqref{(2.8)}, we can write 
\begin{align*}
|h_u(x,y)-1|&\leq\int_u^\infty |\partial _sh_s(x,y)|ds\leq Ce^{R(x)}(1+|y|)^2\int_u^\infty e^{-cs}ds\\
&=Ce^{R(x)}(1+|y|)^2e^{-cu},\quad x,y\in \mathbb{R}^d\mbox{ and }u>0.
\end{align*}
We use again \eqref{3.2'} and that $|y|\leq 1+|x|$, when $(x,y)\in L_{2A}$, to get
\begin{equation*}
e^{-R(y)} \int_1^\infty |h_u(x,y) -1|
\frac{du}{u} \leq Ce^{R(x)-R(y)}(1+|x|)^2\leq C(1+|x|)^d\leq C \frac{1+|x|}{|x-y|^{d-1}},
        \quad (x,y)\in L_{2A},
\end{equation*}
which, jointly with \eqref{S}, leads to
\begin{equation}\label{3.41}
    0\leq R_\infty(x,y)  \leq C \frac{1+|x|}{|x-y|^{d-1}},
        \quad (x,y)\in L_{2A}.
\end{equation}

From \eqref{5.38} and \eqref{3.41} we obtain \eqref{R}.

All the above estimations allow us to conclude that
$$
V_\rho (t\mapsto P_{t,{\rm loc}}^\nu (f)(x)-S_{\rm loc}^\nu (f)(x,t))\leq C\mathscr{S}_{2A}(f)(x),\quad x\in \mathbb{R}^d,
$$
and therefore, that $V_\rho (\{\mathcal{D}_{t,{\rm loc}}^\nu\}_{t>0})$ defined by 
$$
V_\rho (\{\mathcal{D}_{t,{\rm loc}}^\nu\}_{t>0})(f)(x)=V_\rho (t\mapsto P_{t,{\rm loc}}^\nu (f)(x)-S_{\rm loc}^\nu (f)(x,t)),\quad x\in \mathbb R^d,
$$
is bounded from $L^1(\mathbb{R}^d,\gamma_\infty)$ into itself. 

To conclude that the local variation operator $V_{\rho,{\rm loc}}(\{P_t^\nu \}_{t>0})$ is bounded from $L^1(\mathbb{R}^d,\gamma_\infty)$ into $L^{1,\infty}(\mathbb{R}^d,\gamma_\infty)$ we can proceed now as in Section \ref{S3.1.1}. We need the estimation $\int_0^\infty |\partial _t[\mathfrak{g}_\nu \big(\frac{t}{2\sqrt{u}}\big)]|dt\leq C$, $u>0$, and the appropriate $L^p$-boundedness properties for the operator $V_\rho (\{\mathbb{P}_t^\nu \}_{t>0})$.

According to~\cite[Corollary 2.7]{CJRW1} the variation operator 
$V_{\rho}(\{ \mathbb{W}_t\}_{t>0})$ is bounded from $L^p(\mathbb{R}^d,dx)$ into itself, for every $1<p<\infty$, and from $L^1(\mathbb{R}^d,dx)$ into $L^{1,\infty}(\mathbb{R}^d,dx)$. Then, arguing as in \eqref{4.0}, $V_\rho(\{\mathbb{P}^\nu_t\}_{t>0})$ is bounded  from $L^p(\mathbb{R}^d,dx)$ into itself, for every $1<p<\infty$, and from $L^1(\mathbb{R}^d,dx)$ into $L^{1,\infty}(\mathbb{R}^d,dx)$.

We also have that
$$
V_{\rho}(\{\mathbb{P}^\nu_t(z)\}_{t>0})\leq C \int_0^\infty \int_0^{\infty} 
\Big|\partial_t\Big[\mathfrak{g}_\nu \Big(\frac{t}{2\sqrt{u}}\Big)\Big]\Big|
\mathbb{W}_u(z) \frac{du}{u}dt \leq C\int_0^\infty\frac{e^{-\frac{|z|^2}{2u}}}{u^{\frac{d}{2}+1}} du\leq \frac{C}{|z|^d}, \quad z\in\mathbb{R}^d\setminus\{0\},
$$
and, in a similar way, we get
$$
        V_{\rho}(\{\partial_{z_k}\mathbb{P}^\nu_t(z)\}_{t>0})
        \leq \frac{C}{|z|^{d+1}}, \quad z\in\mathbb{R}^d\setminus\{0\} \text{ and }
        k=1,\dots,d.
$$

We have all the ingredients to argue as in Section \ref{S3.1.1} and to conclude that the operator $V_{\rho,{\rm loc}}(\{P^\nu_t\}_{t>0})$ is bounded  from $L^1(\mathbb{R}^d,\gamma_\infty)$ into $L^{1,\infty}(\mathbb{R}^d,\gamma_\infty)$. Theorem \ref{Th1.3} for $\alpha =0$ is then established.


\end{document}